\documentclass[dvips,12pt]{article}
\usepackage{amsmath,amssymb,rotating,amsfonts,pstricks,pstricks-add
}
\def\bt{\begin{tabular}}
\def\te{\end{tabular}}

\def\BM{\begin{pmatrix}}
\def\EM{\end{pmatrix}}

\def\cit{\text{\it I\hskip -6ptC\/}}
\def\ptcit{\hbox{${\scriptstyle I\hskip -4ptC\/}$}}

\def\square{\hfill\hbox{\vrule height .9ex width .8ex depth -.1ex}}
\def\rit{\text{\it I\hskip -2pt  R}}

\def\zit{\text{\it Z\hskip -4pt  Z}}
\def\qit{\text{\it I\hskip -5.5pt  Q}}

\def\nit{\text{\it I\hskip -2pt  N}}

\def\Ad{{\text A}}
\def\Bd{{\text B}}

\def\Ed{{\text E}}

\def\Fs{{\cal F}}

\def\Zs{{\cal Z}}

\def\be{\begin{equation}}
\def\ee{\end{equation}}
\def\beqn{\begin{eqnarray}}
\def\eeqn{\end{eqnarray}}
\def\nobeqn{\begin{eqnarray*}}
\def\noeeqn{\end{eqnarray*}}
\def\ba{\left(\begin{array}}
\def\ea{\end{array} \right) }

\def\bpr{\paragraph{Proof.}}

\def\epr{\square\vskip 6pt}
\def\eop{\hbox{\vrule height .9ex width .8ex depth -.1ex}}

\def\o{\overline}
\def\and{\; \mbox{and} \;}

\newcommand{\half}{\frac{1}{2}}

\def\Ker{\mathop{\rm Ker}\nolimits}

\def\mod{\mathop{\rm mod}\nolimits}

\def\ch{\mathop{\rm ch}\nolimits}

\def\Be{\begin{enumerate}}
\def\Ee{\end{enumerate}}

\def\Bena{\begin{enumerate}
\def\labelenumi{\theenumi)}
\def\theenumi{\arabic{enumi}}
\def\labelenumii{\theenumii)}
\def\theenumii{\alph{enumii}}}

\def\Bean{\begin{enumerate}
\def\labelenumii{\theenumii)}
\def\theenumii{\arabic{enumii}}
\def\labelenumi{\theenumi)}
\def\theenumi{\alph{enumi}}}

\def\Bero{\begin{enumerate}
\def\labelenumii{\theenumii)}
\def\theenumii{\arabic{enumii}}
\def\labelenumi{(\theenumi)}
\def\theenumi{\roman{enumi}}}

\def\BeRo{\begin{enumerate}
\def\labelenumii{\theenumii)}
\def\theenumii{\arabic{enumii}}
\def\labelenumi{(\theenumi)}
\def\theenumi{\Roman{enumi}}}

\def\Bi{\vskip 11pt\begin{itemize}\itemsep=18pt}
\def\Ei{\end{itemize}}

\def\Bd{\begin{description}}
\def\Ed{\end{description}}

\def\R{\right}
\def\L{\left}

\def\bigoplus{\mathop{\oplus}\limits}

\def\prod{\mathop{\Pi}\limits}
\def\sum{\mathop{\Sigma}\limits}

\def\Bsm{\begin{smallmatrix}}
\def\Esm{\end{smallmatrix}}

\def\resp#1{(resp. #1)}
\def\rresp#1{\qquad \mbox{(resp.} \quad #1\ )}

\def\bbf{\boldmath\bf}
\def\o{\overline}
\def\wt{\widetilde}

\def\Bi{\begin{itemize}}
\def\Ei{\end{itemize}}

\newcommand{\ZZ}{\mathbb{Z}\,}

\newcommand{\HH}{{\rm{H}}\,}
\def\tr{\operatorname{tr}}

\def\Gal{\operatorname{Gal}}

\def\mod{\operatorname{mod}}

\def\cor{\operatorname{cor}}
\def\Ad{\operatorname{Ad}}
\def\AdF{\operatorname{AdF}}

\def\Irr{\operatorname{Irr}}

\def\Rep{\operatorname{Rep}}

\def\REPSP{\operatorname{REPSP}}
\def\FREPSP{\operatorname{FREPSP}}
\def\fREPSP{\operatorname{(F)REPSP}}
\def\UfREPSP{\operatorname{U(F)REPSP}}
\def\TAN{\operatorname{TAN}}

\def\RED{\operatorname{RED}}

\def\Int{\operatorname{Int}}
\def\Aut{\operatorname{Aut}}
\def\Out{\operatorname{Out}}

\def\BGL{\operatorname{BGL}}
\def\BFGL{\operatorname{BFGL}}
\def\BFG{\operatorname{BFG}}
\def\BG{\operatorname{BG}}
\def\LGC{\operatorname{LGC}}

\def\EG{\operatorname{EG}}
\def\EGL{\operatorname{EGL}}
\def\GL{\operatorname{GL}}
\def\FGL{\operatorname{FGL}}
\def\EFGL{\operatorname{EFGL}}
\def\GD{\operatorname{GD}}

\def\Tr{\operatorname{tr}}

\def\bM{\begin{matrix}}
\def\eM{\end{matrix}}

\def\lr{left (resp. right) }

\def\bpr{\noindent{\bf{Proof\/}}:\;\;}
\def\To{\longrightarrow }

\def\RL{_{R\times L}}

\begin{document}

\setcounter{page}{0}
{\pagestyle{empty}
\null\vfill
\begin{center}
{\LARGE Higher algebraic $K$-theories related to the global program of  Langlands}
\vfill
{\sc C. Pierre\/}
\vskip 11pt


\vfill
\vfill

\begin{abstract}
The paper revisits concretely the algebraic $K$-theory in the light of the global program of Langlands by taking into account the new algebraic interpretation of homotopy viewed as deformation(s) of Galois representations given by compactified algebraic groups.

More concretely,  we introduce higher algebraic bilinear $K$-theories referring to homotopy and cohomotopy and related to the reducible bilinear global program of Langlands as well as mixed higher bilinear $KK$-theories related to dynamical geometric bilinear global program of Langlands.
\end{abstract}
\vfill
\eject
\end{center}

\tableofcontents
\vfill\eject}
\setcounter{page}{1}
\def\thepage{\arabic{page}}
{\parindent=0pt 
\setcounter{section}{0}
\section*{Introduction}
\addcontentsline{toc}{section}{Introduction}

The higher algebraic $K$-theory of rings developed by D. Quillen \cite{Quil} constitutes an abstract powerful tool \cite{A-K-W} of algebraic topology which generalizes the lower $K$-groups among which the topological $K$-theory \cite{Ati1}, \cite{Ati2} introduced by A. Grothendieck is very popular among mathematicians and physicists.
\vskip 11pt

It is the aim of this paper {\bbf to revisit concretely the algebraic $K$-theory \cite{Kar1},\cite{Kar2} in the light of the recent developments \cite{Pie2}, \cite{Pie6} of the global program of Langlands\/} by taking into account the new algebraic interpretation of homotopy introduced here as a deformation of Galois representation given by compact(ified) algebraic groups \cite{Dol}.

More concretely, {\bbf we introduce\/}:
\Bean
\item {\bbf lower (algebraic) bilinear $K$-theories\/} (in the context of topological $K$-theory) {\bbf referring to homotopy and cohomotopy and related to the irreducible bilinear global program of Langlands\/};

\item {\bbf higher algebraic bilinear $K$-theories referring to homotopy and cohomotopy and related to the  reducible 
bilinear global program of Langlands\/};

\item {\bbf mixed lower and higher bilinear algebraic $KK$-theories related to the dynamical geometric bilinear global program of Langlands\/}.
\Ee
\vskip 11pt

{\bbf Chapter 1 deals with the universal algebraic structures of the Langlands global program.  These algebraic structures are abstract bisemivarieties $G^{(2n)}(F_{\o v}\times F_v)$\/} over products, right by left, $F_{\o v}\times F_v$, of sets of increasing archimedean completions.  These abstract bisemivarieties, in the heart of the Langlands global program, are universal in the sense that:
\Bean
\item they are (functional) representation spaces of the algebraic bilinear semigroups\linebreak $\GL_{2n}(\wt F_{\o v}\times \wt F_v)$ being the $2n$-dimensional representations of the products, right by left, 
$W^{ab}_{\wt F_{\o v}}\times
W^{ab}_{\wt F_{v}}$ of global Weil semigroups;

\item they have open coverings \cite{Har} by affine bisemivarieties 
$G^{(2n)}(\wt F_{\o v}\times \wt F_v)$
where $\wt F_{\o v}$ and $\wt F_{v}$ are symmetric increasing sets of Galois extensions;

\item they  generate cuspidal representations of 
$\GL_{2n}(\wt F_{\o v}\times \wt F_v)$ after a suitable toroidal compactification of these.
\Ee
\vskip 11pt

{\bbf Chapter 2 introduces bilinear versions \cite{Pie3} of the ``lower'' (algebraic) $K$-theory related to the bilinear global program of Langlands\/}.

As these $K$-theories are contravariant functors whose objects are abstract bisemivarieties (or bisemifields) and as they are defined with respect to homotopy groups, {\bbf it is logical to want to give an algebraic interpretation of homotopy allowing it to result from algebraic geometry\/}.

{\bbf The fundamental group can then be expressed in terms of deformations of Galois representations\/}.  Indeed, the equivalence classes of maps between the coefficient semiring $F_v$ and the real linear (semi)variety $G^{(2n)}(F_v)$ are {\bbf the homotopy classes \cite{Bott} corresponding to the classes of the quantum homomorphism $Qh_{v+\ell\to v}:F_{v+\ell}\to F_v$ sending $F_v$ into the deformed global coefficient semirings $F_{v+\ell}$ obtained from $F_v$ by adding ``$\ell$'' transcendental quanta\/} covered by irreducible closed algebraic subsets according to section~2.6.

Let then $fh_\ell:F_{v+\ell } \to G^{(2n)}(F_v)$ be a continuous map in such a say that:
\[
FH : \quad F_v\times I\To G^{(2n)}(F_v)\; , \qquad I=[0,1]\;,
\]
be the homotopy map of $fh=FH(x,0)$ with $FH(x,1)=fh_\ell $.

{\bbf The homotopy classes\/}, corresponding to the classes of the quantum homomorphism 
$Qh_{v+\ell} \to v$, {\bbf are characterized by integers ``$\ell $'' which are in one-to-one correspondence with the values of the parameter $t\in[0,1]$ of the homotopy\/}.

Taking into account the existence of {\bbf a cohomotopy of which classes are the inverse equivalence classes of the corresponding homotopy\/}, 
we see that the set of homotopy classes, being equivalence classes of maps between the set 
$\Omega (L_{v^1},G^{(2n)}(F_v))$ of oriented paths and the semivariety 
$G^{(2n)}(F_v)$, forms {\bbf a group, noted $\Pi _1(G^{(2n)}(F_v),L_{v^1})$ in the big point $L_{v^1}$, depending on deformations of Galois compact representations of these paths corresponding to the increase of these  by a(n) (in)finite number of transcendental (or algebraic) quanta\/}.

{\bbf The semigroup $\Pi _{2i}(G^{(2n)}(F_v),L_{v^1})$ of homotopy classes of maps\/}
\[
{}_sfh^{2i}_\ell : \qquad s^{2i}_{(\ell )}\quad \To G^{(2n)}\quad (F_v)\;,
\]
sending the base point of the $2i$-sphere $S^{2i}$ to the base point 
$L_{v^i_{(j)}}$ of $G^{(2n)}(F_v)$, 

or equivalently, of maps
\[
{}_cfh^i_\ell :\qquad[0,1]\quad \To \quad G^{(2n)}(F_v)\;,
\]
from the $i$-cube $[0,1]^i$ to the semivariety $G^{(2n)}(F_v)$, 
{\bbf results from the deformations of the Galois compact representation of the semigroup
$\GL_{(2i)}(\wt F_v)$ given  by the kernels $G^{(2i)}(\delta F_{v+\ell })$ of the maps\/}:
\begin{eqnarray*}
\GD ^{2i}_\ell : \qquad G^{2i}(F_{v+\ell }) &\quad \To \quad &G^{(2i)}(F_v)\;, \qquad \forall\ \ell \ , \; 1\leq \ell \leq \infty \;,\\
t^{2i} \quad &  \to & \quad \ell ^{2i}\;,
\end{eqnarray*}
in such a way that the $2i$-th powers if the integers ``$\ell $'' be in one-to-one correspondence with the $2i$-th powers of the values of the parameter $t\in [0,1]$.

Similarly, {\bbf the cohomotopy semigroup, noted $\Pi ^{2i}(G^{(2n)}(F_v),L_{v^1})$, is defined by classes resulting from inverse deformations
$(\GD ^{(2i)}_\ell )^{-1}$ of the Galois representations of $\GL_{(2i)}(\wt F_v)$\/}.

If $G^{(2n)}(F_{\o v})$ denotes the semivariety dual of $G^{(2n)}(F_v)$, then the bilinear homotopy \resp{cohomotopy} semigroup will be given by
$\Pi _{2i}(G^{(2n)}(F_{\o v}\times F_v))$
\resp{$\Pi ^{2i}(G^{(2n)}(F_{\o v}\times F_v))$}
in such a way that its classes of (bi)maps:
\[
{}_cfh^{2i}_\ell \times _{(D)}
{}_cfh^{2i}_h :\qquad
[0,1]^{2i}_\ell \times_{(D)} [0,1]^{2i}_\ell \quad \To\quad 
G^{(2n)}(F_{\o v}\times F_v)
\]
result from the \resp{inverse} deformations of the Galois representations of the bisemivariety $G^{(2i)}(\wt F_{\o v}\times  \wt F_v)$.

It is then natural to {\bbf associate to 
$\Pi _{2i}(G^{(2n)}(F_{\o v}\times F_v))$
\resp{to $\Pi ^{2i}(G^{(2n)}(F_{\o v}\times F_v))$} 
the $\Pi $-cohomology
\resp{the $\Pi $-homology} corresponding to the group homomorphism of Hurewicz\/}:
\begin{align*}
hH\RL:\quad \Pi _{2i}(G^{(2n)}(F_{\o v}\times F_v))
&\To
H^{2i}(G^{(2n)}(F_{\o v}\times F_v),\zit\times_{(D)}\zit )\\
\rresp{hcH\RL:\quad \Pi ^{2i}(G^{(2n)}(F_{\o v}\times F_v))
&\To
H _{2i}(G^{(2n)}(F_{\o v}\times F_v),\zit\times_{(D)}\zit )},
\end{align*}
{\bbf where the entire bilinear cohomology\/}
$H^{2i}(G^{(2n)}(F_{\o v}\times F_v),\zit\times_{(D)}\zit)$
\resp{homology\linebreak $H_{2i}(G^{(2n)}(F_{\o v}\times F_v),\zit\times_{(D)}\zit)$} 
{\bbf refers to a bisemilattice deformed by the homotopy \resp{cohomotopy}
classes of maps of
$\Pi _{2i}(G^{(2n)}(F_{\o v}\times F_v))$
\resp{$\Pi ^{2i}(G^{(2n)}(F_{\o v}\times F_v))$}\/}.

The topological (bilinear) $K$-theory
$K^{2i}(G^{(2n)}(F_{\o v}\times F_v))$ of vector bibundles with base
$G^{(2n)}(F_{\o v}\times F_v)$
and (bi)fibre
$G^{(2n-2i+1)}(F_{\o v}\times F_v)$, introduced in section~2.19, leads to set up the Chern character, restricted to the class $c^i$, in the {\bbf bilinear $K$-cohomology\/} by the homomorphism \cite{Sus}, \cite{Wal}:
\[
c^i(G^{(2n)}(F_{\o v}\times F_v)):\quad
K^{2i}(G^{(2n)}(F_{\o v}\times F_v))\To
H^{2i}(G^{(2n)}(F_{\o v}\times F_v),G^{(2i)}(F_{\o v}\times F_v))\;.
\]
Then,
{\bbf the lower bilinear (algebraic) $K$-theory referring to homotopy \resp{cohomotopy} will be given by the equality \resp{homomorphism}\/}:
\begin{align*}
K^{2i}(G^{(2n)}(F_{\o v}\times F_v)) &\underset{(\to)}{=} \Pi _{2i}(G^{(2n)}(F_{\o v}\times F_v))\\
\rresp{K_{2i}(G^{(2n)}(F_{\o v}\times F_v)) &\underset{(\to)}{=} \Pi ^{2i}(G^{(2n)}(F_{\o v}\times F_v))}
\end{align*}
{\bbf in such a way that the homotopy \resp{cohomotopy} classes of maps of
$\Pi _{2i}(\cdot)$
\resp{$\Pi ^{2i}(\cdot)$}
are \resp{correspond to} (inverse) liftings of quantum deformations of the Galois representation $\GL_{2i}(\wt F_{\o v}\times \wt F_v)$\/}.
\vskip 11pt

{\bbf Chapter 3 introduces bilinear versions of the higher algebraic $K$-theory \cite{Mil} related to the reducible bilinear global program of Langlands\/}.

{\bbf An infinite bilinear semigroup 
$\GL (F_{\o v}\times F_v)$\/}, depending on the geometric dimensions ``$i$'', {\bbf is given by 
\[
\GL (F_{\o v}\times F_v)=\lim_{\overrightarrow {\;\; i \;\;}} \GL_{2i}(F_{\o v}\times F_v)
\]
and corresponds to the (partially) reducible (functional) representation space
$\RED(F)\REPSP (\GL_{2n=2+\dots+2i+\dots+2n_s}(F_{\o v}\times F_v))$ of the bilinear semigroup of matrices $\GL_{2n}(F_{\o v}\times F_v)$ with $n\to\infty $\/}.

An equivalent quantum infinite bilinear semigroup given by the set
\[
\L\{\GL^{Q}(F^{2i}_{\o v^1}\times F^{2i}_{v^1})=\lim\limits_{j=1\to r\to\infty }\GL^{(Q)}_j(F^{2i}_{\o v^1}\times F^{2i}_{v^1})\R\}_i\;,
\]
 depending primarily on the algebraic dimension ``$j$' and based on the unitary representation space of
$\GL_{2n}(F_{\o v}\times F_v)$, is also introduced in sections~3.4 to 3.6.

{\bbf The classifying bisemispace
$\BGL(F_{\o v}\times F_v)$ of 
$\GL(F_{\o v}\times F_v)$, associated with the partition $2n=2+\dots+2i+\dots 2n_s$, $n\to\infty $\/}, of the integer $2n$, is defined as the base bisemispace of all equivalence classes of deformations of the Galois representation of 
$\GL(\wt F_{\o v}\times \wt F_v)$ given by the kernels
$\GL(\delta F_{\o {v+\ell }}\times \delta F_{v+\ell })$ of the maps:
\[
\GD _\ell :\qquad
\GL ( F_{\o {v+\ell }}\times  F_{v+\ell })\quad \To\quad 
\GL ( F_{\o {v}}\times  F_{v})\;, \qquad 1\leq \ell \leq \infty \;.
\]
{\bbf The ``plus'' constructed of Quillen\/}, adapted to the bilinear case of the Langlands global program, leads to consider the map:
\[
\BG (1):\qquad\BGL ( F_{\o {v}}\times  F_{v})\quad \To\quad 
\BGL ( F_{\o {v}}\times  F_{v})^+\;,
\]
in such a way that the classifying bisemispace 
$\BGL( F_{\o {v }}\times  F_{v })^+$ is the base bisemispace of all equivalence classes of one-dimensional deformations of the Galois compact representation of
$\GL(\wt F_{\o {v }}\times \wt F_{v})$ given by the kernels $\L\{\GL^{(1)}(\delta F_{\o {v+\ell }}\times \delta F_{v+\ell })\R\}_\ell $ of the maps $\GD (1)_\ell $.

{\bbf The bilinear version of the algebraic $K$-theory \cite{Blo}, \cite{Ger}, \cite{Gil2}, of Quillen related to the Langlands global program is:
\[
K^{2i} (G^{(2n)}_{\rm red}(F_{\o {v}}\times  F_{v}))=\Pi _{2i}(\BGL (F_{\o {v}}\times  F_{v})^+)
\]
where $G^{(2n)}_{\rm red}(F_{\o {v}}\times  F_{v})=\RED(F)\REPSP(\GL_{2n}(F_{\o {v}}\times  F_{v}))$ is the (partially) reducible (functional) representation space of the bilinear semigroup of matrices $\GL_{2n}(F_{\o {v}}\times  F_{v})$\/} in such a way that the partition
$2n=2+\dots+2i+\dots+2n_s$ of the geometric dimension $2n$, $
n\leq \infty $, refers to the reducibility of $GL_{2n}(F_{\o {v}}\times  F_{v})$.

This higher version of the Langlands global program implies 
by the homotopy bisemigroup $\Pi _{2i}(\cdot\times \cdot)$
that the equivalence classes of $2i$-dimensional deformations of the Galois representations of the reducible bilinear semigroup $\GL_{2n}(F_{\o {v}}\times  F_{v})$ result from quantum homomorphisms of the global coefficient bisemiring $F_{\o {v}}\times  F_{v}$.

This higher algebraic $K$-theory, referring to homotopy, implies {\bbf the commutative diagram\/}:
\[\begin{psmatrix}[colsep=1cm,rowsep=.6cm]
K^{2i}(G^{(2n)}_{\rm red}(F_{\o v}\times F_v)) & & \Pi _{2i}(\BGL(F_{\o v}\times F_v)^+)\\[22pt]
& H^{2i}(G^{(2n)}_{\rm red}(F_{\o v}\times F_v),\zit\times_{(D)}\zit) & 
\psset{arrows=->,nodesep=5pt}
\everypsbox{\scriptstyle}
\ncline{1,1}{1,3}
\ncline{1,1}{2,2}<{\Bsm \text{Chern higher}\\ \text{restricted character}\\ \text{in $K$-cohomology}\Esm}
\ncline{1,3}{2,2}>{\Bsm \text{inverse restricted higher}\\ \text{$\Pi $-cohomology}\Esm}
\end{psmatrix}
\]
in such a way that:
\Be
\item {\bbf the classes of the entire bilinear cohomology 
$H^{2i}(G^{(2n)}_{\rm red}(F_{\o v}\times F_v),\zit\times_{(D)}\zit)$ refers to a bisemilattice deformed by the homotopy classes of maps of
$\Pi _{2i}(\BGL(F_{\o v}\times F_v)^+)$, corresponding to lifts of quantum deformations of the Galois representations of 
$\GL_{2n}^{\rm red}(\wt F_{\o v}\times \wt F_v)$\/}.

\item the restricted higher $K$-cohomology implies the restricted higher $\Pi $-cohomology.
\Ee

{\bbf The total higher algebraic $K$-theory relative to homotopy is\/}:
\[
K^*(G^{(2n)}_{\rm red}(F_{\o v}\times F_v))=\Pi _*( \BGL (F_{\o v}\times F_v)^+)
\]
where ``$*$'' refers to the partition $2n=2+\dots+2i+\dots+2n_s$ of $2n$.

Similarly,  the higher bilinear algebraic $K$-theory relative to cohomotopy is given by the equality:
\[
K_{2i}(G^{(2n)}_{\rm red}(F_{\o v}\times F_v))=\Pi ^{2i}( \BGL ( F_{\o v}\times F_v)^+)\;,
\]
where $\Pi ^{2i}( \BGL ( F_{\o v}\times F_v)^+)$ are the cohomotopy equivalence classes of $2i$-dimensional deformations of the Galois reducible representations of
$\GL_{2n} ( \wt F_{\o v}\times \wt F_v)$ and its total version is:
\[
K_*(G^{(2n)}_{\rm red}(F_{\o v}\times F_v))=\Pi ^{*}( \BGL ( F_{\o v}\times F_v)^+)\;.
\]
\vskip 11pt

{\bbf Chapter 4 deals with  mixed higher bilinear algebraic $KK$-theories \cite{B-D-F}, \cite{Kas}, \cite{Jan} related to the Langlands dynamical bilinear global program and referring to the existence of $K_*K^*$ functors on the categories of elliptic bioperators and (reducible) bisemisheaves $FG^{(2n)}_{\rm (red)} (F_{\o v}\times F_v)$\/}.

A prerequisite is {\bbf the introduction of a bilinear contracting fibre $\Fs^{2k}\RL(\TAN)$ in the tangent bibundle 
$\TAN(FG^{(2n)}(F_{\o v}\times F_v))$\/} implying the homology:
\[
H_{2k} ( FG^{(2i[2k])} ( F_{\o v}\times F_v),\Fs^{2k}\RL (\TAN ))
\simeq \Ad (F)\REPSP (\GL_{2k}(\rit\times\rit ))\]
in such a way that the homology \cite{Nov}, \cite{Kas} of the $2i$-dimensional bisemisheaf\linebreak
$FG^{(2i[2k])} ( F_{\o v}\times F_v)$ shifted on $2k$ dimensions, $k\leq i$, be given by the adjoint functional representation space 
$\Ad (F)\REPSP (\GL_{2k}(\rit\times\rit ))$ of 
$\GL_{2k}(\rit\times\rit )$.

In order to introduce a mixed homotopy bisemigroup, we have to precise what must be the cohomotopy bisemigroup corresponding to the homology
$H_{2k} (FG^{(2i[2k])} ( F_{\o v}\times F_v),\linebreak \Fs^{2k}\RL (\TAN ))$.

As a ``Galois'' cohomotopy bisemigroup refers to classes resulting from inverse deformations of Galois representations, {\bbf the searched cohomotopy bisemigroup\linebreak 
$\Pi ^{2k}( FG^{(2i[2k])} ( F_{\o v}\times F_v))$
must be described by classes resulting from inverse deformations of the differential Galois \cite{Car} representations of $\GL_{2k}(\rit\times\rit )$ and depending on the classes of deformations of the Galois representations of $\GL_{2k}(\wt  F_{\o v}\times \wt F_v)$\/}.

Consequently, {\bbf the mixed homotopy bisemigroup
$\Pi _{2i[2k]} ( FG^{(2n[2k])}( F_{\o v}\times F_v))$ of the shifted bisemisheaf
$FG^{(2n[2k])}( F_{\o v}\times F_v)$\/} under the action of differential bioperators
{\bbf will be given by the product\/}:
\[
\Pi _{2i[2k]}( FG^{(2n[2k])} ( F_{\o v}\times F_v))
= \Pi ^{2k}(F G^{(2i[2k])} ( F_{\o v}\times F_v))\times
\Pi _{2i}( FG^{(2n[2k])} ( F_{\o v}\times F_v))\]
where
$\Pi _{2i}( FG^{(2n[2k])} ( F_{\o v}\times F_v))$ is the homotopy bisemigroup of the bisemisheaf 
$( FG^{(2n)} ( F_{\o v}\times F_v))$ shifted in $2k$ dimensions.

The mixed bilinear semigroup homomorphim of Hurewicz, introducing a restricted $\Pi $-homology-$\Pi $-cohomology, will be given by:
\[ mhH: \quad
\Pi _{2i[2k]}( FG^{(2n[2k])} ( F_{\o v}\times F_v))\To
H^{2i-2k}( FG^{(2n[2k])} ( F_{\o v}\times F_v),\zit\times_{(D)}\zit )
\]
where $H^{2i-2k}( \cdot)$ is the entire mixed bilinear cohomology defined from:
\begin{multline*}
H^{2i-2k}( FG^{(2n[2k])} ( F_{\o v}\times F_v),FG^{(2i[2k])} ( F_{\o v}\times F_v))\\
=
H_{2k}(F G^{(2i[2k])} ( F_{\o v}\times F_v), \Fs^{2k} \RL ( \TAN ))
\times
\L[ H^{2k}( FG^{(2n[2k])} ( F_{\o v}\times F_v), FG^{(2k)} (F_{\o v}\times F_v ))\R]\\
\oplus H^{2i-2k}(F G^{(2n[2k])} ( F_{\o v}\times F_v), FG^{(2i-2k)} (F_{\o v}\times F_v))
\end{multline*}
as developed in sections~4.1 to 4.3.

Similarly, the Chern mixed restricted character in the $K$-homology-$K$-cohomology, corresponding to a bilinear version of the index theorem, is given by:
\[
c_{[k]}\cdot
c^i(G^{(2n[2k])} ( F_{\o v}\times F_v)):\quad
K^{2i-2k} (FG^{(2n)} ( F_{\o v}\times F_v))\To
H^{2i-2k} (FG^{(2n)} ( F_{\o v}\times F_v))
\]
where 
$K^{2i-2k} (FG^{(2n)} ( F_{\o v}\times F_v))$ is {\bbf the mixed topological (bilinear) $K$-theory expanded according to\/}:
\[
K^{2i-2k} (FG^{(2n)} ( F_{\o v}\times F_v))=
K_{2k} (FG^{(2i[2k])} ( F_{\o v}\times F_v))\times
K^{2i} (FG^{(2n)} ( F_{\o v}\times F_v))
\]
with $K_{2k} (FG^{(2i[2k])} ( F_{\o v}\times F_v))$ a topological bilinear contracting $K$-theory of contracting tangent bibundles.

{\bbf A mixed lower bilinear (algebraic) $K$-theory can then be defined by the equality \resp{homomorphism}\/}:
\[
K^{2i-2k} (FG^{(2n)} ( F_{\o v}\times F_v))\underset{(\to)}{=}
\Pi _{2i[2k]} (FG^{(2n)} ( F_{\o v}\times F_v))\;.\]

{\bbf The total Chern mixed character in the $K$-homology-$K$-cohomology is\/}:
\[
\ch^*_* (FG^{(2n[2k])} ( F_{\o v}\times F_v)):\quad
K_{2*} K^{2*} (FG^{(2n)} ( F_{\o v}\times F_v))\To
H_{2*} H^{2*} (FG^{(2n)} ( F_{\o v}\times F_v))\]
leading to {\bbf a simplified version of the (bilinear) version of the Riemann-Roch theorem\/} introduced in  proposition~4.9.

In order to develop a bilinear version of the higher algebraic mixed $KK$-theory, we have to introduce {\bbf a higher bilinear algebraic operator $K$-theory relative to cohomotopy\/}.

This implies the definition of the infinite bilinear classifying semisheaf $\BFGL (F_{\o v}\times F_v)$ over the 
classifying bisemispace $\BGL (F_{\o v}\times F_v)$ as the base bisemispace of all equivalence classes of inverse deformations of the Galois differential representation of 
\[
\GL ( \rit \times\rit ) =\lim_{\overrightarrow{\; k\;}}\GL_{2k}(\rit \times\rit )
\] acting on
$\BGL ( F_{\o v}\times F_v)$.

The mixed classifying bisemisheaf $\BFGL((F_{\o v}\otimes\rit)\times(F_v\otimes\rit))$ then results from the biaction of $\BFGL(\rit\times\rit)$ on $\BFGL(F_{\o v}\times F_v)$.

{\bbf The bilinear version of the mixed higher algebraic $KK$-theory related to the Langlands dynamical bilinear global program is\/}:
\begin{multline*}
K_{2k} ( FG^{(2n)}_{\rm red} ( \rit \times\rit ))\times
K^{2i} (FG^{(2n)}_{\rm red} (F_{\o v}\times F_v))\\
=
\Pi ^{2k} ( \BFGL ( \rit \times\rit )^+)\times
\Pi _{2i} (\BFGL (F_{\o v}\times F_v)^+)
\end{multline*}
written in condensed form according to:
\[
K^{2i-2k} ( FG^{(2n)}_{\rm red} (F_{\o v}\otimes \rit ) \times ( F_{v}\otimes \rit ))
=
\Pi _{2i[2k]} ( \BFGL ( F_{\o v}\otimes \rit ) \times ( F_{v}\otimes \rit )^+ )
\]
in such a way that the bilinear contracting $K$-theory
$K_{2k} (F G^{(2n)}_{\rm red} ( \rit \times\rit ))$, responsible for a  differential biaction, acts on the $K$-theory
$K^{2i} (FG^{(2n)}_{\rm red} (F_{\o v}\times F_v))$ of the reducible functional representation space 
$ FG^{(2n)}_{\rm red} ( F_{\o v}\times F_v )$ of the bilinear semigroup
$\GL_{2n} (F_{\o v}\times F_v)$ in one-to-one correspondence with the biaction of the cohomotopy bisemigroup\linebreak
$\Pi ^{2k} ( \BFGL ( \rit \times\rit )^+)$ of the ``$+$'' classifying bisemispace
$\BFGL ( \rit \times\rit )^+$.

Finally, {\bbf the bilinear version of the total mixed higher algebraic $KK$-theory related to the dynamical reducible global program of Langlands is\/}:
\begin{multline*}
\qquad K_* (F G^{(2n)}_{\rm red} ( \rit  \times \rit )) \times K^* (F G^{(2n)}_{\rm red} ( F_{\o v}\times F_v ))\\
=
\Pi ^* ( \BFGL ( \rit  \times \rit )^+) \times \Pi _* ( \BFGL ( F_{\o v}\times F_v )^+)\;.\qquad 
\end{multline*}

\section[Universal algebraic structures of the Langlands global program]{Universal algebraic structures of the Langlands\linebreak global program}

\subsection{Pseudoramified infinite archimedean places of number fields}

\Bean
\item Let $\wt F$ denote a set of finite extensions of a number field $k$ of characteristic $0$: {\bbf $\wt F$ is assumed to be a set of symmetric splitting fields\/} composed of the right and left algebraic extension semifields $\wt F_R$ and $\wt F_L$ being in one-to-one correspondence.
$\wt F_L$ \resp{$\wt F_R$} is then composed of the set of complex \resp{conjugate complex} simple roots of the polynomial ring $k[x]$.

If the algebraic extension fields are real, then the symmetric splitting fields $\wt F^+$ are composed of the left and right symmetric splitting semifields  $\wt F^+_L$ and $\wt F^+_R$ being given respectively by the set of positive and symmetric negative simple real roots.

\item {\bbf The left and right equivalence classes of infinite archimedean completions\/} of
$\wt F_L$ \resp{$\wt F_R$} are the left and right infinite complex places
$\omega =\{\omega _1,\dots,\omega _j,\dots,\linebreak \omega _r\}$
\resp{$\o \omega =\{\o\omega _1,\dots,\o\omega _j,\dots,\o\omega _r\}$}.
In the real case, the infinite places are similarly
$v =\{v_1,\dots,v _j,\dots,v _r\}$
\resp{$\o v =\{\o v _1,\dots,\o v _j,\dots,\o v _r\}$}.

\item All these (pseudoramified) completions, corresponding to transcendental extensions, proceed from the associated algebraic extensions by a suitable isomorphism of compactification \cite{Pie7} and are built from irreducible subcompletions 
$F_{v^1_j}$
\resp{$F_{\o v^1_j}$}
characterized by a transcendence degree
\[
\tr\cdot d\cdot F_{v^1_j}/k=
\tr\cdot d\cdot F_{\o v^1_j}/k=N\]
equal to
\[[\wt F_{v^1_j}:k]=
[\wt F_{\o v^1_j}:k]=N
\]
which is the Galois extension degree of the associated algebraic closed subsets
$\wt F_{v^1_j}$ and $\wt F_{\o v^1_j}$.

{\bbf All these irreducible subcompletions \resp{(sub)extensions} are assumed to be transcendental \resp{algebraic} quanta \cite{Pie1}\/}.

\item {\bbf The pseudoramified real extensions are characterized by degrees\/}:
\[
[\wt F_{v_j}:k]=
[\wt F_{\o v_j}:k]=*+j\ N\;,\qquad 1\leq j\leq r\leq \infty \;,\]
which are integers modulo $N$, $\zit/N\ \zit$,

where:
\Bi
\item $\wt F_{v_j}$ and $\wt F_{\o v_j}$ are extensions corresponding to  completions respectively of the $v_j$-th and $\o v_j$-th symmetric real infinite pseudoramified places;
\item $*$ denotes an integer inferior to $N$.
\Ei

Similarly, {\bbf the pseudoramified complex extensions 
$\wt F_{\omega _j}$ and $\wt F_{\o \omega _j}$
corresponding to the completions 
$F_{\omega _j}$ and $  F_{\o \omega _j}$ at the infinite places $\omega _j$ and $\o\omega _j$\/}
are characterized by extension degrees
\[
[\wt F_{\omega _j}:k]=
[\wt F_{\o \omega _j}:k]=(*+j\ N)\ m^{(j)}\;,
\]
where $m^{(j)}=\sup(m_j+1)$ is the multiplicity of the $j$-th real extension covering its $j$-th complex equivalent.

Let then $\wt F_{v_{{j,m_j}}}$
\resp{$\wt F_{\o v_{{j,m_j}}}$} denote a \lr pseudoramified real extension equivalent to
$\wt F_{v_j}$ 
\resp{$\wt F_{\o v_j}$}.

\item {\bbf The corresponding pseudounramified real extensions\/}
$\wt F^{nr}_{v_{{j,m_j}}}$ and
$\wt F^{nr}_{\o v_{{j,m_j}}}$ are characterized by their global class residue degrees:
\[
f_{v_j}=[\wt F^{nr}_{v_{{j,m_j}}}:k]=j \qquad \text{and} \qquad
f_{\o v_j}=[\wt F^{nr}_{\o v_{{j,m_j}}}:k]=j\;.\]

\item Let $\wt F_{v_{j,m_j}}$ \resp{$\wt F_{\o v_{j,m_j}}$}
be the real pseudoramified extension corresponding to the completion
$  F_{v_{j,m_j}}$ and let
$\wt F^{nr}_{v_{j,m_j}}$
\resp{$\wt F^{nr}_{\o v_{j,m_j}}$}
be the corresponding pseudounramified extension.

Let $\Gal (\wt F_{v_{j,m_j}}/k)$ \resp{$\Gal (\wt F_{\o v_{j,m_j}}/k)$} and
$\Gal (\wt F^{nr}_{v_{j,m_j}}/k)$ \resp{$\Gal (\wt F^{nr}_{\o v_{j,m_j}}/k)$} be the associated Galois subgroups.

The {\bbf corresponding global Weil subgroups\/}
$W (\dot{\wt F}_{v_{j,m_j}}/k)=\Gal (\dot{\wt F}_{v_{j,m_j}}/k)$
\resp{$W (\dot{\wt F}_{\o v_{j,m_j}}/k)=\Gal (\dot{\wt F}_{\o v_{j,m_j}}/k)$}
are the Galois subgroups of the real pseudoramified extensions
$\dot{\wt F}_{v_{j,m_j}}/k$ \resp{$\dot{\wt F}_{\o v_{j,m_j}}/k$} characterized by extension degrees $d=0\ \mod N=jN$.

The global inertia subgroup
$I_{\dot{\wt F}_{v_{j,m_j}}}$ 
\resp{$I_{\dot{\wt F}_{\o v_{j,m_j}}}$}
of $W(\dot{\wt F}_{v_{j,m_j}}/k)$
\resp{of $W(\dot{\wt F}_{\o v_{j,m_j}}/k)$}
is defined by
\begin{align*}
I_{\dot{\wt F}_{v_{j,m_j}}} &= W(\dot{\wt F}_{v_{j,m_j}}/k) \Big/ W(\dot{\wt F}^{nr}_{v_{j,m_j}}/k)\\
\rresp{I_{\dot{\wt F}_{\o v_{j,m_j}}} &= W(\dot{\wt F}_{\o v_{j,m_j}}/k) \Big/ W(\dot{\wt F}^{nr}_{\o v_{j,m_j}}/k)}
\end{align*}
and, having an order $N$, is considered as the subgroup of inner automorphisms of Weil (or Galois).  Remark that 
$W(\wt F^{nr}_{v_{j,m_j}}/k)=
\Gal(\dot{\wt F}^{nr}_{v_{j,m_j}}/k)$.

Finally, the global Weil (semi)group
\[
W^{ab}_{\dot{\wt F}_v}=\{ W(\dot{\wt F}_{v_{j,m_j}}/k)\}_{{j,m_j}}\quad
\rresp{W^{ab}_{\dot{\wt F}_{\o v}}=\{ W ( \dot{\wt F}_{\o v_{j,m_j}}/k)\}_{{j,m_j}}},
\]
is the semigroup of all global Weil subsemigroups of real pseudoramified extensions
$\dot{\wt F}_{v_{j,m_j}}$
\resp{$\dot{\wt F}_{v_{j,m_j}}$} where
\[
\dot{\wt F}_{v}=\{ \dot{\wt F}_{v_1},\dots,\dot{\wt F}_{v_{j,m_j}},\dots\}\quad
\rresp{\dot{\wt F}_{\o v}=\{ \dot{\wt F}_{\o v_1},\dots,\linebreak \dot{\wt F}_{\o v_{j,m_j}},\dots\}}.\]
\Ee
\vskip 11pt

\subsection {Proposition ((Semi)groups of automorphisms of archimedean completions)}

{\em 
{\bbf The Galois sub(semi)groups 
$\Gal (\wt F_{v_{j,m_j}}/k)$
\resp{$\Gal (\wt F_{\o v_{j,m_j}}/k)$} of the extensions
$\wt F_{v_{j,m_j}}$
\resp{$\wt F_{\o v_{j,m_j}}$}
 are in one-to-one correspondence with the sub(semi)\-groups of automorphisms
 $\Aut_k( F_{v_{j,m_j}})$
\resp{$\Aut_k( F_{\o v_{j,m_j}})$}
 of the corresponding completions\/} (or transcendental extensions)
 $F_{v_{j,m_j}}$
\resp{$F_{\o v_{j,m_j}}$} in such a way that
\Bean

\item the completions $F_{v_{j,m_j}}$
\resp{$F_{\o v_{j,m_j}}$} are characterized by a transcendence degree
$\tr\cdot d\cdot F_{v_{j,m_j}}/k$
\resp{$\tr\cdot d\cdot F_{\o v_{j,m_j}}/k$}
 verifying
 \begin{align*}
 \tr\cdot d\cdot F_{v_{j,m_j}}&=[\wt F_{v_{j,m_j}}:k]=*+j\cdot N\\
 \rresp{\tr\cdot d\cdot F_{\o v_{j,m_j}}&=[\wt F_{\o v_{j,m_j}}:k]=*+j\cdot N}
 \end{align*}
 which is the cardinal number of the transcendence base of 
 $ F_{v_{j,m_j}}$
 \resp{$ F_{\o v_{j,m_j}}$} over $k$.
 
 \item there is a one-to-one correspondence between the set of all  transcendental extension subfields and the set of all sub(semi)groups of automorphisms of these.
 \Ee
 }
 \vskip 11pt
 
 \bpr
 \Bi
 \item The completions 
$ F_{v_{j,m_j}}$
\resp{$ F_{\o v_{j,m_j}}$}
 are transcendental extensions since they are generated from the corresponding algebraic extensions
 $ \wt F_{v_{j,m_j}}$
\resp{$ \wt F_{\o v_{j,m_j}}$}
by isomorphisms of compactifications
\[
c_{v_{j,m_j}}: \quad \wt F_{v_{j,m_j}}\To  F_{v_{j,m_j}}\qquad
\rresp{c_{\o v_{j,m_j}}: \quad \wt F_{\o v_{j,m_j}}\To  F_{\o v_{j,m_j}}}\]
sending
$\wt F_{v_{j,m_j}}$
\resp{$ \wt F_{\o v_{j,m_j}}$} by embedding into their compact isomorphic images
$ F_{v_{j,m_j}}$
\resp{$ F_{\o v_{j,m_j}}$}
which are closed compact subsets of $\rit_+$ \resp{$\rit_-$}.

As a result, certain points of these completions do not belong to algebraic extensions and correspond to transcendental extensions.

\item 		If the degrees of the Galois sub(semi)groups correspond to the class zero of the integers modulo $N$, i.e. if they are equal to $d=0\ \mod N$
, then, these Galois sub(semi)groups are global Weil sub(semi)groups of extensions
$\dot{\wt F}_{v_{j,m_j}}$
\resp{$\dot{\wt F}_{\o v_{j,m_j}}$}
 constructed from sets of ``$j$'' algebraic quanta (which are algebraic closed subsets characterized by an extension degree equal to $N$).
 
 By the isomorphism of compactification
 $c_{v_{j,m_j}}$
 \resp{$c_{\o v_{j,m_j}}$},
 the ``$j$'' (non compact) algebraic quanta of
 $\dot{\wt F}_{v_{j,m_j}}$
\resp{$\dot{\wt F}_{\o v_{j,m_j}}$} are sent into {\bbf the corresponding compactified ``$j$'' transcendental (compact) quanta forming the completions
$\dot{ F}_{v_{j,m_j}}$
\resp{$\dot{F}_{\o v_{j,m_j}}$} also noted simply
$ F_{v_{j,m_j}}$
\resp{$F_{\o v_{j,m_j}}$}\/}.

\item {\bbf The compact archimedean completion
$ F_{v_{j,m_j}}$
\resp{$F_{\o v_{j,m_j}}$} can also be viewed as resulting from the sub(semi)-group of automorphisms\linebreak
$\Aut_k( F_{v_{j,m_j}})$
\resp{$\Aut_k(F_{\o v_{j,m_j}})$}\/} in such a way that:
\Bean
\item $\Aut_k( F_{v_{j,m_j}})$
\resp{$\Aut_k(F_{\o v_{j,m_j}})$} is the compact sub(semi)group of the automorphisms of order ``$j$'' of a transcendental quantum
$F_{v^1_{j,m_j}}$
\resp{$F_{\o  v^1_{j,m_j}}$};

\item $\Aut_k( F_{v_{j,m_j}})$
\resp{$\Aut_k(F_{\o v_{j,m_j}})$} is a semigroup of reflections \cite{Dol}
of a transcendental quantum.
\Ee

\item It is the evident that:
\Bean
\item $\Aut_k( F_{v_{j,m_j}})\simeq \Gal (\wt  F_{v_{j,m_j}}/k)$
\resp{$\Aut_k( F_{\o v_{j,m_j}})\simeq \Gal (\wt  F_{\o v_{j,m_j}}/k)$};

\item as in the Galois case, there is a one-to-one correspondence between all  transcendental extension subfields
\begin{align*}
&F_{v_1}\subset \dots \subset F_{v_{j,m_j}}\subset \dots\subset  F_{v_{r,m_r}}\\
\rresp{&F_{\o v_1}\subset \dots \subset F_{\o v_{j,m_j}}\subset\dots \subset F_{\o v_{r,m_r}}}
\end{align*}
and the set of all normal sub(semi)groups of automorphisms of these:
\begin{align}
&\Aut_k(F_{v_1})\subset \dots \subset\Aut_k(F_{v_{j,m_j}})\subset \dots \subset\Aut_k(F_{v_{r,m_r}})\notag\\
\rresp{&\Aut_k(F_{\o v_1})\subset \dots \subset\Aut_k(F_{\o v_{j,m_j}})\subset \dots \subset\Aut_k(F_{\o v_{r,m_r}})}.\tag*{\eop}
\end{align}
\Ee
\Ei
\vskip 11pt

\subsection{Real  algebraic bilinear semigroups}

\Bean
\item Let $B_{\wt F_v}$
\resp{$B_{\wt F_{\o v}}$}
 be a \lr division semialgebra of real dimension $2n$ over the set
  ${\wt F_v}$
\resp{${\wt F_{\o v}}$} of increasing real pseudoramified extensions
 ${\wt F_{v_{j,m_j}}}$
 \resp{${\wt F_{\o v_{j,m_j}}}$} of $k$.
 
 Then, $B_{\wt F_v}$
\resp{$B_{\wt F_{\o v}}$}, which is a \lr vector semispace restricted to the upper \resp{lower} half space, is isomorphic to the semialgebra of Borel upper \resp{lower} triangular matrices:
\[
B_{\wt F_v}\simeq T _{2n}(\wt F_v)\qquad
\rresp{B_{\wt F_{\o v}}\simeq T^t _{2n}(\wt F_{\o v})}.\]

{\bbf This allows to define the algebraic bilinear semigroup of matrices\linebreak 
$\GL_{2n}(\wt F_v\times\wt F_{\o v})$ by:
\[
B_{\wt F_v}\otimes B_{\wt F_v}\simeq 
T ^t_{2n}(\wt F_{\o v})\times
T _{2n}(\wt F_v)\equiv \GL_{2n}(\wt F_{\o v}\times\wt F_{v})\]
}
in such a way that its representation (bisemi)space is given by the tensor product
$\wt M_{v_R}\otimes\wt M_{v_L}$ of a right $T ^t_{2n}(\wt F_{\o v})$-semimodule
$\wt M_{v_R}$ by a left $T _{2n}(\wt F_v)$-semimodule $\wt M_{v_L}$.

The $\GL_{2n}(\wt F_{\o v}\times\wt F_{v})$-bisemimodule
$\wt M_{v_R}\otimes\wt M_{v_L}$ is an algebraic bilinear (affine) semigroup noted
$G^{(2n)}(\wt F_{\o v}\times\wt F_{v})$ verifying the commutative diagram:
\[\begin{psmatrix}[colsep=1cm,rowsep=.6cm]
\GL_{2n}(\wt F_{\o v}\times \wt F_{v}) & &\wt M_{v_R}\otimes\wt M_{v_L}\equiv G^{(2n)}(\wt F_{\o v}\times \wt F_{v})\\[15pt]
\GL(\wt M_{v_R}\otimes\wt M_{v_L})& &
\psset{arrows=->,nodesep=5pt}
\ncline{1,1}{1,3}
\ncline{1,1}{2,1}
\ncline{2,1}{1,3}
\end{psmatrix}
\]
where $\GL(\wt M_{v_R}\otimes\wt M_{v_L})$ is the bilinear semigroup of automorphisms of
$\wt M_{v_R}\otimes\wt M_{v_L}$.

Then, {\bbf $\GL(\wt M_{v_R}\otimes\wt M_{v_L})$ constitutes the $2n$-dimensional equivalent of the product
$W^{ab}_{\dot{\wt F}_{\o v}}\times
W^{ab}_{\dot{\wt F}_{v}}$ of the global Weil semigroups\/}  and the  bilinear algebraic semigroup $G^{(2n)}(\wt F_{\o v}\times\wt F_{v})$ becomes naturally the $2n$-dimensional (irreducible) representation space
$\Irr\Rep^{2n}_{W_{F^+\RL}}(W^{ab}_{\dot{\wt F}_{\o v}}\times W^{ab}_{\dot{\wt F}_{v}})$ of
$(W^{ab}_{\dot{\wt F}_{\o v}}\times W^{ab}_{\dot{\wt F}_{v}})$ in such a way that
\[
\Irr\Rep^{2n}_{W_{F^+\RL}}(W^{ab}_{\dot{\wt F}_{\o v}}\times W^{ab}_{\dot{\wt F}_{v}}):\quad
\GL (\wt M_{v_R}\otimes\wt M_{v_L})\To
G^{(2n)}(\wt F_{\o v}\times\wt F_{v})\]
implies the monomorphism:
\[
\wt \sigma _{v_R}\times\wt \sigma _{v_L}:\qquad
W^{ab}_{\dot{\wt F}_{\o v}}\times W^{ab}_{\dot{\wt F}_{v}}\quad \To\quad 
\GL_{2n}(\wt F_{\o v}\times\wt F_{v})\;.\]

\item The isomorphisms
\begin{align*}
 \Aut_k (F_{v_{j,m_j}}) &\simeq \Gal (\wt F_{v_{j,mj}}/k)\\
\rresp{\Aut_k (F_{\o v_{j,m_j}})&\simeq \Gal (\wt F_{\o v_{j,mj}}/k)},\qquad \forall\ j,m_j
\;,
\end{align*}
between the subgroups of automorphisms of the completions
$F_{v_{j,m_j}}$ \resp{$F_{\o v_{j,m_j}}$} and the corresponding Galois subgroups of the extensions $\wt F_{v_{j,m_j}}$ \resp{$\wt F_{\o v_{j,m_j}}$}, as developed in proposition~1.1.2, naturally leads to the commutative diagram:
\[\begin{psmatrix}[colsep=1cm,rowsep=.6cm]
W^{ab}_{\dot{\wt F}_{\o v}}\times W^{ab}_{\dot{\wt F}_{v}} 
&  G^{(2n)}(\wt F_{\o v}\times\wt F_{v})\\[15pt]
\Aut_k (F_{\o v})\times \Aut_k (F_{v})
& G^{(2n)}( F_{\o v}\times F_{v})
\psset{arrows=->,nodesep=5pt}
\everypsbox{\scriptstyle}
\ncline{1,1}{1,2}^{\wt \sigma _{v_R}\times\wt \sigma _{v_L}}
\ncline{2,1}{2,2}^{\sigma _{v_R}\times\sigma _{v_L}}
\ncline{1,1}{2,1}>{\rotatebox{90}{$\sim$}}
\ncline{1,2}{2,2}>{\rotatebox{90}{$\sim$}}
\end{psmatrix}
\]
where $\sigma _{v_R}\times\sigma _{v_L}$ is the monomorphism between the product, right by left, of the semigroups of automorphisms of the set of completions $F_{\o v}$ and $F_v$ and
{\bbf the complete locally compact (algebraic) bilinear semigroup 
$G^{(2n)}( F_{\o v}\times  F_{v})$\/}
defining an abstract bisemivariety.

Thus, the isomorphism $W^{ab}_{\dot{\wt F}_{\o v}}\times W^{ab}_{\dot{\wt F}_{v}} 
\overset{\sim }{\to} \Aut_k (F_{\o v})\times \Aut_k (F_{v})$ implies
{\bbf the  homomorphism 
$G^{(2n)}(\wt F_{\o v}\times\wt F_{v})\simeq G^{(2n)}(F_{\o v}\times F_{v})$ between bilinear semigroups such that the abstract bisemivariety
$G^{(2n)}( F_{\o v}\times F_{v})$ be covered by the algebraic (affine) semigroup
$G^{(2n)}(\wt F_{\o v}\times\wt F_{v})$\/}.

\item Let $G^{(2n)}(\wt F^{nr}_{\o v}\times\wt F^{nr}_{v})$ be the algebraic bilinear semigroup over the product of the sets of increasing pseudounramified extensions with
$\wt F^{nr}_{v}=\{ \wt F^{nr}_{v_1},\dots,\wt F^{nr}_{v_{j,m_j}},\dots\}$ and
$\wt F^{nr}_{\o v}=\{ \wt F^{nr}_{\o v_1},\dots,\wt F^{nr}_{\o v_{j,m_j}},\dots\}$.

Then, the kernel $\Ker (G^{(2n)}_{\wt F\to\wt F^{nr}})$ of the map:
\[
G^{(2n)}_{\wt F\to\wt F^{nr}}:\qquad
G^{(2n)}(\wt F_{\o v}\times\wt F_{v})\quad \To\quad 
G^{(2n)}(\wt F^{nr}_{\o v}\times\wt F^{nr}_{v})
\]
is {\bbf the smallest bilinear normal pseudoramified subgroup of
$G^{(2n)}(\wt F_{\o v}\times\wt F_{v})$\/}:
\[\Ker (G^{(2n)}_{\wt F\to\wt F^{nr}})=P^{(2n)}(\wt F_{\o v^1}\times\wt F_{v^1})\:,\]
i.e. {\bbf the minimal bilinear parabolic subsemigroup
$P^{(2n)}(\wt F_{\o v^1}\times\wt F_{v^1})$\/}
over the product
$(\wt F_{\o v^1}\times\wt F_{v^1})$ of the sets
\[ \wt F_{\o v^1}=\{ \wt F_{\o v^1_1},\dots, \wt F_{\o v^1_{j,m_j}},\dots,\wt F_{\o v^1_{r,m_r}}\} \quad \text{and} \quad
 \wt F_{v^1}=\{ \wt F_{v^1_1},\dots, \wt F_{v^1_{j,m_j}},\dots,\wt F_{v^1_{r,m_r}}\}\]
 of unitary archimedean pseudoramified extensions
 $\wt F_{\o v^1_{j,m_j}}$ and $\wt F_{v^1_{j,m_j}}$  in 
 $\wt F_{\o v_{j,m_j}}$ and  $\wt F_{v_{j,m_j}}$ respectively.
 
 \item At every infinite biplace $\o v_j\times v_j$ of $F_{\o v}\times F_v$ corresponds 
 {\bbf a conjugacy class
 $g^{(2n)}_{v\RL}[j]$ of the algebraic bilinear semigroup
 $G^{(2n)}(\wt F_{\o v}\times\wt F_{v})$\/}.  The number of representatives of
 $g^{(2n)}_{v\RL}[j]$ corresponds to the number of equivalent extensions of
 $\wt F_{\o v_j}\times\wt F_{v_j}$.
 
 So, we have the injective morphism:
 \[ I_{F-G_v}:\qquad \{ \wt F_{\o v_{j,m_j}}\times\wt F_{v_{j,m_j}} \}_{m_j}\To
\{G^{(2n)}( \wt F_{\o v_{j,m_j}}\times\wt F_{v_{j,m_j}}) \}_{m_j}\]
leading to the homeomorphism:
\[ \prod_{2n} ( \wt F_{\o v_{j,m_j}}\times\wt F_{v_{j,m_j}}) \simeq 
G^{(2n)}( \wt F_{\o v_{j,m_j}}\times\wt F_{v_{j,m_j}}) \]
where $G^{(2n)}( \wt F_{\o v_{j,m_j}}\times\wt F_{v_{j,m_j}}) $ is the $(j,m_j)$-th conjugacy class representative of\linebreak 
$G^{(2n)}( \wt F_{\o v_{j,m_j}}\times\wt F_{v_{j,m_j}}) $

\item {\bbf $G^{(2n)}( \wt F_{\o v}\times\wt F_{v}) $ acts on the bilinear parabolic subsemigroup
$P^{(2n)}( \wt F_{\o v^1}\times\wt F_{v^1}) $ by conjugation\/} in such a way that the number of conjugates of
$P^{(2n)}( \wt F_{\o v^1_{j}}\times\wt F_{v^1_{j}}) $ in
$G^{(2n)}( \wt F_{\o v_{j,m_j}}\times\wt F_{v_{j,m_j}}) $ is the index of the normalizer
$P^{(2n)}( \wt F_{\o v^1_{j}}\times\wt F_{v^1_{j}}) $ in $G^{(2n)}( \wt F_{\o v_{j,m_j}}\times\wt F_{v_{j,m_j}})$:
\[
\L|
G^{(2n)}( \wt F_{\o v_{j,m_j}}\times\wt F_{v_{j,m_j}}) :
P^{(2n)}( \wt F_{\o v^1_{j}}\times\wt F_{v^1_{j}}) \R|=j\;.\]

\item Let $\Out (G^{(2n)}( \wt F_{\o v}\times\wt F_{v}))= \Aut (G^{(2n)}( \wt F_{\o v}\times\wt F_{v}))\ \big /\ \Int (G^{(2n)}( \wt F_{\o v_{j}}\times\wt F_{v_{j}}))$ be the
(bisemi)group of (Galois) automorphisms of the algebraic bilinear semigroup\linebreak 
$G^{(2n)}( \wt F_{\o v}\times\wt F_{v}) $ where 
$\Int (G^{(2n)}( \wt F_{\o v}\times\wt F_{v}))$ is the (bisemi)group of (Galois) inner automorphisms.

As we have that
$\Int (G^{(2n)}( \wt F_{\o v}\times\wt F_{v}))=
\Aut (P^{(2n)}( \wt F_{\o v^1}\times\wt F_{v^1}))$, the bilinear parabolic semigroup
$P^{(2n)}( \wt F_{\o v^1}\times\wt F_{v^1})$ can be considered as the unitary irreducible representation  space of the algebraic bilinear semigroup
$\GL_{2n}(\wt F_{\o v}\times\wt F_{v})$ of matrices \cite{Pie2}.
\Ee
\vskip 11pt

\subsection{Complex  algebraic bilinear semigroups}

Similarly as it was done for the real case in section~1.3, let us consider the complex case and, especially:

\Bean
\item Let $B_{\wt F_{\omega }}$
\resp{$B_{\wt F_{\o \omega }}$} be a \lr division semialgebra of complex dimension $n$ over the set 
$\wt F_{\omega }=\{\wt F_{\omega_1},\dots, F_{\omega _{j,m_j}},\dots\}$
\resp{$\wt F_{\o\omega }=\{\wt F_{\o\omega_1},\dots, F_{\o\omega _{j,m_j}},\dots\}$}
of increasing complex pseudoramified extensions of $k$.

This allows to define the algebraic bilinear semigroup 
$\GL_n(\wt F_{\o \omega }\times\wt F_\omega )$ by:
\[
B_{\wt F_{\o\omega }}\otimes B_{\wt F_{\omega }}
\simeq 
\GL_n (\wt F_{\o \omega }\times\wt F_\omega )
\equiv T^t_n (\wt F_{\o\omega })\times T_n(\wt F_{\omega })
\]
in such a way that its representation space be given by the
$\GL_n(\wt F_{\o \omega }\times\wt F_\omega )$-bisemi\-module
$\wt M_{\omega _R}\otimes\wt M_{\omega _L}$ which is also a complex affine algebraic bilinear semigroup 
$G^{(2n)}(\wt F_{\o \omega }\times\wt F_\omega )$ homeomorphic to the complete (algebraic) bilinear semigroup\linebreak 
$G^{(2n)}(F_{\o \omega }\times F_\omega )$
over the sets of completions $F_{\o\omega }$ and $F_\omega $.

Let $\GL(\wt M_{\omega _R}\otimes\wt M_{\omega_L})$ denote the bilinear semigroup of automorphisms of
$(\wt M_{\omega _R}\otimes\wt M_{\omega_L})$ verifying:
\[\begin{psmatrix}[colsep=1cm,rowsep=.6cm]
\GL_{n}(\wt F_{\o \omega }\times \wt F_{\omega }) & &\wt M_{\omega _R}\otimes\wt M_{\omega _L}
\equiv G^{(2n)}(\wt F_{\o \omega }\times \wt F_{\omega })\\[15pt]
\GL(\wt M_{\omega _R}\otimes\wt M_{\omega _L})& &
\psset{arrows=->,nodesep=5pt}
\ncline{1,1}{1,3}
\ncline{1,1}{2,1}
\everypsbox{\scriptstyle}
\ncline{2,1}{1,3}^{\sim }
\end{psmatrix}
\]
because $\GL(\wt M_{\omega _R}\otimes\wt M_{\omega _L})$ constitutes the $n$-dimensional complex equivalent of the product
$W^{ab}_{\dot{\wt F}_{\o \omega }}\times
W^{ab}_{\dot{\wt F}_{\omega }}$
of the corresponding global Weil semigroups.

{\bbf So, $G^{(2n)}(\wt F_{\o \omega }\times\wt F_\omega )$ is the $n$-dimensional (or $2n$-dimensional real) complex (irreducible) representation space of
$W^{ab}_{\dot{\wt F}_{\o \omega }}\times
W^{ab}_{\dot{\wt F}_{\omega }}$\/} given by:
\[
\Irr\Rep^{2n}_{W_{F\RL}}(W^{ab}_{\dot{\wt F}_{\o \omega }}\times W^{ab}_{\dot{\wt F}_{\omega }}):\qquad
\GL (\wt M_{\omega _R}\otimes\wt M_{\omega _L})\quad \overset{\sim}{\To}\quad 
G^{(2n)}(\wt F_{\o \omega }\times\wt F_{\omega })\]
implying the morphism:
\[
\sigma _{\o \omega _R}\times\sigma _{\o \omega _L}:\qquad
W^{ab}_{\dot{\wt F}_{\o \omega }}\times W^{ab}_{\dot{\wt F}_{\omega }}\quad \To\quad 
\GL_{n}(\wt F_{\o \omega }\times\wt F_{\omega })\;.\]

\item At every biplace $(\o\omega _j\times\omega _j)$ of $(F_{\o\omega }\times F_\omega )$ corresponds {\bbf a conjugacy class $g^{(2n)}_{\wt\omega \RL}[j]$ of
$G^{(2n)}(\wt F_{\o \omega }\times\wt F_{\omega })$\/} leading to the injective morphism:
\[ I_{F-G_\omega }:\qquad \{ \wt F_{\o \omega _{j,m_j}}\times\wt F_{\omega _{j,m_j}} \}_{m_j}\To
\{G^{(2n)}( \wt F_{\o \omega _{j,m_j}}\times\wt F_{\omega _{j,m_j}}) \}_{m_j}\]
where $G^{(2n)}( \wt F_{\o \omega _{j,m_j}}\times\wt F_{\omega _{j,m_j}}) $ is the $(j,m_j)$-th conjugacy class representative of\linebreak 
$G^{(2n)}( \wt F_{\o \omega} \times\wt F_{\omega }) $.

\item $G^{(2n)}( \wt F_{\o \omega }\times\wt F_{\omega }) $ acts by conjugation on the
{\bbf bilinear parabolic semigroup\linebreak  
$P^{(2n)}( \wt F_{\o \omega^1 }\times\wt F_{\omega^1 }) $ which can be considered as the unitary irreducible representation space of the complex algebraic bilinear semigroup
$\GL_{n}( \wt F_{\o \omega }\times\wt F_{\omega }) $\/} of matrices because the (bisemi)group of (Galois) inner automorphisms of
$G^{(2n)}( \wt F_{\o \omega }\times\wt F_{\omega }) $ verifies:
\[
\Int (G^{(2n)}( \wt F_{\o \omega }\times\wt F_{\omega }))=
\Aut (P^{(2n)}( \wt F_{\o \omega ^1}\times\wt F_{\omega ^1}))\;.\]
\Ee
\vskip 11pt

\subsection{Inclusion of real (algebraic) bilinear semigroups into their complex equivalents}

The complex $\GL_{n}( \wt F_{\o \omega }\times\wt F_{\omega }) $-bisemimodule
$\wt M_{\omega _R}\otimes\wt M_{\omega _L}$ is the representation space of the algebraic bilinear semigroup of matrices  $\GL_{n}( \wt F_{\o \omega }\times\wt F_{\omega }) $.

Assume that each conjugacy class representative
$G^{(2n)}( \wt F_{\o \omega_j }\times\wt F_{\omega _j}) $ of
$G^{(2n)}( \wt F_{\o \omega }\times\wt F_{\omega }) \equiv \wt M_{\omega _R}\otimes\wt M_{\omega _L}$ is unique in the $j$-th class.

Then, the set $\{G^{(2n)}( \wt F_{\o \omega_j }\times\wt F_{\omega _j})\}^r_{j=1} $ of conjugacy class representatives of 
$G^{(2n)}( \wt F_{\o \omega }\times\wt F_{\omega }) $ is the representation (bisemi)space of {\bbf the restricted complex algebraic bilinear semigroup $\GL^{(\rm res)}_{n}( \wt F_{\o \omega }\times\wt F_{\omega }) $\/}.

As a result, each complex conjugate class representative
$G^{(2n)}( \wt F_{\o \omega_j }\times\wt F_{\omega _j})$ of
$G_{(\rm res)}^{(2n)}( \wt F_{\o \omega }\times\wt F_{\omega }) $ is covered by the $m_j$ real conjugacy class representatives
$G^{(n)}( \wt F_{\o v_{j,m_j} }\times\wt F_{v _{j,m_j}})$ of
$G^{(n)}( \wt F_{\o v }\times\wt F_{v })$ \cite{Pie2}.

So, the complex bipoints of
$G_{(\rm res)}^{(2n)}( \wt F_{\o \omega }\times\wt F_{\omega }) $  are in one-to-one correspondence with the real bipoints of 
$G^{(2n)}( \wt F_{\o v}\times\wt F_{v}) $ and we have the inclusion:
\[\begin{psmatrix}[colsep=1cm,rowsep=.6cm]
G^{(2n)}( \wt F_{\o v}\times\wt F_{v}\big/ G^{(n)}( \wt F_{\o v}\times\wt F_{v}) \simeq \wt M_{v _R}\otimes\wt M_{v _L}
&&
G_{(\rm res)}^{(2n)}( \wt F_{\o \omega }\times\wt F_{\omega })\equiv 
\wt M^{(\rm res)}_{\omega  _R}\otimes\wt M^{(\rm res)}_{\omega  _L}
\psset{arrows=->,nodesep=5pt}
\ncline[hookwidth=-1mm]{H->}{1,1}{1,3}
\end{psmatrix}
\]
where $\wt M^{(\rm res)}_{\omega  _L}$ \resp{$\wt M^{(\rm res)}_{\omega  _R}$} is the \lr restricted 
$T^{(\rm res)}_{n}(\wt F_\omega )$-semimodule 
\resp{$T^{t(\rm res)}_{n}(\wt F_{\o\omega} )$-semimodule}.

$G_{(\rm res)}^{(2n)}( \wt F_{\o \omega }\times\wt F_{\omega })$ is then said to be covered by
$G^{(2n)}( \wt F_{\o v }\times\wt F_{v})$.
\vskip 11pt

\subsection{Cuspidal representation of complex algebraic bilinear semigroups}

\Bean
\item Providing {\bbf a cuspidal representation\/} of the complex bilinear algebraic semigroup
$G_{(\rm res)}^{(2n)}( \wt F_{\o \omega }\times\wt F_{\omega })$ consists in finding a cuspidal form on
$G_{(\rm res)}^{(2n)}( \wt F_{\o \omega }\times\wt F_{\omega })${\bbf by summing the cuspidal subrepresentations of its conjugacy class representatives
$G^{(2n)}( \wt F_{\o \omega_j }\times\wt F_{\omega _j})$\/}.

Let then
\[
\gamma^T _{F_{\omega _j}}:\quad\wt F_{\omega _j}\To F^T_{\omega _j}\qquad
\rresp{\gamma^T _{F_{\o\omega _j}}:\quad\wt F_{\o\omega _j}\To F^T_{\o\omega _j}}
\]
be the {\bbf toroidal isomorphism\/} mapping each \lr extension
$\wt F_{\omega _j}$ \resp{$\wt F_{\o\omega _j}$} into its toroidal compact equivalent
$F^T_{\omega _j}$ \resp{$F^T_{\o\omega _j}$} which is a complex one-dimensional semitorus localized in the upper \resp{lower} half space.

Then, the morphism:
\begin{align*}
T_n(F^T_{\omega _j}): \qquad F^T_{\omega _j}\quad &\To \quad T^{(2n)}(F^T_{\omega _j})=T^{2n}_L[j]\\
\rresp{T_n(F^T_{\o\omega _j}): \qquad F^T_{\o\omega _j}\quad &\To \quad T^{(2n)}(F^T_{\o\omega _j})=T^{2n}_R[j]}
\end{align*}
of the respective fibre bundle sends $F^T_{\omega _j}$ \resp{$F^T_{\o\omega _j}$}
into the $n$-dimensional complex semitorus $T^{2n}_L[j]$ \resp{$T^{2n}_R[j]$} corresponding to the upper \resp{lower} conjugacy class representative 
$T^{(2n)}(F^T_{\omega _j})=G^{(2n)}(F^T_{\omega _j})$
\resp{$T^{(2n)}(F^T_{\o\omega _j})=G^{(2n)}(F^T_{\o\omega _j})$}.

So, we have a homeomorphism
$G^{(2n)}( \wt F_{\o \omega_j }\times\wt F_{\omega_j })\simeq G^{(2n)}(  F^T_{\o \omega_j }\times F^T_{\omega _j})$ between the conjugacy class representative
$G^{(2n)}( \wt F_{\o \omega _j}\times\wt F_{\omega_j })$ of
$G^{(2n)}( \wt F_{\o \omega }\times\wt F_{\omega })$ and the conjugacy class representative
$G^{(2n)}(  F^T_{\o \omega _j}\times F^T_{\omega _j})$ of
$G^{(2n)}(  F^T_{\o \omega }\times F^T_{\omega })$ where
\[
F^T_{\omega }=\{ F^T_{\omega _1},\dots,F^T_{\omega _j},\dots\}\quad
\rresp{F^T_{\o \omega }= \{ F^T_{\o\omega _1},\dots, F^T_{\omega _j},\dots \}}.\]

\item {\bbf Every \lr function on the conjugacy class representative
$G^{(2n)}(  F^T_{\omega _j})$
\resp{$G^{(2n)}(  F^T_{\o \omega _j})$} is a function \resp{cofunction}\/}
$\phi _L(T^{2n}_L[j])$
\resp{$\phi _R(T^{2n}_R[j])$}
{\bbf on the complex semitorus $T^{2n}_L[j]$ \resp{$T^{2n}_R[j]$} having the analytic development\/}:
\[
\phi _L(T^{2n}_L[j])=\lambda ^\half(2n,j)\ e^{2\pi ijz}\quad
\rresp{\phi _R(T^{2n}_R[j])=\lambda ^\half(2n,j)\ e^{-2\pi ijz}}\]
where:
\Bi
\item $\vec z=\sum\limits^{2n}_{d=1}z_d\ \vec{e_d}$ is a complex point of $G^{(2n)}(  F^T_{\omega _j})$;

\item $\lambda (2n,j)$ is a product of the eigenvalues of the $j$-th coset representative of the product, right by left, of Hecke operators \cite{Pie2}.
\Ei
\item This \lr function $\phi _L(T^{2n}_L[j])$ \resp{$\phi _R(T^{2n}_R[j])$} constitutes the cuspidal representation
$\Pi ^{(j)}(G^{(2n)}( \wt F_{\omega _j}))$
\resp{$\Pi ^{(j)}(G^{(2n)}( \wt F_{\o\omega _j}))$} of the $j$-th conjugacy class representative of
$G^{(2n)}_{(\rm res)}( \wt F_{\omega })$
\resp{$G^{(2n)}_{(\rm res)}( \wt F_{\o\omega })$}
 in such a way that {\bbf the cuspidal biform of
 $\GL_{n}^{(\rm res)}( F_{\o\omega }\times  F_{\omega })$ is given by the Fourier biseries\/}:
 \[
 \Pi  (\GL_{n}^{(\rm res)}(\wt F_{\o\omega _\oplus}\times_D\wt F_{\omega _\oplus}))=
 \bigoplus_{j=1}^r \Pi ^{(j)} (\GL_{n}^{(\rm res)}(\wt F_{\o\omega _j}\times_D\wt F_{\omega _j})) \ , \; 1\le r\le\infty \ ,
 \]
 where:
\Bi
\item $\wt F_{\omega _\oplus}=\sum\limits_j \wt F_{\omega _j}$;

\item $\GL_{n}^{(\rm res)}(\wt F_{\o\omega }\times_D\wt F_{\omega })$ is a bilinear ``diagonal'' algebraic semigroup.
\Ei
\Ee
\vskip 11pt

\subsection{Proposition (Langlands global correspondence on $\GL_{n}(\wt F_{\o\omega }\times_D\wt F_{\omega })$)}

{\em 
The Langlands global correspondence on the complex (diagonal) bilinear algebraic semigroup
$\GL_{n}^{(\rm res)}(\wt F_{\o\omega }\times_D\wt F_{\omega })$
is given by the isomorphism:
\[
\LGC_{\ptcit}: \qquad
\sigma ^{(\rm res)}_{\wt \omega \RL}(W^{ab}_{\dot{\wt F}_{\o \omega }}\times_D
W^{ab}_{\dot{\wt F}_{\omega }})
\quad\To\quad \Pi  (\GL_{n}^{(\rm res)}(\wt F_{\o\omega }\times_D\wt F_{\omega }))
\] between the set $\sigma ^{(\rm res)}_{\wt \omega \RL}(W^{ab}_{\dot{\wt F}_{\o \omega }}\times_D
W^{ab}_{\dot{\wt F}_{\omega }})$ of the $n$-dimensional complex conjugacy class representatives of the diagonal products, right by left, of global Weil subgroups given by the diagonal algebraic bilinear semigroup 
$G_{(\rm res)}^{(2n)}( \wt F_{\o \omega }\times_D\wt F_{\omega }) $ and its cuspidal representation given by $\Pi  (\GL_{n}^{(\rm res)}(\wt F_{\o\omega }\times_D\wt F_{\omega }))$ in such a way that
{\bbf
\[
\Pi  (\GL_{n}^{(\rm res)}(\wt F_{\o\omega _\oplus}\times_D\wt F_{\omega _\oplus}))
= \sum_j \L(
\lambda ^\half(2n,j)\ e^{-2\pi ijz}\times_D
\lambda ^\half(2n,j)\ e^{2\pi ijz} \R)
\]
be the cuspidal biform of 
$\GL_{n}^{(\rm res)}(\wt F_{\o\omega }\times_D\wt F_{\omega })$\/}.
}
\vskip 11pt

\bpr
From section 1.1.4, it results that:
\[
\sigma ^{(\rm res)}_{\wt \omega \RL}(W^{ab}_{\dot{\wt F}_{\o \omega }}\times_D
W^{ab}_{\dot{\wt F}_{\omega }})
=\GL_{n}^{(\rm res)}(\wt F_{\o\omega }\times_D\wt F_{\omega })\;,\]
where $\sigma ^{(\rm res)}_{\wt \omega \RL}=\sigma ^{(\rm res)}_{\omega _R}\times_D
\sigma ^{(\rm res)}_{\omega _L}$, for the restricted case introduced in section~1.1.5.

According to section~1.1.6, the $j$-th cuspidal representation of the $j$-th conjugacy class of
$G_{(\rm res)}^{(2n)}( \wt F_{\o \omega }\times_D\wt F_{\omega }) $ is given by $
\Pi  _j(G_{(\rm res)}^{(2n)}( \wt F_{\o \omega_j }\times_D\wt F_{\omega_j}))$ (or by
$\Pi  _j(\GL^{(\rm res)}_{n}( \wt F_{\o \omega_j }\times_D\wt F_{\omega_j}))$).

So we get the commutative diagram:
\[\begin{psmatrix}[colsep=1cm,rowsep=.6cm]
\sigma ^{(\rm res)}_{\wt \omega \RL}(W^{ab}_{\dot{\wt F}_{\o \omega }}\times_D
W^{ab}_{\dot{\wt F}_{\omega }})
& & \Pi  (\GL^{(\rm res)}_{n}( \wt F_{\o \omega }\times_D\wt F_{\omega}))\\[15pt]
\{G_{(\rm res)}^{(2n)}( \wt F_{\o \omega_j }\times_D\wt F_{\omega_j }) \}_{j=1}^r
& &
\psset{arrows=->,nodesep=5pt}
\everypsbox{\scriptstyle}
\ncline{1,1}{1,3}^{\sim}
\ncline{2,1}{1,3}>{\hspace{-15mm}\Bsm \mbox{}\\[15pt] \text{toroidal isomorphisms}\\
\{ \gamma^T _{F_{\o\omega _j}}\times
 \gamma^T _{F_{\omega _j}}\}_j\Esm}
\psset{doubleline=true,arrows=-}
\ncline{1,1}{2,1}
\end{psmatrix}
\]
where $\Pi  (\GL^{(\rm res)}_{n}( \wt F_{\o \omega }\times_D\wt F_{\omega}))=
\{ \Pi  _j(G^{(2n)}( \wt F_{\o \omega_j }\times_D\wt F_{\omega _j} )\}_j=
\{ \phi _R(T^{2n}_L[j])\otimes_D \phi _L(T^{2n}_L[j])\}_j$
as resulting from section~1.6.\epr
\vskip 11pt

\subsection{Cuspidal representation of real algebraic bilinear semigroups}

\Bean
\item {\bbf
A real cuspidal representation, covering the complex cuspidal representation
$\Pi  (\GL^{(\rm res)}_{n}( \wt F_{\o \omega }\times_D\wt F_{\omega}))$\/}, can be obtained for the real diagonal bilinear algebraic semigroup
$G^{(2n)}( \wt F_{\o v }\times_D\wt F_{v_j})$
{\bbf by summing the cuspical subrepresentations of its conjugacy class representatives taking into account the inclusion\/} (which is also a covering)
\[\begin{psmatrix}[colsep=1cm,rowsep=.6cm]
G^{(2n)}( \wt F_{\o v }\times_D\wt F_{v})
&&
G^{(2n)}_{(\rm res)}( \wt F_{\o \omega }\times_D\wt F_{\omega })
\psset{arrows=->,nodesep=5pt}
\ncline[hookwidth=-1mm]{H->}{1,1}{1,3}
\end{psmatrix}
\]
mentioned in section~1.1.5.

Let then:
\[
\gamma ^T_{F_{v _{j,m_j}}}:\quad \wt F_{v _{j,m_j}}\To F^T_{v_{j,m_j}}\qquad
\rresp{\gamma ^T_{F_{\o v _{j,m_j}}}:\quad \wt F_{\o v _{j,m_j}}\To F^T_{\o v_{j,m_j}}}
\]
be the toroidal isomorphism mapping each \lr extension
$\wt F_{v _{j,m_j}}$ \resp{$\wt F_{\o v _{j,m_j}}$} into its toroidal compact equivalent
$F^T_{v_{j,m_j}}$ \resp{$F^T_{\o v_{j,m_j}}$} which is a semicircle localized in the upper \resp{lower} half space.

The fibre bundle morphism:
\begin{align*}
T_{2n}(F^T_{v_{j,m_j}}): \qquad F^T_{v_{j,m_j}} \quad &\To \quad T^{(2n)}(F^T_{v_{j,m_j}})\\
\rresp{T_{2n}(F^T_{\o v_{j,m_j}}): \qquad F^T_{\o v_{j,m_j}} \quad &\To \quad T^{(2n)}(F^T_{\o v_{j,m_j}})}
\end{align*}
sends $F^T_{v_{j,m_j}}$ \resp{$F^T_{\o v_{j,m_j}}$} into the $2n$-dimensional real semitorus
$T^{(2n)}(F^T_{v_{j,m_j}})$ \resp{$T^{(2n)}(F^T_{\o v_{j,m_j}})$}.

\item Every \lr function on the conjugacy class representative
$G^{(2n)}(  F^T_{v_{j,m_j}})$
\resp{$G^{(2n)}(  F^T_{\o v_{j,m_j}})$} is a function \resp{cofunction}
$\phi _L(T^{2n}_L[j,mj])$
\resp{\linebreak $\phi _R(T^{2n}_R[j,m_j])$} on the real semitorus
$T^{2n}_L[j,m_j]$
\resp{$T^{2n}_R[j,m_j]$} having the analytic development:
\begin{align*}
\phi _L(T^{2n}_L[j,m_j]) & = \lambda ^\half(2n,j,m_j)\ e^{2\pi ijz}\;, \quad x\in \rit^{2n}\;,\\
\rresp{\phi _R(T^{2n}_R[j,m_j]) & = \lambda ^\half(2n,j,m_j)\ e^{-2\pi ijz}}
\end{align*}
where $\lambda (2n,j,m_j)$ is the product of the eigenvalues of the $(j,m_j)$-th coset representative of the product, right by left, of Hecke operators.

This function \resp{cofunction} $\phi _L(T^{2n}_L[j,m_j]) $
\resp{$\phi _R(T^{2n}_R[j,m_j]) $} is the cuspidal representation
$\Pi  ^{(j,m_j)}(G^{(2n)}(\wt F_{v_{j,m_j}}))$
\resp{$\Pi  ^{(j,m_j)}(G^{(2n)}(\wt F_{\o v_{j,m_j}}))$} of the $(j,m_j)$-th conjugacy class  representative of 
$G^{(2n)}(\wt F_{v})$
\resp{$G^{(2n)}(\wt F_{\o v})$} because\linebreak 
$\Pi  (\GL_{2n}(\wt F_{\o v_{\otimes}}\times_D\wt F_{v_{\otimes}}))$
is given by the Fourier biseries
\[
\Pi  (\GL_{2n}(\wt F_{\o v_{\otimes}}\times_D\wt F_{v_{\otimes}}))= \bigoplus_{j,m_j}
\Pi ^{(j,m_j)} (\GL_{2n}(\wt F_{\o v_{j,m_j}}\times_D\wt F_{v_{j,m_j}}))\;,\]
where 
$\wt F_{v_{\oplus}}=\sum\limits_{j,m_j} \wt F_{v_{j,m_j}}$, and corresponds to a cuspidal biform.
\Ee
\vskip 11pt

\subsection{Proposition (Langlands global correspondence on $\GL_{2n}(\wt F_{\o v}\times_D\wt F_{v})$)}

{\em
The Langlands global correspondence on the real diagonal bilinear algebraic semigroup
$\GL_{2n}(\wt F_{\o v}\times_D\wt F_{v})$ is given by the isomorphism:
\[
\LGC_\rit: \qquad 
\sigma _{\wt v \RL}(W^{ab}_{\dot{\wt F}_{\o v }}\times_D
W^{ab}_{\dot{\wt F}_{v }})
\quad\To\quad \Pi  (\GL_{2n}(\wt F_{\o v }\times_D\wt F_{v}))
\]
between the set $\sigma _{\wt v \RL}(W^{ab}_{\dot{\wt F}_{\o v }}\times_D
W^{ab}_{\dot{\wt F}_{v }})$ of the $2n$-dimensional real conjugacy class representatives of the diagonal products, right by left, of global Weil subgroups given by the algebraic bilinear semigroup $G^{(2n)}(\wt F_{\o v }\times_D\wt F_{v})$ and its cuspidal representation
$\Pi  (\GL_{2n}(\wt F_{\o v}\times_D\wt F_{v}))$ in such  a way that $\Pi  (\GL_{2n}(\wt F_{\o v_\oplus}\times_D\wt F_{v_\oplus}))$ be a  ``cuspidal biform'' on
$\GL_{2n}(\wt F_{\o v}\times_D\wt F_{v})$.
}
\vskip 11pt

\bpr
As in proposition~1.7, the proposition results from the commutative diagram:
\be\begin{psmatrix}[colsep=1cm,rowsep=.6cm]
\sigma _{\wt v \RL}(W^{ab}_{\dot{\wt F}_{\o v }}\times_D
W^{ab}_{\dot{\wt F}_{v}})
& & \Pi  (\GL_{2n}( \wt F_{\o v }\times_D\wt F_{v}))\\[15pt]
\{ G^{(2n)}( \wt F_{\o v_{j,m_j} }\times_D\wt F_{v_{j,m_j} }) \}_{j,m_j}
& &
\psset{arrows=->,nodesep=5pt}
\everypsbox{\scriptstyle}
\ncline{1,1}{1,3}^{\LGC_\rit}_{\sim}
\ncline{2,1}{1,3}>{\hspace{-15mm}\Bsm \mbox{}\\[15pt] \text{toroidal isomorphisms}\\
\{ \gamma^T _{F_{\o v_{j,m_j}}}\times
 \gamma^T _{F_{v _{j,m_j}}}\}_{j,m_j}\Esm}
\psset{doubleline=true,arrows=-}
\ncline{1,1}{2,1}
\end{psmatrix}\tag*{\eop}
\ee
\vskip 11pt


\subsection{Corollary (Inclusion of the Langlands real global correspondence into the complex global correspondence)}

{\em
The commutative diagram:
\[\begin{psmatrix}[colsep=1cm,rowsep=.6cm]
\sigma _{\wt \omega \RL}(W^{ab}_{\dot{\wt F}_{\o \omega }}\times_D
W^{ab}_{\dot{\wt F}_{\omega }})
& & \Pi  (\GL^{(\rm res)}_{n}( \wt F_{\o \omega }\times_D\wt F_{\omega}))\\[15pt]
\sigma _{\wt v \RL}(W^{ab}_{\dot{\wt F}_{\o v }}\times_D
W^{ab}_{\dot{\wt F}_{v }})
& & \Pi  (\GL_{2n}( \wt F_{\o v }\times_D\wt F_{v}))
\psset{arrows=->,nodesep=5pt}
\everypsbox{\scriptstyle}
\ncline{1,1}{1,3}^{\LGC_\cit}
\ncline{2,1}{2,3}^{\LGC_\rit}
\ncline[hookwidth=-1mm]{H->}{2,1}{1,1}
\ncline[hookwidth=-1mm]{H->}{2,3}{1,3}
\end{psmatrix}
\]
implies that the cuspidal biform
$\Pi  (\GL^{(\rm res)}_{n}( \wt F_{\o \omega_\oplus }\times_D\wt F_{\omega_\oplus}))$
is covered by the product, right by left, 
$\Pi  (\GL_{2n}( \wt F_{\o v_\oplus }\times_D\wt F_{v_\oplus}))$ of Fourier series over real archimedean completions.
}
\vskip 11pt

\bpr This is a consequence of the propositions~1.7 and 1.9.\epr
\vskip 11pt


\subsection{Bilinear context for the Langlands functoriality conjecture}

The functoriality conjecture introduced by R. Langlands deals with the product of cuspidal representations of algebraic linear groups over adele rings.  Transposed in this bilinear context, this problem is easily solvable by taking into account the cross binary operation of bilinear (algebraic) semigroups introduced in \cite{Pie3}.  Indeed, the Langlands functoriality conjecture then results from the reducibility of representations of bilinear (algebraic) semigroups \cite{Pie4}, covering their linear equivalents.

This new bilinear approach useful in the decomposition of the  bilinear cohomology of the bilinear (algebraic) semigroups $ \GL_{2n}( \wt F_{\o v }\times_D\wt F_{v})$ can be stated as follows:
\vskip 11pt


\subsection{Proposition (Reducibility of representations of bilinear (algebraic) semigroups)}

{\em
{\bbf 
The cuspidal (and holomorphic) representation 
$\Pi (\GL_{2n}( \wt F_{\o v }\times_D\wt F_{v}))$ of the (bilinear (algebraic) semigroup
$\GL_{2n}( \wt F_{\o v }\times_D\wt F_{v})$ is (non orthogonally) completely reducible if it decomposes:\/}
\Bean
\item diagonally according to the direct sum
\[ \bigoplus_{\ell=1}^n \Pi  ^{(2_\ell)} (\GL_{2_\ell}( \wt F_{\o v }\times_D\wt F_{v}))\]
of irreducible cuspidal (and holomorphic) representations of the (algebraic) bilinear semigroups
$\GL_{2_\ell}( \wt F_{\o v }\times_D\wt F_{v})$;

\item and off-diagonally according to the direct sum
\[ \bigoplus_{k\neq\ell=1}^r \L(\Pi  ^{(2_k)} (\GL_{2_k}( \wt F_{\o v }))\otimes
\Pi  ^{(2_\ell)} (\GL_{2_\ell}( \wt F_{ v }))\R)\]
of the products of irreducible cuspidal (and holomorphic) representations of cross linear (algebraic)  semigroups
$\GL_{2_k}( \wt F_{\o v })\times\GL_{2_\ell}( \wt F_{ v } )\equiv 
T^t_{2_k}( \wt F_{\o v })\times T_{2_\ell }( \wt F_{ v })$.
\Ee
}
\vskip 11pt

\bpr The thesis directly results from the definition of a bilinear semigroup introduced in \cite{Pie3} and was developed in \cite{Pie3}.\epr

\section{Lower bilinear  $K$-theory based on homotopy semigroups viewed as deformations of Galois representations}

\subsection{Main tool of the global program of Langlands}

It results from chapter~1 that the Langlands global program refers mainly to the (functional) representation space
$\fREPSP(\GL_n(F_{\o\omega }\times F_\omega ))\equiv G^{(2n)}(F_{\o\omega }\times F_\omega )$ of the complex complete (algebraic) bilinear semigroup
$\GL_n(F_{\o\omega }\times F_\omega )$, covered by its real equivalent
$\fREPSP(\GL_{2n}(F_{\o v }\times F_v ))\equiv G^{(2n)}(F_{\o v}\times F_v )$, because these (bisemi)spaces are representations of the products, right by left, of global Weil semigroups.
\vskip 11pt

\subsection{General bilinear cohomology}

Related to the reducibility of representations of bilinear (algebraic) semigroups (which are abstract bisemivarieties), recalled in proposition~1.12, {\bbf a general bilinear cohomology theory was defined\/} in section~3.2 of \cite{Pie2} {\bbf as a contravariant bifunctor\/}:
\begin{multline*}
\HH^{2*}:\L\{
\text{smooth abstract (algebraic) bisemivarieties}\R.\\
\L. G^{(2n)}(F_{\o\omega }\times F_\omega )=\fREPSP ( \GL_n ( F_{\o\omega }\times F_\omega ))\R\}\hspace{3cm}\\
\To \L\{ \text{graded (functional) representation spaces of the}\R.\\
\hspace{3cm} \text{complete (algebraic) bilinear semigroups }
\L.\GL_* ( F_{\o\omega }\times F_\omega )\R\}
\end{multline*}
written in the conventional form:
\begin{multline*}
H^{2*}(G^{(2n)} ( F_{\o\omega }\times F_\omega ),\fREPSP ( \GL_* ( F_{\o\omega }\times F_\omega )))\\
=\bigoplus_i
H^{2i} ( G^{(2n)}(F_{\o\omega }\times F_\omega ),\fREPSP ( \GL_i (F _{\o\omega }\times F_\omega )))\;.\end{multline*}
Taking into account the inclusion
\[\begin{psmatrix}[colsep=1cm,rowsep=.6cm]
G^{(2n)} ( F_{\o v }\times F_v )
&&
G^{(2n)}_{(\rm res)} ( F_{\o\omega }\times F_\omega )
\psset{arrows=->,nodesep=5pt}
\ncline[hookwidth=-1mm]{H->}{1,1}{1,3}
\end{psmatrix}
\]
of the real (algebraic) bilinear semigroup
$G^{(2n)} ( F_{\o v }\times F_v )$, which is an abstract real bisemivariety, into the corresponding complex (algebraic) abstract bisemivariety
$G^{(2n)}_{(\rm res)} ( F_{\o\omega }\times F_\omega )$ as developed in section~1.5,
{\bbf the general bilinear cohomology can be rewritten in function of rational (bi)coefficients (algebraic case) or in function of real (bi)coefficients (abstract (complete) general case)\/}:
\begin{multline*}
H^{2*}(G^{(2n)} ( F_{\o\omega }\times F_\omega ),\FREPSP ( \GL_{2*} ( F_{\o v }\times F_v )))\\
=\bigoplus_{\Bsm i\\ i\le n\Esm}
H^{2i} ( G^{(2n)}(F_{\o\omega }\times F_\omega ),\FREPSP ( \GL_{2i} (F _{\o v }\times F_v )))\;.\end{multline*}
This corresponds to Hodge bisemicycles \cite{D-M-O-S}, \cite{Riv}, sending the abstract (and, thus, also the algebraic) complex bisemivariety
$G^{(2n)}(F_{\o\omega }\times F_\omega )$ into the abstract (and, thus, also algebraic) real bisemivarieties
$G^{(2i)} (F _{\o v }\times F_v )\equiv
\FREPSP ( \GL_{2i} (F _{\o v }\times F_v ))$ \cite{Pie2} in such a way that there is a bifiltration
$F^p\RL$ on the right and left cohomology semigroups of
$H^{2i} ( G^{(2n)}(\cdot \times \cdot  ),-)$ given by
\begin{multline*}
\qquad F^p\RL H^{2i}  ( G^{(2n)}(F_{\o\omega }\times F_\omega ), G^{(2i)}(F_{\o v }\times F_v ))\\
=\bigoplus_{i=p+q}
 H^{2(p+q)}  ( G^{(2n)}(F_{\o\omega }\times F_\omega ), G^{2(p+q)}(F_{\o v }\times F_v ))\;.\qquad \end{multline*}
 \vskip 11pt


\subsection{Main properties of the general bilinear cohomology}

In addition to the bifiltration on Hodge bisemicycles, the general bilinear cohomology is characterized by the following properties.

\Bean
\item {\bbf A bisemicycle map\/} \cite{Mur}, \cite{Mor}:
\[
\gamma ^i_{G^{(2n)}_{\o\omega \times\omega }}: \qquad 
\Zs^i ( G^{(2n)}(F_{\o\omega }\times F_\omega ))\quad \To\quad 
H^{2i} ( G^{(2n)}(F_{\o\omega }\times F_\omega ), G^{(2i)}(F_{\o v }\times F_v ))
\]
\Bi
\item from the bilinear semigroup 
$\Zs^i ( G^{(2n)}(F_{\o\omega }\times F_\omega ))$ of compactified \resp{noncompactified} bisemicycles of codimension $i$, in the abstract \resp{algebraic} case, on the bilinear complete \resp{algebraic} semigroup
$G^{(2n)}(F_{\o\omega }\times F_\omega )$
\resp{$G^{(2n)}(\wt F_{\o\omega }\times \wt F_\omega )$}

\item into the bilinear cohomology
$H^{2i} ( G^{(2n)}(F_{\o\omega }\times F_\omega ), G^{(2i)}(F_{\o v }\times F_v ))$ in such a way the the embedding
\[\begin{psmatrix}[colsep=1cm,rowsep=.6cm]
 G^{(2i)}(F_{\o v }\times F_v )
&&
G^{(2i)}(F_{\o \omega  }\times F_\omega  ) 
\psset{arrows=->,nodesep=5pt}
\ncline[hookwidth=-1mm]{H->}{1,1}{1,3}
\end{psmatrix}
\]
 of the real bisemivariety
$ G^{(2i)}(F_{\o v }\times F_v )$ into its complex equivalent
$ G^{(2i)}(F_{\o \omega  }\times F_\omega  )$ is directly related to the Hodge bisemicycles according to section~2.2.
\Ei

\item {\bbf A K\"unneth isomorphism\/}:
\begin{multline*}
H^{2i} ( G^{(2n)} ( F_{\o\omega }), G^{(2i)}(F_{\o v } ) )
\otimes_{F_{\o v}\times F_v}
H^{2i} ( G^{(2n)} ( F_{\omega }), G^{(2i)}(F_{v } ) )\\
\To \qquad 
H^{2i} ( G^{(2n)}(F_{\o\omega }\times F_\omega ), G^{(2i)}(F_{\o v }\times F_v ))\;, \quad n\ge i\;,
\end{multline*}
associated with the existence of an abstract real bisemivariety
$G^{(2i)}(F_{\o v }\times F_v )$ of real dimension $2i$, covering its complex equivalent
$G^{(2i)}(F_{\o \omega  }\times F_\omega )$, in the (tensor) product of a right complex abstract semivariety
$G^{(2n)} ( F_{\o \omega } )$ of complex dimension $n$ by its left equivalent
$G^{(2n)} (F_{ \omega } )$.
\Ee
\vskip 11pt

\subsection{New algebraic interpretation of homotopy}

The next step consists in finding the $K$-theory associated with the general bilinear cohomology and in defining the corresponding Chern classes.  But, as $K$-theories are related to homotopy and as the proposed bilinear cohomology  is essentially a motivic (bilinear) cohomology theory or a Weil (bilinear) cohomology theory \cite{Pie2}, {\bbf the homotopy must be proved to result from algebraic geometry in order that this general context be coherent\/}.

It will then be proved that the concept of homotopy in topology corresponds to a deformation of Galois representation as introduced in \cite{Pie7} and briefly recalled now.
\vskip 11pt

\subsection{Deformations of Galois representations \cite{Maz}}

{\bbf Two kinds of deformations of $(2)n$-dimensional representations of global Weil (or Galois) (semi)groups given by bilinear (algebraic) semigroups over complete global Noetherian bisemirings were envisaged \/} \cite{Pie7}, \cite{Pie1}.

\Bean
\item {\bbf global bilinear quantum deformations\/} leaving invariant the orders of inertia subgroups;

\item {\bbf global bilinear deformations\/} inducing the invariance of their bilinear residue (i.e. pseudounramified) semifields.
\Ee

Case a) will be only taken into account in this paper because the inertia subgroups, being the subgroups of automorphisms of algebraic space quanta, are supposed to be stable.
\vskip 11pt

\subsection{Uniform quantum homomorphism between global coefficient semirings}

Then, a global quantum deformation results from a global coefficient semiring quantum homomorphism.

A \lr global (compactified) coefficient semiring $F_v$ \resp{$F_{\o v}$} is given by the set of infinite pseudoramified archimedean embedded completions:
\[ F_{v_1}\subset \dots \subset F_{v_{j,m_j}} \subset \dots\subset F_{v_{r,m_r}}\qquad 
\rresp{F_{\o v_1}\subset \dots \subset F_{\o v_{j,m_j}} \subset \dots\subset F_{\o v_{r,m_r}}},
\]
as developed in section~1.1, where two neighbouring completions $F_{v_j}$ and $F_{v_{j+1}}$ differ by a transcendental quantum $F_{v^1_j}$ characterized by a transcendence degree $\Tr\cdot d\cdot F_{v^1_j} \big/ k=N$.

Let $F_{v+\ell}$
\resp{$F_{\o{v+\ell}}$} denote another \lr global coefficient semiring of which completions are those of $F_v$ \resp{$F_{\o v}$} increased by ``$\ell$'' transcendental quanta:
\[ F_{v_{1+\ell}}\subset \dots \subset F_{v_{j+\ell}} \subset \dots\subset F_{v_{r+\ell}}\qquad 
\rresp{F_{\o v_{1+\ell}}\subset \dots \subset F_{\o v_{j+\ell}} \subset \dots\subset F_{\o v_{r+\ell}}}.
\]

{\bbf A uniform quantum homomorphism between global coefficient semirings is given by\/}:
\[
Qh_{F_{v+\ell}\to F_v}:\quad F_{v+\ell}\To F_v\qquad 
\rresp{Qh_{F_{\o{v+\ell}}\to F_{\o v}}:\quad F_{\o{v+\ell}}\To F_{\o v}}
\]
in such a way that:
\Bena
\item the kernel $K(Qh_{F_{v+\ell}\to F_v} )$ \resp{$K(Qh_{F_{\o{v+\ell}}\to F_{\o v}})$} of the quantum homomorphism
 $Qh_{F_{v+\ell}\to F_v}$ \resp{$Qh_{F_{\o{v+\ell}}\to F_{\o v}}$}, inducing an isomorphism on their global inertia subgroups, is characterized by a transcendence degree 
\[
\Tr\cdot d\cdot F_{v+\ell}\big/ k-\Tr\cdot d\cdot F_{v}\big/ k=N\times \ell\times \sum\limits_j m_{1+j}\] (if $m_{1+j+\ell}=m_{1+j}$).
 
 \item {\bbf this quantum homomorphism corresponds to a base change \/} from
 $F_v$ \resp{$F_{\o v}$} into $F_{v+\ell}$ \resp{$F_{\o{v+\ell}}$} of which transcendence extensions degree is
\[
\Tr\cdot d\cdot F_{v+\ell}\big/ k-\Tr\cdot d\cdot F_v\big/ k\]
which means an increment of $\ell$ quanta on each completion of the coefficient semiring
$F_v$ \resp{$F_{\o v}$};
\Ee
\vskip 11pt

\subsection{Quantum deformations of Galois representations over global bisemirings}

{\bbf A global bilinear quantum deformation representative, resulting from a global bilinear coefficient semiring quantum homomorphism\/}
\[
Qh_{F_{\o{v+\ell}}\times F_{v+\ell}\to F_{\o v}\times F_v}:
\quad F_{\o{v+\ell}}\times F_{v+\ell} \quad \To\quad F_{\o v}\times F_v\;,
\]
is an equivalence class representative $\rho _{F_\ell}$ of lifting
\[
\begin{psmatrix}[colsep=1cm,rowsep=.6cm]
\Gal ( \dot{\wt F}_{\o{v+\ell}}/k ) \times \Gal ( \dot{\wt F}_{v+\ell}/k )
&&
\Gal ( \dot{\wt F}_{\o{v}}/k ) \times \Gal ( \dot{\wt F}_{v}/k )\\[15pt]
\GL_n ( {F}_{\o{v+\ell}} \times F_{v+\ell} )&&
\GL_n ( {F}_{\o{v}} \times F_{v} )
\psset{arrows=->,nodesep=5pt}
\everypsbox{\scriptstyle}
\ncline{1,1}{1,3}^{Qh_{F_\ell\to F}}
\ncline{2,1}{2,3}^{Qh_{G_\ell\to G}}
\ncline{1,1}{2,1}>{\rho _{F_\ell}}
\ncline{1,3}{2,3}>{\rho _F}
\end{psmatrix}
\]
with the notations of sections~1.1 and 1.3.
\vskip 11pt

{\bbf A $n$-dimensional global bilinear quantum deformation of $\rho _F$ is an equivalence class of liftings
$\{\rho _{F_\ell}\}_\ell$, $1\le\ell\le\infty $\/}, described by the following diagram:
\[
\begin{psmatrix}[colsep=1cm,rowsep=.6cm]
1 & 
 \Gal ( \delta \dot{\wt F}_{\o{v+\ell}}/k ) \qquad  &
\Gal ( \dot{\wt F}_{\o{v+\ell}}/k)\qquad  & 
\Gal ( \dot{\wt F}_{\o{v}}/k)\qquad  & 1 \\[-6pt]
 & 
\quad \times \Gal (\delta \dot{\wt F}_{{v+\ell}}/k) &
\quad \times \Gal ( \dot{\wt F}_{v+\ell}/k) & 
\quad \times \Gal ( \dot{\wt F}_{v}/k) &  \\[15pt]
1 &\GL_n (\delta { F}_{\o{v+\ell}}\times\delta {F}_{v+\ell})&
\GL_n ( { F}_{\o{v+\ell}}\times {F}_{v+\ell})&
\GL_n ( { F}_{\o{v}}\times {F}_{v})&1
\psset{arrows=->,nodesep=5pt}
\everypsbox{\scriptstyle}
\ncline{1,1}{1,2}
\ncline{1,2}{1,3}
\ncline{1,3}{1,4}
\ncline{1,4}{1,5}
\ncline{3,1}{3,2}
\ncline{3,2}{3,3}
\ncline{3,3}{3,4}
\ncline{3,4}{3,5}
\ncline{2,2}{3,2}>{\delta \rho _{F_\ell}}
\ncline{2,3}{3,3}>{\rho _{F_\ell}}
\ncline{2,4}{3,4}>{\rho _{F}}
\end{psmatrix}
\]
of which ``Weil kernel'' is 
$\Gal(\delta \dot{\wt F}_{\o{v+\ell}}/k)  \times \Gal(\delta \dot{\wt F}_{{v+\ell}}/k)$
and ``$\GL_n(\cdot\times\cdot)$'' kernel is\linebreak
$\GL_n(\delta { F}_{\o{v+\ell}}\times\delta {F}_{{v+\ell}})$.

This equivalence class of liftings $\{\rho _{F_\ell}\}_\ell$ is then given by
\[ \rho _{F_\ell}=\rho _F+\delta \rho _{F_\ell}\;, \qquad \forall\ \ell\;, \quad 1\le\ell\le\infty \;,\]
in such a way that two liftings 
$\rho _{F_{\ell_1}}$ and $\rho _{F_{\ell_2}}$ are strictly equivalent if they can be transformed one into another by conjugation by bielements of
$\GL_n( { F}_{\o{v+\ell}}\times {F}_{{v+\ell}})$ in the kernel of
$Qh_{G_\ell\to\ G}$.
\vskip 11pt

\subsection{Proposition}

{\em
The transformation of kernels
\[
\GL_n (\delta { F}_{\o{v+\ell_1}}\times\delta {F}_{v+\ell_1})\quad \To \quad
\GL_n (\delta { F}_{\o{v+\ell_2}}\times\delta {F}_{v+\ell_2})\]
corresponds to a base change from
$\GL_n( { F}_{\o{v+\ell_1}}\times {F}_{v+\ell_1})$ into
$\GL_n( { F}_{\o{v+\ell_2}}\times {F}_{v+\ell_2})$ of which dimension is given by the difference of ranks
\[
\delta r_{G_n (\ell_2-\ell_1 )}=N^{n^2}(f^{n^2}_{v+\ell_2}-f^{n^2}_{v+\ell_1})\]
where $f_{v+\ell_1}$ is the sum of all global residue degrees corresponding to the conjugacy class representatives of
$\GL_n( { F}_{\o{v+\ell_1}}\times {F}_{v+\ell_1})$.
}
\vskip 11pt

\bpr
Referring to section~2.7, it is clear that the liftings
$\rho _{F_{\ell_1}}$ and $\rho _{F_{\ell_2}}$ are respectively characterized by the kernels 
$\GL_n(\delta { F}_{\o{v+\ell_1}}\times\delta {F}_{v+\ell_1})$ and
$\GL_n(\delta { F}_{\o{v+\ell_2}}\times\delta {F}_{v+\ell_2})$.

The kernel
$\GL_n(\delta { F}_{\o{v+\ell_1}}\times\delta {F}_{v+\ell_1})$ is characterized by a rank
$r_{\delta _{G_{n{\ell_1}}}}=f^{n^2}_{\ell_1}\times N^{n^2}$ and the kernel
$\GL_n(\delta { F}_{\o{v+\ell_2}}\times\delta {F}_{v+\ell_2})$ is characterized by a rank
$r_{\delta _{G_{n{\ell_2}}}}=f^{n^2}_{\ell_2}\times N^{n^2}$.

These ranks
$r_{\delta _{G_{n{\ell_1}}}}$ and
$r_{\delta _{G_{n{\ell_2}}}}$ describe the increase of the algebraic dimensions respectively of all the conjugacy class representatives of
$\GL_n( { F}_{\o{v+\ell_1}}\times {F}_{v+\ell_1})$ and
$\GL_n( { F}_{\o{v+\ell_2}}\times {F}_{v+\ell_2})$.

So, the difference of ranks
$(r_{\delta _{G_{n{\ell_2}}}}-r_{\delta _{G_{n{\ell_1}}}})$
characterizes the difference of liftings
$(\rho _{F_{\ell_2}}-\rho _{F_{\ell_1}})$ and describes the base change from
$\GL_n( { F}_{\o{v+\ell_1}}\times {F}_{v+\ell_1})$ to
$\GL_n( { F}_{\o{v+\ell_2}}\times {F}_{v+\ell_2})$.\epr
\vskip 11pt

\subsection{Galois homotopy}

Let $Qh_{v+\ell\to v} : F_{v+\ell}\to F_v$ be a uniform quantum homomorphism sending the global coefficient semiring $F_v$ into the deformed global coefficient semiring $F_{v+\ell}$ obtained from $F_v$ by adding ``$\ell$'' transcendental quanta according to section~2.6.

Let $fh_\ell:F_{v+\ell}\to G^{(2n)}(F_v)$ be a continuous map from $F_{v+\ell}$ into the real abstract linear (semi)variety $G^{(2n)}(F_v)$ over the set $F_v$ of archimedean completions.

Then, there exists a continuous map
\[
FH: \qquad F_v\times I \quad \To \quad G^{(2n)}(F_v)\;, \qquad I=[0,1]\;, \]
such that $FH(x,0)=fh$ and $FH(x,1)=fh_\ell$, $\forall\ x\in F_v$, $x$ being a point or a big point (i.e. a quantum), where $fh$ is the continuous map: $fh:F_v\to G^{(2n)}(F_v)$.

{\bbf This continuous map $FH$ is thus the homotopy of $fh$, and will be called the Galois homotopy of $fh$\/}.
\vskip 11pt

\subsection{Lemma}

{\em
{\bbf The Galois homotopy of the continuous map $fh:F_v\to G^{(2n)}(F_v)$ results from a quantum homomorphism between global coefficient semirings\/}.
}
\vskip 11pt

\bpr It is sufficient to prove that the homotopy classes for all the functions $fh_\ell$ correspond to the classes of the quantum homomorphism $Qh_{v+\ell\to v}$.

Let 
\[\cor FH:\quad \begin{array}[t]{ccc}
F_v\times I & \To & F_{v+\ell}\\[-9pt] t &\to & \ell \end{array}\;, \qquad t\in[0,1]\;,\]
denote the one-to-one correspondence between the product of the basic coefficient semiring $F_v$ by the unit interval $[0,1]$ and the deformed global coefficient semiring $F_{v+\ell}$ in such a way that to any $t\in[0,1]$ corresponds an integer $\ell$ labelling the number of quanta added to $F_v$.

Then, the homotopy $FH:F_v\times I\to G^{(2n)}(F_v)$, interpreted as a family of continuous maps $fh_{(t)}:F_v\to G^{(2n)}(F_v)$ by the relation
$fh_{(t)}(x)=FH(x,t)$, $0\le t\le 1$, allows to associate its homotopy classes for every value of the parameter ``$t$'' with the places $v+\ell$ of $F_{v+\ell}$ which are the classes of the deformed global coefficient semiring $F_{v+\ell}$.

This homotopy, resulting from a quantum deformation of the global coefficient semiring, will be called a Galois homotopy because the deformed coefficient semiring $F_{v+\ell}$ is homeomorphic to the algebraic coefficient semiring
\[ \dot{\wt F}_{v+\ell}=\L\{ \dot{\wt F}_{v_{1+\ell}}, \dot{\wt F}_{v_{j+\ell}},
 \dot{\wt F}_{v_{r+\ell}}\R\}\]
 given by this set of real pseudoramified extensions according to sections~1.1 and 2.6.
 \epr
 \vskip 11pt

\subsection{Galois cohomotopy}

{\bbf The Galois cohomotopy\/} is the inverse Galois homotopy defined by the relation
\[CFH(x,t)=FH(x,1-t)\]
 and corresponding to the homotopy between $fh_\ell$ and $fh$ in such a way that it {\bbf results from the inverse quantum homomorphism\/}
$Qh^{-1}_{v+\ell\to v}$ between global coefficient semirings.
\vskip 11pt

\subsection{Retract semirings}
Let $F_{v^1}=\{F_{v^1},\dots,F_{v^1_r}\}$ denote the unit subset of $F_v$ composed of one quantum in each completion of $F_v$.

{\bbf The global coefficient semiring $F_{v+\ell}$ is said to be retract\/} if the Galois homotopy $CFH(x,1)=FH(x,0)$ corresponds to the constant homotopy, i.e. if $F_{v+\ell}$ is sent to its subsemiring $F_v$.

{\bbf The global coefficient semiring $F_{v+\ell}$ is said to be strongly retract\/} if $F_{v+\ell}$ is sent to the unit subset $F_{v^1}$ of $F_v$.
\vskip 11pt

\subsection{Fundamental group in terms of deformations of Galois representations}

The equivalence classes of maps between a fixed basic coefficient semiring $F_v$ and the real linear (semi)variety $G^{(2n)}(F_v)$ are the homotopy classes corresponding to the classes of the quantum homomorphism $Q_{v+\ell\to v}$ characterized by the integers ``$\ell$'' which are in one-to-one correspondence with the values of the parameter $t\in [0,1]$ of the homotopy.

Let $fh_\ell:F_{v+\ell}\to G^{(2n)}(F_v)$ and $fh_{\ell,d}:F_{(v+\ell)+d}\to G^{(2n)}(F_v)$ be two maps relative respectively to $fh:F_v\to G^{(2n)}(F_v)$ and $fh_\ell$.

They belong to two difference equivalence classes of maps characterized respectively by the integers $\ell$ and $d$.

Being homotopic is then an equivalence relation compatible with the product of equivalence classes.

That is to say, if $\{fh_\ell\}$ denotes the set of continuous maps from $F_{v+\ell}$ into $G^{(2n)}(F_v)$ with respect to $F_v$, and $\{fh_{\ell,d}\}$ denotes the set of continuous maps from $F_{(v+\ell)+d}$ into $G^{(2n)}(F_v)$ with respect to $F_{v+\ell}$,
$\{fh_\ell\}\times\{fh_{\ell,d}\}$ will correspond to the product of the equivalence classes
$\{fh_\ell\}\times\{fh_{\ell,d}\}$.

Taking into account the existence of
\Bean
\item the Galois cohomotopy of which classes are the inverse equivalence classes of the corresponding Galois homotopy,

\item the null homotopy associated with identity homotopy maps,
\Ee
we see that the set of equivalence classes of Galois homotopy forms a group noted\linebreak 
$\Pi(F_v,G^{(2n)}(F_v))$.

If $L_{v^1}$ is the image in $G^{(2n)}(F_v)$ of $F_{v^1}$, we get the fundamental group\linebreak 
$\Pi_1(G^{(2n)}(F_v),L_{v^1})$ in the big point $L_{v^1}$, which is one quantum or the center of blowup of this one \cite{Pie5}.

Remark that {\bbf the equivalence classes of maps between the coefficient semiring
$F_v$ and the real linear abstract (semi)variety $G^{(2n)}(F_v)$ are the equivalence classes of maps between the set $\Omega (L_{v^1},G^{(2n)}(F_v))$ of  oriented paths, which are the set $F_v$ of archimedean completions\/}, and $G^{(2n)}(F_v)$.

These equivalence classes thus depend on the deformations of the Galois compact representations of these paths corresponding to the increase of these by a(n) (in)finite number of transcendental or algebraic quanta.
\vskip 11pt

\subsection{Homotopy (semi)groups in terms of deformations of Galois representations}

The definition of the fundamental (Galois) homotopy \cite{D-N-F} group 
$\Pi_1 ( G^{(2n)} ( F_v),L_{v^1})$ in terms of deformations of Galois representations of paths (or loops) can be easily generalized to the {\bbf $i$-th homotopy group $\Pi_i ( G^{(2n)} ( F_v),L_{v^1_{(j)}})$\/} of the given (semi)variety $G^{(2n)}(F_v)$ with base point $L_{v^1_{(j)}}$.

The set of homotopy classes of maps
\[
{}_Sfh^i_\ell : \qquad S^i\quad \To \quad G^{(2n)}(F_v)\;, \]
sending the base point $b$ of the $i$-sphere $S^i$ to the base point $L_{v^i_{(j)}}$ of $G^{(2n)}(F_v)$,

are equivalently described by maps
\[
{}_Cfh^i_\ell: \qquad [0,1]^i \quad \To \quad G^{(2n)}(F_v)\]
from the $i$-cube $[0,1]^i$ to $G^{(2n)}(F_v)$ by taking its boundary $\delta [0,1]^i$ to $L_{v^i_{(j)}}$.
\vskip 11pt

\subsection{Proposition}

{\em
{\bbf The (semi)group $\Pi_{2i}(G^{(2n)}(F_v),L_{v^1_{(j)}})$ of homotopy classes of maps 
\begin{align*}
{}_Sfh^{2i}_\ell:\qquad S^{2i}_{(\ell)}\quad &\To \quad G^{(2n)}(F_v)\\
 \text{or} \qquad 
{}_Cfh^{2i}_\ell:\qquad [0,1]^{2i}_{\ell}\quad &\To \quad G^{(2n)}(F_v)\end{align*}
 results from the deformations of the Galois compact representation of the semigroup
$\GL_{2i}(\wt F_v)$\/} of real dimension $2i$ {\bbf given by the kernels $G^{(2i)}(\delta F_{v+\ell})$\/} of the maps:
\[
\GD_\ell ^{2i}: \qquad \begin{array}[t]{ccc}
G^{(2i)}(F_{v+\ell}) & \To & G^{(2i)}(F_v)\;, \\[-9pt]
t^{2i} & \to & \ell^{2i} \end{array} \qquad \forall\ \ell\;, \quad 1\le \ell\le\infty \;, \]
in such a way that the $2i$-th powers of the integers ``$\ell$'' be in one-to-one correspondence with the $2i$-th powers of the values of the parameter $t\in[0,1]$.
}
\vskip 11pt

\bpr
The structure of (semi)group of
$\Pi_{2i}(G^{(2n)}(F_v),L_{v^1_{(j)}})$ results from the composition of its homotopy classes.

Let ${}_Cfh^{2i}_\ell:[0,1]^{2i}_{\ell}\to G^{(2n)}(F_v)$ be the homotopy class of maps
${}_Cfh^{2i}_\ell$ characterized by the value(s) $t^{2i}_\ell\in[0,1]^{2i}$ of the parameter $t^{2i}$ of the $2i$-cube $[0,1]^{2i}$.

Let ${}_Cfh^{2i}_d:[0,1]^{2i}_{d}\to G^{(2n)}(F_v)$ be another homotopy class of maps 
 ${}_Cfh^{2i}_d$ characterized by the value $t^{2i}_d\in[0,1]^{2i}$ of the parameter $t^{2i}$.
 
 Then, the composition (i.e. sum) of these two homotopy classes is given by the homotopy class of maps:
\[
{}_Cfh^{2i}_{\ell+d}:\qquad [0,1]^{2i}_{\ell+d}\quad \To \quad G^{(2n)}(F_v)\]
characterized by the value $t^{2i}_{\ell+d}\in[0,1]^{2i}$ of the parameter $t^{2i}$.

Referring to lemma~2.10 and section~2.7, it appears that the homotopy class of maps
${}_Cfh^{2i}_\ell:[0,1]^{2i}_{\ell}\to G^{(2n)}(F_v)$ is a deformation of the semivariety
$G^{(2i)}(F_v)\subset G^{(2n)}(F_v)$ and corresponds to the deformation of the Galois compact representation of the linear semigroup $\GL_{2i}(\wt F_v)$ given by the kernel
$G^{(2i)}(\delta F_{v+\ell})$ of the map:
\[
\GD^{2i}_\ell : \qquad G^{(2i)}(F_{v+\ell}) \quad \To \quad G^{(2i)}(F_v)\;.\]
This kernel $(G^{(2i)}(F_{v+\ell}))$ is then characterized by the integer $\ell^{2i}$ which
\Bean
\item denotes the number of quanta added to $G^{(2i)}(F_v)$ by the envisaged deformation;
\item is in one-to-one correspondence with the parameter $t^{2i}_\ell$.
\Ee

It is then clear that the composition
${}_Cfh^{2i}_{\ell+d}$ of the two homotopy classes of maps
${}_Cfh^{2i}_{\ell}$ and ${}_Cfh^{2i}_{d}$ results from the kernel
$G^{(2i)}(\delta F_{v+\ell+d})$ of the composition $\GD^{2i}_d\circ \GD ^{2i}_\ell$ of the maps $\GD ^{2i}_\ell$ and $\GD ^{2i}_d$:
\be
\GD^{2i}_d\circ \GD ^{2i}_\ell: \qquad
G^{(2i)}(F_{v+\ell+d})\quad \To \quad G^{(2i)}(F_v)\;.\tag*{\eop}\ee
\vskip 11pt

\subsection{Cohomotopy (semi)groups in terms of inverse deformations}

If the homotopy group
$\Pi_{2i}(G^{(2n)}(F_v),L_{v^1_{(j)}})$ lacks for inverse homotopy classes (and null-homotopy class), it becomes a homotopy semigroup whose dual is the cohomotopy semigroup noted
$\Pi^{2i}(G^{(2n)}(F_v),L_{v^1_{(j)}})$.

Thus, {\bbf the cohomotopy semigroup
$\Pi^{2i}(G^{(2n)}(F_v),L_{v^1_{(j)}})$ is defined by classes resulting from inverse deformations
$(\GD ^{2i}_\ell)^{-1} : G^{(2i)}(F_v)\to G^{(2i)}(F_{v+\ell})$ of the Galois representations of
$\GL_{2i}(\wt F_v)$\/} \cite{Pie7}.
\vskip 11pt

\subsection{Bilinear (co)homotopy in terms of (inverse) deformations}

Let $G^{(2n)}(F_{\o v})$ denote the semivariety dual of $G^{(2n)}(F_v)$ and let $L_{\o v^1_{(j)}}$ be the base point of $G^{(2n)}(F_{\o v})$.

Then, {\bbf the (Galois) bilinear homotopy (semi)group will be given by\linebreak
$\Pi_{2i}(G^{(2n)}(F_{\o v}\times F_v),L_{\o v^1_{(j)}}\times L_{v^1_{(j)}})$\/} in such a way that its classes of (bi)maps:
\[
{}_Cfh^{2i}_\ell\times_{(D)}{}_Cfh^{2i}_\ell: \qquad
[0,1]^{2i}_\ell\times_{(D)}[0,1]^{2i}_\ell \quad \To \quad
G^{(2n)}(F_{\o v}\times F_v)\]
result from the deformations of the Galois (compact) representations of the bisemivariety
$G^{(2i)}(\wt F_{\o v}\times \wt F_v)$ given by the (bi)kernels
$G^{(2i)}(\delta F_{\o v+\ell}\times \delta F_{v+\ell})$ of the (bi)maps
\[
\GD ^{2i}_{\ell\RL}: \qquad
G^{(2i)}(F_{\o v+\ell}\times F_{v+\ell})\quad \To \quad
G^{(2i)}(F_{\o v}\times F_v)\;, \quad \forall\ \ell\;.\]
Similarly, the bilinear (Galois) cohomotopy (semi)group will be given by
$\Pi^{2i}(G^{(2n)}(F_{\o v}\times F_v),L_{\o v^1_{(j)}}\times L_{v^1_{(j)}})$ in such a way that its classes of (bi)maps are inverse of those of the bilinear (Galois) homotopy (semi)group
$\Pi_{2i}(G^{(2n)}(F_{\o v}\times F_v),L_{\o v^1_{(j)}}\times L_{v^1_{(j)}})$.
\vskip 11pt

\subsection{Proposition}

{\em
Taking into account the group homomorphism of  Hurewicz:
\[
hH: \qquad
\Pi_{2i}(G^{(2n)}(F_v),L_{v^1_{(j)}}) \quad \To \quad
H_{2i} ( G^{(2n)}(F_v),\ZZ )\;, \]
we can specialize it to the Galois bilinear homotopy and cohomotopy semigroups according to:
\[
hH\RL : \qquad \Pi_{2i}(G^{(2n)}(F_{\o v}\times F_v),L_{\o v^1_{(j)}}\times L_{v^1_{(j)}}) \quad \To \quad 
H^{2i} ( G^{(2n)}(F_{\o v}\times F_v),\zit\times_{(D)} \zit )\]
\[
hCH\RL : \qquad \Pi^{2i}(G^{(2n)}(F_{\o v}\times F_v),L_{\o v^1_{(j)}}\times L_{v^1_{(j)}}) \quad \To \quad 
H_{2i} ( G^{(2n)}(F_{\o v}\times F_v),\zit\times_{(D)} \zit )\]
where:
\Bena
\item $hH\RL$ is the bilinear semigroup homomorphism from the bilinear homotopy
$\Pi_{2i}(\cdot )$ into the bilinear cohomology $H^{2i}(\cdot )$;

\item $hCH\RL$ is the bilinear semigroup homomorphism from the bilinear cohomotopy
$\Pi^{2i}(\cdot )$ into the bilinear homology $H_{2i}(\cdot )$.
\Ee
}
\vskip 11pt

\bpr
As $\Pi_{2i}(G^{(2n)}(F_{\o v}\times F_v),L_{\o v^1_{(j)}}\times L_{v^1_{(j)}})$
\resp{$\Pi^{2i}(G^{(2n)}(F_{\o v}\times F_v),L_{\o v^1_{(j)}}\times L_{v^1_{(j)}})$}
 is a bilinear homotopy \resp{cohomotopy} semigroup resulting from deformations \resp{inverse deformations} of the Galois (compact) representations of the bilinear semigroup
 $\GL_{2i}(\wt F_{\o v}\times \wt F_v)\subset\GL_{2n}(\wt F_{\o v}\times \wt F_v)$, it is natural to associate to it by the homomorphism $hH\RL$ \resp{$hCH\RL$} the entire bilinear cohomology \resp{homology}
 $H^{2i}(G^{(2n)}(F_{\o v}\times F_v),\zit\times_{(D)} \zit)$
 \resp{$H_{2i}(G^{(2n)}(F_{\o v}\times F_v),\zit\times_{(D)} \zit)$}
 where $\zit\times_{(D)}\zit$ refers to  a bisemilattice deformed by the classes of deformations \resp{inverse deformations} of Galois representations of
$\GL_{2i}(\wt F_{\o v}\times_{(D)} \wt F_v)$ in one-to-one correspondence with the classes of homotopy \resp{cohomotopy}.

Indeed, the cohomology \resp{the homology} is defined with respect to a coboundary \resp{boundary} homomorphism increasing \resp{decreasing} the dimension of one unit.\epr
\vskip 11pt

 \subsection{Topological bilinear $K$-theory}
 
 As the universal cohomology theory is bilinear \cite{Pie2} referring to the Tannakian category \cite{Riv} of representations of affine group schemes, the (topological) $K$-theory of the compact real \resp{complex} bisemivariety
 $G^{(2i)}(F_{\o v}\times F_v)$
 \resp{$G^{(2i)}(F_{\o \omega }\times F_\omega )$} can be naturally introduced and is proved to correspond to the classical definition of the $K$-theory.
 
 Indeed, classically, if $M$ denotes the abelian semigroup (or monoid) of classes of isomorphism of $k$-vector bundles over a compact space $X$, the topological $K$-theory
 $K^0_{\rm top}(X)$ is the symmetrized group of $M$, i.e. the quotient of $M\times M$ by the equivalence relation identifying $(x,y)$ to $(x',y')$, or is the set of cosets of $\Delta (M)$ in $M\times M$ where $\Delta :M\to M\times M$ is a diagonal homomorphism of semigroups \cite{Ati1}.
 
 The locally constant function $r:X\to\nit$ given by $r(x)=\dim E_x$ (where $E$ is the vector bundle over $X$) defines the group homomorphism $K^0_{\rm top}(X) \to H^0(X;\zit)$, where $H^0(X;\zit)$ is the first \v Cech cohomology group of $X$ given by locally constant functions over $X$ with values in $\zit$ \cite{Kar1}.
 
 In this context, let $K^0_{{\rm top}_L}(G^{(2i)}(F_v))$ denote the abelian semigroup (or monoid) of classes of isomorphism of $k$-vector bundles over the compact left semivariety $G^{(2i)}(F_v)$.
 
 Then, {\bbf the topological (bilinear) $K$-theory, noted 
$K^0_{{\rm top}\RL}(G^{(2i)}(F_{\o v}\times F_v))$ or simply
$K^0_{\rm top}(G^{(2i)}(F_{\o v}\times F_v))$, is the set of cosets of
$K^0_{{\rm top}_R} (K^0_{{\rm top}_L}(G^{(2i)}(F_v)))$ where
$K^0_{{\rm top}_R} : K^0_{{\rm top}_L}(G^{(2i)}(F_v))\to
K^0_{{\rm top}\RL}(G^{(2i)}(F_{\o v}\times F_v)))$\/} is the diagonal homomorphism sending the left abelian semigroup
$K^0_{{\rm top}_L}(G^{(2i)}(F_v))$ on the left semivariety
$G^{(2i)}(F_v)$ into the diagonal bilinear semigroup \cite{Pie3}
$K^0_{{\rm top}\RL}(G^{(2i)}(F_{\o v}\times F_v))$ on the product, right by left, of the semivarieties
$G^{(2i)}(F_{\o v})$ and $G^{(2i)}(F_v)$.

Remark that, in the diagonal bilinear semigroup, the cross products are not considered.

Taking into account the derived functors of the $K$-theory of the variety $X$, 
\[ K^{-n}(X)=K^0 ( X\times \rit^n )\]
introduced by Atiyah and Hirzebruch, and the periodicity of the Clifford algebra $C^n$, the group $K^n(X)$ can be introduced from the category of $k$-fibre bundles in graded modules on $C^n$ as a general cohomology theory $n\to K^n(X)$ on the category of pointed compact spaces (or locally compact spaces) \cite{Kar2}.

In this context, let $G^{(2i)}(F_{\o v}\times F_v)\to (F_{\o v}\times F_v)$ be a vector (bi)bundle with (bi)fibre
$\GL_{2i-1}(F_{\o v}\times F_v)$ (or $\rit^{2i-1}$).

Let $K^{2i}(G^{(2n)}(F_{\o v}\times F_v))$ denote the topological (bilinear) $K$-theory of vector (bi)bundles with base 
$G^{(2n)}(F_{\o v}\times F_v)$ and bifibre
$G^{(2n-2i+1)}(F_{\o v}\times F_v)$.

$K^{2i}(G^{(2n)}(F_{\o v}\times F_v))$ is then equivalent to 
$K^{2i}(F_{\o v}\times F_v)$ since the total space of these two vector (bi)bundles is the bisemivariety
$G^{(2i)}(F_{\o v}\times F_v)$.

The cohomology class of obstruction of these fibre bundles is
\[
\alpha _{2i}\in H^{2n-2i+1}(F_{\o v}\times F_v,
\Pi_{2n-2i+1}(G^{(2i)}(F_{\o v}\times F_v)))\]
leading to the Stiefel-Whitney class $W_{2i}=\alpha _{2n-2i+1}\in H^{2i}(F_{\o v}\times F_v,\zit\times_{(D)}\zit )$ of which polynomial is
\[
W(x)=1+W_1\ x+\dots+W_i\ x^{2i}+\dots\]

Similarly, in the complex case, the Chern classes
\[
C_i = \beta _{n-i+1}\in H^{2i}(F_{\o\omega }\times F_\omega ,\zit\times_{(D)}\zit )\]
of the fibre bundles
$G^{(2i)}(F_{\o\omega }\times F_\omega)\to (F_{\o\omega }\times F_\omega)$ with base 
$ (F_{\o\omega }\times F_\omega)$ are included into the Chern polynomial:
\[
C(X)=1+C_1\ x+\dots+C_i\ x^i+\dots\;.\]
\vskip 11pt

\subsection[Chern and Stiefel-Whitney classes in the bilinear $K$-cohomology]{Chern and Stiefel-Whitney classes in the bilinear $K$-\linebreak cohomology}

{\bbf The Stiefel-Whitney character restricted to the class $W_{2i}$ in the bilinear $K$-cohomology\/} of the abstract real bisemivariety
$G^{(2n)}(F_{\o v}\times F_v)$  given by the homomorphism:
\begin{multline*}
SW^{2i}(G^{(2n)}(F_{\o v}\times F_v)): \quad
K^{2i} (G^{(2n)}(F_{\o v}\times F_v)) \\
\quad \To \quad H^{2i}(G^{(2n)}(F_{\o v}\times F_v),G^{(2i)}(F_{\o v}\times F_v))\;, \quad \forall\ i\le n\;, \end{multline*}
and corresponds to {\bbf the Chern character restricted to the class $C_i$ in the bilinear $K$-cohomology\/} of the abstract complex bisemivariety
$G^{(2n)}(F_{\o \omega }\times F_\omega )$  given by the homomorphism
\[
C^i(G^{(2n)}(F_{\o \omega }\times F_\omega )): \qquad
K^{2i} (G^{(2n)}(F_{\o \omega }\times F_\omega ))\quad \To \quad
H^{2i}(G^{(2n)}(F_{\o \omega }\times F_\omega ),G^{(2i)}(F_{\o \omega }\times F_\omega ))\;.\]
\vskip 11pt

\subsection{Proposition}

{\em
Let 
\[
hH\RL: \quad \Pi _{2i}(G^{(2n)}(F_{\o v}\times F_v),L_{\o v^1_{(j)}}\times L_{v^1_{(j)}})
\To 
H^{2i} ( G^{(2n)} (F_{\o v}\times F_v),\zit \times_{(D)}\zit )\]
be the bilinear semigroup homomorphism of Hurewicz from the ``Galois'' bilinear homotopy $\Pi _{2i}(\cdot )$ into the entire bilinear cohomology $H^{2i}(\cdot )$:
{\bbf it can be called (restricted) $\Pi $-cohomology with reference to $K$-cohomology\/}.

Let
\[ C^i(G^{(2n)}(F_{\o v}\times F_v)):\quad
K^{2i}(G^{(2n)}(F_{\o v}\times F_v))
 \To 
H^{2i}(G^{(2n)}(F_{\o v}\times F_v),G^{(2i)}(F_{\o v}\times F_v))\]
be the restricted Chern character in the bilinear $K$-cohomology of the abstract real bisemivariety $G^{(2n)}(F_{\o v}\times F_v)$.

Then, {\bbf the lower bilinear (algebraic) $K$-theory\/} will be given by the equality \resp{homomorphism}
\[
K^{2i}(G^{(2n)}(F_{\o v}\times F_v))
\underset{(\to)}=
\Pi _{2i}(G^{(2n)}(F_{\o v}\times F_v),L_{\o v^1_{(j)}}\times L_{v^1_{(j)}})\]
in such a way {\bbf that the homotopy classes of maps of $\Pi _{2i}(\cdot)$ are \resp{correspond to} liftings of quantum deformations of the Galois representation $\GL_{2i}(\wt F_{\o v}\times \wt F_v)$\/}.
}
\vskip 11pt

\bpr Referring to proposition~1.9, the functional representation space of the product, right by left, of global Weil (semi)groups is given by the real abstract bisemivariety $G^{(2n)}(F_{\o v}\times F_v)$ in the frame of the Langlands global program.

So, the relations between the bilinear cohomology, homotopy and topological $K$-theory of
$G^{(2n)}(F_{\o v}\times F_v)$ are given according to sections~2.2, 2.15 and 2.19 by the commutative diagram
\[
\begin{psmatrix}[colsep=1cm,rowsep=.6cm]
K^{2i}(G^{(2n)}(F_{\o v}\times F_v)) & & \Pi _{2i}(G^{(2n)}(F_{\o v}\times F_v),L_{\o v^1{(j)}}\times L_{v^1_{(j)}})\\[15pt]
& \hspace*{-2cm} H^{2i}(G^{(2n)}(F_{\o v}\times F_v),G^{(2i)}(F_{\o v}\times F_v))\hspace*{-2cm} 
\psset{arrows=->,nodesep=5pt}
\everypsbox{\scriptstyle}
\ncline{1,1}{1,3}^{\Bsm\text{lower (algebraic)}\\ \text{$K$-theory}\Esm\qquad\quad}
\ncline{2,2}{1,1}<{\Bsm \text{inverse restricted}\\ \text{Chern character}\Esm}>{\Bsm \text{inverse restricted}\\ \text{$K$-cohomology}\Esm}
\ncline{1,3}{2,2}>{\text{Hurewicz homomorphism}}
\end{psmatrix}
\]
in such a way that {\bbf the classes of the entire bilinear cohomology
$H^{2i}(G^{(2n)}(F_{\o v}\times F_v),\zit\times_{(D)}\zit)$\/} as well as those of the bilinear $K$-theory $K^{2i}(G^{(2n)}(F_{\o v}\times F_v))$ {\bbf are the homotopy classes of maps of
$\Pi _{2i}(G^{(2n)}(F_{\o v}\times F_v))$ corresponding to the lifts of quantum deformations of the Galois compact representations of
$\GL_{2i}(\wt F_{\o v}\times \wt F_v)$\/} according to proposition~2.15.

It then results that:
\[ 
K^{2i}(G^{(2n)}(F_{\o v}\times F_v))=\Pi _{2i}(G^{(2n)}(F_{\o v}\times F_v),L_{\o v^1_{(j)}}\times L_{v^1_{(j)}})\]
{\bbf defining a lower bilinear (algebraic) $K$-theory\/} with reference of the higher (bilinear) algebraic $K$-theory (of Quillen) reexamined afterwards according to the Langlands global program.\epr
\vskip 11pt

\subsection{Corollary}

{\em
Let
\[ hCH\RL: \quad
\Pi ^{2i}(G^{(2n)}(F_{\o v}\times F_v),L_{\o v^1_{(j)}}\times L_{v^1_{(j)}})
 \To 
H_{2i}(G^{(2n)}(F_{\o v}\times F_v),\zit\times_{(D)}\zit )\]
be the bilinear semigroup homomorphism of Hurewicz from the bilinear ``Galois'' cohomotopy $\Pi ^{2i}(\cdot )$ into the entire bilinear homology $H_{2i}(\cdot )$.

It can be {\bbf called restricted $\Pi $-homology\/} with reference to the $K$-homology.

Let 
\[
C_i(G^{(2n)}(F_{\o v}\times F_v)): \quad
K_{2i}(G^{(2n)}(F_{\o v}\times F_v)) \To 
H_{2i}(G^{(2n)}(F_{\o v}\times F_v),G^{(2i)}(F_{\o v}\times F_v))\]
be the restricted Chern character in the bilinear $K$-homology.

Then, {\bbf the lower bilinear (algebraic) $K$-theory, referring to the cohomotopy\/}, will be given by the equality \resp{homomorphism}:
\[
K_{2i}(G^{(2n)}(F_{\o v}\times F_v))
\underset{(\to)}=
\Pi ^{2i}(G^{(2n)}(F_{\o v}\times F_v),L_{\o v^1_{(j)}}\times L_{v^1_{(j)}})\]
in such a way that the cohomotopy classes of maps
$\Pi ^{2i}(\cdot )$ are \resp{correspond to} inverse liftings of inverse quantum deformations of the Galois (compact) representation
$\GL_{2i}(\wt F_{\o v}\times \wt F_v)$.
}
\vskip 11pt

\bpr
The proof of proposition~2.21, transposed to the lower bilinear algebraic $K$-theory referring to the cohomotopy, is evident here if section~2.16 is taken into account as well as the commutative diagram:
\[
\begin{psmatrix}[colsep=1cm,rowsep=.6cm]
K_{2i}(G^{(2n)}(F_{\o v}\times F_v)) & & \Pi ^{2i}(G^{(2n)}(F_{\o v}\times F_v),L_{\o v^1{(j)}}\times L_{v^1_{(j)}})\\[15pt]
& \hspace*{-2cm} H_{2i}(G^{(2n)}(F_{\o v}\times F_v),G^{(2i)}(F_{\o v}\times F_v))\hspace*{-2cm} 
\psset{arrows=->,nodesep=5pt}
\everypsbox{\scriptstyle}
\ncline{1,1}{1,3}
\ncline{2,2}{1,1}<{\Bsm \text{inverse restricted}\\ \text{$K$-homology}\Esm}
\ncline{1,3}{2,2}>{\text{restricted $\Pi $-homology}}
\end{psmatrix}
\]

\section{Higher bilinear algebraic $K$-theories related to the reducible bilinear global program of Langlands}


\subsection{Prerequisite}

It was noticed in section~2.1 that the main tool of the Langlands global program is the (functional) representation space 
$\fREPSP(\GL_{2n}(F_{\o v}\times F_v))$
\resp{$\fREPSP(\GL_{2n}(F_{\o \omega }\times F_\omega ))$}, of the real \resp{complex} (algebraic) bilinear semigroup
$\GL_{2n}(F_{\o v}\times F_v)$
\resp{$\GL_{n}(F_{\o \omega }\times F_\omega )$}, that is to say a real \resp{complex} abstract bisemivariety
$G^{(2n)}(F_{\o v}\times F_v)$
\resp{$G^{(2n)}(F_{\o \omega }\times F_\omega )$}.

This led us to define a ``lower'' bilinear (algebraic) $K$-theory on the basis of this abstract bisemivariety
$G^{(2n)}(F_{\o v}\times F_v)$
\resp{$G^{(2n)}(F_{\o \omega }\times F_\omega )$}.

In order to introduce a ``higher'' bilinear algebraic $K$-theory referring to the Langlands global program, we have also to take into account the unitary (functional) representation space of the bilinear (algebraic) semigroup
$\GL_{2n}(F_{\o v}\times F_v)$.
\vskip 11pt

\subsection{Parabolic bilinear semigroup}

Let
$P^{(2n)}(F_{\o v^1}\times F_{v^1})$
\resp{$P^{(2n)}(F_{\o \omega^1 }\times F_{\omega^1} )$} be the real \resp{complex} parabolic bilinear semigroup viewed as the smallest bilinear normal pseudoramified subgroup of
$G^{(2n)}(F_{\o v}\times F_v)$
\resp{$G^{(2n)}(F_{\o \omega }\times F_\omega )$}
according to section~1.3 \resp{1.4},

where
$F_{v^1}=\{ F_{v^1_1},\dots,F_{v^1_j},\dots,F_{v^1_r}\}$ is the set of classes of unitary pseudoramified real archimedean completions.

Referring to \cite{Pie2}, the (bisemi)group
$\Int(G^{(2n)}(\wt F_{\o v}\times \wt F_v))$ of Galois inner automorphisms of
$G^{(2n)}(\wt F_{\o v}\times \wt F_v)$ corresponds to the (bisemi)group
$\Aut(P^{(2n)}(\wt F_{\o v^1}\times \wt F_{v^1}))$ of Galois automorphisms of the bilinear parabolic subsemigroup
$P^{(2n)}(\wt F_{\o v^1}\times \wt F_{v^1})$.

It then results that {\bbf $P^{(2n)}(\wt F_{\o v^1}\times \wt F_{v^1})$ can be considered as the unitary representation space of the algebraic bilinear semigroup
$\GL_{2n}(F_{\o v}\times F_v)$\/} because it is the isotropy subgroup of
$G^{(2n)}(F_{\o v}\times F_v)$ fixing its bielements.

On the other hand, as
$G^{(2n)}(\wt F_{\o v}\times \wt F_v)$ is a smooth reductive bilinear affine semigroup, we have that
\[ P^{(2n)}(\wt F_{\o v^1}\times \wt F_{v^1})\approx
(\wt F_{\o v^1})^{2n}\times_{(D)} (\wt F_{v^1})^{2n}\;.\]
\vskip 11pt

\subsection{Lemma}

{\em
The unitary (functional) representation space
$\UfREPSP(\GL_{2n}(F_{\o v}\times F_v))$ of the (algebraic) bilinear semigroup
$\GL_{2n}(F_{\o v}\times F_v)$ is given by
\begin{align*}\UfREPSP ( \GL_{2n} ( F_{\o v}\times F_v ))
&= \fREPSP ( P_{2n}(F_{\o v^1}\times F_{v^1}))\\
&= ( F_{\o v^1})^{2n}\times_{(D)} ( F_{v^1})^{2n}\;.\end{align*}
}
\vskip 11pt

\paragraph{Sketch of proof}: \hspace{5mm} 
As  $G^{(2n)}(\wt F_{\o v}\times \wt F_v)$ is a reductive bilinear semigroup,
$P^{(2n)}(\wt F_{\o v^1}\times \wt F_{v^1})$, being its isotropy subgroup, is (isomorphic to) the product, right by left, of unitary algebraic semitori
$\wt F^{2n}_{\o v^1}$ and
$\wt F^{2n}_{v^1}$.\epr
\vskip 11pt

\subsection{Classical and quantum higher bilinear algebraic $K$-theories}

Two types of equivalent ``higher'' bilinear algebraic $K$-theories on the basis of the global program of Langlands will now be introduced.

\Bena
\item {\bbf The first ``classical'' depends on the geometric dimensions of the classifying bisemispace $\BGL(F_{\o v}\times F_v)$ of
$\GL(F_{\o v}\times F_v)$\/} where
\[
\GL (F_{\o v}\times F_v)=\underrightarrow{\lim} \GL_m (F_{\o v}\times F_v)\]
in such a way that $\GL_m(F_{\o v}\times F_v)$ embeds in
$\GL_{m+1}(F_{\o v}\times F_v)$ and
$\GL_m(F_{\o v}\times F_v)\simeq (F_{\o v})^m \times_{(D)} (F_v)^m$.

\item {\bbf The second, called ``quantum'', refers at first sight to the algebraic dimensions ``$j$\/}, i.e. Galois extension degrees corresponding to  global residue degrees
(see section~1.1), of the classifying bisemispace
$\BGL^{(Q)}(F^{2i}_{\o v^1}\times F^{2i}_{v^1})$ of
$\GL^{(Q)}(F^{2i}_{\o v^1}\times F^{2i}_{v^1})$ where
\[
\GL^{(Q)}(F^{2i}_{\o v^1}\times F^{2i}_{v^1})=
\lim_{j=1\to r} \GL_j^{(Q)}(F^{2i}_{\o v^1}\times F^{2i}_{v^1})\;, \qquad
\forall\ i\ , \; 1\le i\le n\;, \]
in such a way that:
\Be
\item $\GL_1^{(Q)}(F^{2i}_{\o v^1}\times F^{2i}_{v^1})=P_{2i}(F_{\o v^1}\times  F_{v^1})$ is the unitary, i.e. parabolic, bilinear semigroup 	of
$\GL_{2i}(F_{\o v^1}\times  F_{v^1})
$;

\item $\GL_j^{(Q)}(F^{2i}_{\o v^1}\times F^{2i}_{v^1})=
\GL_{2i}(F_{\o v_j}\times F_{v_j}\simeq 
(F^{2i}_{\o v_j}\times F^{2i}_{v_j})$ where the integer ``$j$'' denotes a global residue degree and the integer ``$2i$'' denotes a geometric dimension.

\item $\GL_j^{(Q)}(F^{2i}_{\o v^1}\times F^{2i}_{v^1})\subset\GL_{j+1}^{(Q)}(F^{2i}_{\o v^1}\times F^{2i}_{v^1})$;

\item $\GL_j^{(Q)}(F^{2i}_{\o v^1}\times F^{2i}_{v^1})\subset\GL_{j}^{(Q)}(F^{2i+1}_{\o v^1}\times F^{2i+1}_{v^1})$: geometric inclusion.
\Ee\Ee
\vskip 11pt

\subsection{Lemma}

{\em
The set of ``quantum'' infinite general bilinear semigroups
\[
\L\{\GL^{(Q)}(F^{2i}_{\o v^1}\times F^{2i}_{v^1})\R\}_i=
\L\{\lim\limits_{j=1\to r\le\infty }\GL_j^{(Q)}(F^{2i}_{\o v^1}\times F^{2i}_{v^1})\R\}_i\]
corresponds to the ``classical'' infinite general bilinear semigroups
\[
\GL (F_{\o v}\times F_v)=\underrightarrow{\lim} \GL_m (F_{\o v}\times F_v)\]
 where
$F_v=\{F_{v_1},\dots,F_{v_j},\dots,F_{v_r}\}$
\resp{$F_{\o v}=\{F_{\o v_1},\dots,F_{\o v_j},\dots,F_{\o v_r}\}$} is the set of $r$ classes of archimedean pseudoramified real completions and
$F_{v^1}=\{ F_{v^1_1},\dots,F_{v^1_j},\dots,F_{v^1_r}\}$ is the corresponding set of unitary completions.
}
\vskip 11pt

\bpr
The ``quantum'' infinite bilinear semigroup
$\GL^{(Q)}(F^{2i}_{\o v^1}\times F^{2i}_{v^1})$ generates the set:
\[
G^{(1),(Q)}(F^{2i}_{\o v^1}\times F^{2i}_{v^1})\subset \dots\subset
G^{(j),(Q)}(F^{2i}_{\o v^1}\times F^{2i}_{v^1})\subset \dots\subset
G^{(r),(Q)}(F^{2i}_{\o v^1}\times F^{2i}_{v^1})\;, \quad 1\le i\le n\;, \]
of embedded (abstract) bisemispaces which are respectively (isomorphic to) the classes of products, right by left, of embedded algebraic semitori (increasing algebraic filtration):
\[
F^{2i}_{\o v_1}\times_{(D)}F^{2i}_{v_1}\subset\dots\subset
F^{2i}_{\o v_j}\times_{(D)}F^{2i}_{v_j}\subset\dots\subset
F^{2i}_{\o v_r}\times_{(D)}F^{2i}_{v_r}\]
since $F^{2i}_{v_j}=j\times F^{2i}_{v_1}$.

As the geometric dimension, given by the integer ``$i$'', varies, we have $n=n_1+\dots+i+n_s$ such increasing filtrations with $n_s\to \infty $.

On the other hand, the ``classical'' infinite general bilinear semigroup
$\GL (F_{\o v}\times F_v)$ generates the set:
\[
G^{(1)}(F_{\o v}\times F_v)\subset \dots\subset
G^{(m)}(F_{\o v}\times F_v)\subset \dots\subset
G^{(2n_s)}(F_{\o v}\times F_v)\]
of embedded (abstract) bisemispaces which are $1\dots m\dots 2n_s$-dimensional products, right by left, of symmetric towers of increasing algebraic semitori.

$G^{(2i)}(F_{\o v}\times F_v)$ is then the $2i$-th algebraic filtration:
\[
F^{2i}_{\o v_1}\times_{(D)} F^{2i}_{v_1} \subset \dots\subset
F^{2i}_{\o v_j}\times_{(D)} F^{2i}_{v_j} \subset \dots\subset
F^{2i}_{\o v_r}\times_{(D)} F^{2i}_{v_r} \;,\]
i.e. $\GL^{(Q)}(F^{2i}_{\o v^1}\times F^{2i}_{v^1})$ which corresponds to the (functional) representation space\linebreak $\fREPSP(\GL_{2i}(F_{\o v}\times F_{v}))$ of the bilinear semigroup
$\GL_{2i}(F_{\o v}\times F_{v})$.

So, the quantum infinite general bilinear semigroup
\[
\GL^{(Q)}(F^{2i}_{\o v^1}\times F^{2i}_{v^1})=
\lim\limits_{j=1\to r\to\infty }(\GL^{(Q)}_j(F^{2i}_{\o v^1}\times F^{2i}_{v^1}))\]
is $\GL_{2i}(F_{\o v}\times F_{v})$.

And, the set of ``quantum'' infinite bilinear semigroups
$\{\GL^{(Q)}(F^{1}_{\o v^1}\times F^{1}_{v^1}),\dots,
\GL^{(Q)}(F^{m}_{\o v^1}\times F^{m}_{v^1}),\dots,
\GL^{(Q)}(F^{2n}_{\o v^1}\times F^{2n}_{v^1})\}$
corresponds to the ``classical'' infinite bilinear semigroup
\be
\GL (F_{\o v}\times F_{v})=\underrightarrow{\lim}
\GL_m (F_{\o v}\times F_{v})\;.\tag*{\eop}
\ee
\vskip 11pt

\subsection{Proposition}

{\em
{\bbf The classical (and quantum) infinite bilinear semigroup
\[
\GL(F_{\o v}\times F_{v})=\underset{i}{\underrightarrow{\lim}}\GL_{2i}(F_{\o v}\times F_{v})\] corresponds to the (partially) reducible (functional) representation space\/}\linebreak
$
\RED\fREPSP (\GL_{2n=2+\dots+2i+\dots+2n_s}(F_{\o v}\times F_v))$ of the bilinear semigroup
$\GL_{2n}(F_{\o v}\times F_v)$ with $n\to\infty $.
}
\vskip 11pt

\bpr
The infinite bilinear semigroup
$\GL(F_{\o v}\times F_{v})$ is the disjoint union of the
$\GL_{2i}(F_{\o v}\times F_{v})$ modulo an equivalence relation together with morphisms
\[ mg\ell_{2i} : \qquad \GL_{2i}(F_{\o v}\times F_v) \quad \To \quad 
 \GL_{2n}(F_{\o v}\times F_v) \]
of $ \GL_{2i}(F_{\o v}\times F_v)$ into
$ \GL_{2n}(F_{\o v}\times F_v)$.

So, we have that:
\begin{align*}
\GL(F_{\o v}\times F_{v})
&=
\GL_2(F_{\o v}\times F_{v})\cup\dots\cup
\GL_{2i}(F_{\o v}\times F_{v})\cup\dots\cup
\GL_{2n_s}(F_{\o v}\times F_{v})\\
&=\underrightarrow{\lim} \GL_{2i}(F_{\o v}\times F_{v})\end{align*}
with $ \GL_{2i+2}(F_{\o v}\times F_{v})\subset
 \GL_{2i+2}(F_{\o v}\times F_{v})$.
 
And, thus, $\GL(F_{\o v}\times F_{v})$ generates to the (partially) reducible (functional) representation space
$\RED\fREPSP(\GL_{2n}(F_{\o v}\times F_{v}))$ which decomposes according to the partition $2n=2+\dots+2i+\dots+2n_s$:
\begin{multline*}
\RED\fREPSP (\GL_{2n}(F_{\o v}\times F_{v}))=
\fREPSP (\GL_{2}(F_{\o v}\times F_{v})) \\
\boxplus \dots \boxplus \fREPSP (\GL_{2i}(F_{\o v}\times F_{v})) \\
 \boxplus \dots\boxplus \fREPSP (\GL_{2n_s}(F_{\o v}\times F_{v})) \end{multline*}
as introduced in \cite{Pie2}.

{\bbf Summarizing, we have\/}:
\begin{align*}
\GL(F_{\o v}\times F_{v})
&= \underrightarrow{\lim} \GL_{2i}(F_{\o v}\times F_{v})\simeq \RED\fREPSP ( \GL_{2n}(F_{\o v}\times F_{v}))\\
&= \{ \GL^{(Q)}(F^{2i}_{\o v^1}\times F^{2i}_{v^1})\}_i=\underrightarrow{\lim} \GL^{(Q)}(F^{2i}_{\o v^1}\times F^{2i}_{v^1})\;,
\end{align*}
$\forall\ i\in$ partition $2n=2+\dots+2i+\dots+2n_s$, $n\to\infty $.\epr
 \vskip 11pt

\subsection{The classifying bisemispace $\BGL(F_{\o v}\times F_v)$}

The classifying bisemispace $\BGL(F_{\o v}\times F_v)$ of
$\GL(F_{\o v}\times F_v)$ is the quotient of a weakly contractible bisemispace 
$\EGL(F_{\o v}\times F_v)$ by a free action of
$\GL(F_{\o v}\times F_v)$; that is to say, generalizing the homotopy linear definition of a classifying space, the contractible bisemispace
$\EGL(F_{\o v}\times F_v)$ is the total bisemispace of a universal principal
$\GL(F_{\o v}\times F_v)$-bibundle over the classifying bisemispace
$\BGL(F_{\o v}\times F_v)$ given by the continuous mapping
\[
\GD_\ell: \qquad
\EGL(F_{\o v}\times F_v)\quad \To \quad \BGL(F_{\o v}\times F_v)\;.\]
This approach is more basic than the condition implying classically that the higher homotopy groups are trivial (or vanish).
\vskip 11pt

\subsection{Proposition}

{\em
The (continuous) mapping
\[
\GD: \qquad
\EGL(F_{\o v}\times F_v)\quad \To \quad \BGL(F_{\o v}\times F_v)\]
of the principal 
$\GL(F_{\o v}\times F_v)$-bibundle over the classifying bisemispace
$\BGL(F_{\o v}\times F_v)$ is a homotopy map corresponding to the deformations of the Galois compact representation of
$\BGL(\wt F_{\o v}\times \wt F_v)$ given by the (bi)fibres of $\GD_\ell$, $\forall\ \ell$, $1\le \ell\le \infty $.
}
\vskip 11pt

\bpr
Referring to proposition~2.15 introducing homotopy maps as deformations of Galois representations of linear semigroups, we see that the map $\GD_\ell$ of the principal
$\GL(F_{\o v}\times F_v)$-bibundle corresponds to a deformation of the Galois compact representation of $\GL(\wt F_{\o v}\times \wt F_v)$:
\[
\GD_\ell: \qquad
\GL(F_{\o {v+\ell}}\times F_{v+\ell})\quad \To \quad \GL(F_{\o v}\times F_v)\]
in such a way that the kernel
$\GL(\delta F_{\o {v+\ell}}\times \delta F_{v+\ell})$ of $\GD_\ell$ is responsible for the increase of sets of powers of ``$\ell$'' biquanta to
$\GL(F_{\o {v}}\times F_{v})$.

$\GD_\ell$  then belongs to an equivalence class of homotopy maps, given by deformations of the Galois compact representation of
$\GL(\wt F_{\o {v}}\times \wt F_{v})$.

And, the set $\{\GD_\ell\}_\ell$ of all equivalence classes of homotopy maps is the continuous mapping
\be
\GD: \qquad
\EGL(F_{\o v}\times F_v)\quad \To \quad \BGL(F_{\o v}\times F_v)\;.
\tag*{\eop}\ee
\vskip 11pt

\subsection{Corollary}

{\em
{\bbf The classifying bisemispace
$\BGL(F_{\o {v}}\times F_{v})$ is the base bisemispace of all equivalence classes of deformations of the Galois representation of
$\GL(\wt F_{\o {v}}\times \wt F_{v})$ given by the kernels
$\GL(\delta F_{\o {v+\ell}}\times \delta F_{v+\ell})$ of the maps\/}
\[
\GD_\ell: \qquad
\GL(F_{\o {v+\ell}}\times F_{v+\ell})\quad \To \quad \GL(F_{\o v}\times F_v)\;, \quad 1\le \ell\le \infty \;.\]
}
\vskip 11pt

\subsection{The ``plus'' construction of Quillen}

The ``plus'' construction, adapted to the bilinear case of the Langlands global program, leads to consider a map
\[
\BG(1): \qquad
\BGL(F_{\o {v}}\times F_{v})\quad \To \quad \BGL(F_{\o v}\times F_v)^+\;, \]
unique up to homotopy, such that:
\Bena
\item the kernel of $\Pi _1(\BG(1))$ be one-dimensional deformations of the Galois compact representation of
$\GL(\wt F_{\o {v}}\times \wt F_{v})$;

\item the homotopy fibre of $\BG(1)$ has the same integral homology as a point (or\linebreak $\BGL( F_{\o {v}}\times  F_{v})$ and $\BGL(F_{\o {v}}\times  F_{v})^+$ have the same integral homology).
\Ee
\vskip 11pt

\subsection{Proposition}

{\em
Let $\{\GL^{(1)}(\delta F_{\o {v+\ell}}\times \delta F_{v+\ell})\}_\ell$ denote the set of kernels of the maps:
\[
\GD(1)_\ell: \qquad
\GL^{(1)}(F_{\o {v+\ell}}\times F_{v+\ell})\quad \To \quad \GL^{(1)}(F_{\o v}\times F_v)\;, 
\qquad 1\le\ell\le\infty \;,\]
where $\GL^{(1)}(F_{\o v}\times F_v)$ denote the set of one-dimensional irreducible components of the bisemispace $\GL(F_{\o v}\times F_v)$.

Then, the classifying bisemispace $\BGL(F_{\o v}\times F_v)^+$ is the base bisemispace of all equivalence classes of one-dimensional deformations of the Galois compact representation of
$\GL(\wt F_{\o v}\times \wt F_v)$ given by the kernels
$\{\GL^{(1)}(\delta F_{\o {v+\ell}}\times \delta F_{v+\ell})\}_\ell$ of the maps $\GD(1)_\ell$.
}
\vskip 11pt

\bpr
The classifying bisemispace $\BGL(F_{\o v}\times F_v)$ is the base bisemispace of the principal 
$\GL(F_{\o v}\times F_v)$-bibundle whose map is:
\[
\GD: \qquad
\EGL(F_{\o {v}}\times F_{v})\quad \To \quad \BGL(F_{\o v}\times F_v)\;.\]
Similarly, the ``plus'' classifying bisemispace
$\BGL(F_{\o v}\times F_v)^+$ must be the base bisemispace of the principal
$\GL^{(1)}(F_{\o v}\times F_v)$-bibundle whose map is:
\[
\GD(1): \qquad
\EGL(F_{\o {v}}\times F_{v})^+\quad \To \quad \BGL(F_{\o v}\times F_v)^+\;,\]
where $\EGL(F_{\o {v}}\times F_{v})^+$ is the total bisemispace verifying the equivalent conditions:
\Bean
\item $\EGL(F_{\o {v}}\times F_{v})^+=\Pi _1(\BGL(F_{\o {v}}\times F_{v})^+)$;
\item $\EGL(F_{\o {v}}\times F_{v})^+$ corresponds to all equivalence classes of one-dimensional deformations 
$\GL^{(1)}(F_{\o {v+\ell}}\times F_{v+\ell})$ of the Galois compact representations of
$\GL(\wt F_{\o {v}}\times \wt F_{v})$:
\begin{align}
\EGL(F_{\o {v}}\times F_{v})^+ &= \{ \GL^{(1)}(F_{\o {v+\ell }}\times F_{v+\ell })\}_\ell \notag\\
\and \qquad 
\Pi _1(\BGL(F_{\o {v}}\times F_{v})^+) &= \GL(F_{\o {v}}\times F_{v})\big/ \GL^{(1)}(F_{\o {v}}\times F_{v})\;.\tag*{\eop}
\end{align}
\Ee
\vskip 11pt


\subsection{Corollary}

{\em
{\bbf The ``$+$'' construction leads to the following commutative diagram\/}:
\[
\begin{psmatrix}[colsep=1cm,rowsep=.6cm]
\EGL(F_{\o {v}}\times F_{v}) & & \EGL(F_{\o {v}}\times F_{v})^+ \\[15pt]
\BGL(F_{\o {v}}\times F_{v}) && \BGL(F_{\o {v}}\times F_{v})^+
\psset{arrows=->,nodesep=5pt}
\everypsbox{\scriptstyle}
\ncline{1,1}{1,3}^{\EG(1)}
\ncline{2,1}{2,3}^{\BG(1)}
\ncline{1,1}{2,1}>{\GD}
\ncline{1,3}{2,3}>{\GD(1)}
\end{psmatrix}
\]
where $\EG(1)$ is the map:
\[
\EG(1): \qquad
\{\GL(F_{\o {v+\ell }}\times F_{v+\ell })\}_\ell\quad \To \quad 
\{\GL^{(1)}(F_{\o v+\ell }\times F_{v+\ell })\}_\ell\;,\]
from deformations $\{\GL(F_{\o {v+\ell }}\times F_{v+\ell })\}_\ell$ of the Galois compact representations
$\GL(\wt F_{\o {v }}\times \wt F_{v })$ to one-dimensional deformations
$\{\GL^{1)}(F_{\o {v+\ell }}\times F_{v+\ell })\}_\ell$ of the Galois compact representations of
$\GL(\wt F_{\o {v }}\times \wt F_{v })$.
}
\vskip 11pt

\subsection{Proposition}

{\em 
{\bbf The bilinear version of the algebraic $K$-theory of Quillen adapted to the Langlands global program is\/}:
\[
K^{(2i)}(G^{(2n)}_{\rm red}(F_{\o v}\times F_v))=
\Pi _{2i}(\BGL(F_{\o v}\times F_v)^+)\;,\]
where $G^{(2n)}_{\rm red}(F_{\o v}\times F_v)=\RED\fREPSP(\GL_{2n}(F_{\o v}\times F_v))$ is the (partially) reducible (functional) representation space of the bilinear semigroup of matrices
 $\GL_{2n}(F_{\o v}\times F_v)$, in such a way that:
 \Bean
 \item  the partition $2n=2+\dots+2i+\dots+2n_s$ of the geometric dimension $2n$, $n\le\infty $, refers to the reducibility of
$\GL_{2n}(F_{\o v}\times F_v)$;
\item the dimension $2i$ of the bisemigroup of homotopy
$\Pi _{2i}(\BGL(F_{\o v}\times F_v)^+)$ must be inferior or equal to each term of the partition of $2n$ in order that this homotopy bisemigroup be non trivial.
\Ee
}
\vskip 11pt

\bpr
Referring to proposition~3.6, the infinite bisemigroup $\GL(F_{\o v}\times F_v)$ is
\[
\GL(F_{\o v}\times F_v)
=\underrightarrow{\lim} \GL_{2i}(F_{\o v}\times F_v)
\simeq \RED\fREPSP(\GL_{2n}(F_{\o v}\times F_v))\;, \]
i.e. the decomposition of the partially reducible bisemivariety
$G^{(2n)}_{\rm red}(F_{\o v}\times F_v)$ into
\[ 
G^{(2n)}_{\rm red}(F_{\o v}\times F_v)= G^{(2)}(F_{\o v}\times F_v)\oplus \dots\oplus
G^{(2i)}(F_{\o v}\times F_v)\oplus \dots\oplus G^{(2n_s)}(F_{\o v}\times F_v)\;.\]

It is then evident that the homotopy bisemigroup
$\Pi _{2i}(\BGL(F_{\o v}\times F_v)^+)$ is null for the bisemivarieties
$G^{(2h)}_{\rm red}(F_{\o v}\times F_v)$ whose geometric dimension $h<i$.\epr
\vskip 11pt

\subsection{Corollary}

{\em The bilinear version of the algebraic $K$-theory
\[
K^{2i}(G^{(2n)}_{\rm red}(F_{\o v}\times F_v))
=\Pi _{2i}(\BGL(F_{\o v}\times F_v)^+)\;, \]
relative to the Langlands global program and corresponding to a higher version of this global program, is in one-to-one correspondence with the ``quantum'' bilinear version of the algebraic $K$-theory:
\[
K^{2i}(G^{(2n)}_{\rm red}(F_{\o v}\times F_v))
=\Pi _{2i}(\BGL^{(Q)}(F^{2i}_{\o v^1}\times F^{2i}_{v^1})^+)\;. \]
}
\vskip 11pt

\bpr
Indeed, according to lemma~3.5 and proposition~3.6, we have that the ``classical'' infinite bisemigroup
\[
\GL(F_{\o v}\times F_v)=\underset{i}{\underrightarrow{\lim}}\GL_{2i}(F_{\o v}\times F_v)\]
is equal to its ``quantum'' version given by:
\[
\GL(F_{\o v}\times F_v)=\underset{i}{\underrightarrow{\lim}}\GL^{(Q)}(F^{2i}_{\o v^1}\times F^{2i}_{v^1})\]
for every $i$ belonging to the partition of $2n$ associated with the dimensions of the reducibility of the bisemivariety
\[
G^{(2n)}_{\rm red}(F_{\o v}\times F_v)
=\RED\fREPSP(\GL_{2n}(F_{\o v}\times F_v))\;. \]

The ``quantum'' version works explicitly with the algebraic dimensions ``$j$'' by the mapping:
\[\underrightarrow{\GL^{(Q)}_j} : \qquad
F^{2i}_{\o v^1}\times F^{2i}_{v^1} \quad \To \quad
\underset{j}{\underrightarrow{\lim}}\GL^{(Q)}_j(F^{2i}_{\o v^1}\times F^{2i}_{v^1})
\begin{aligned}[t]
&= \GL^{(Q)}(F^{2i}_{\o v^1}\times F^{2i}_{v^1})\\
&= \GL_{2i}(F_{\o v}\times F_{v})\end{aligned}
\]
while the classical version is based on the fibre bundle:
\[
\GL_{2i}: \qquad
F_{\o {v}}\times F_{v}\quad \To \quad 
\GL_{2i}(F_{\o v }\times F_{v})\;,\]
with ``geometric'' bifibre $\GL_{2i-1}(F_{\o v }\times F_{v})$.\epr
\vskip 11pt

\subsection{Proposition}

{\em
{\bbf The higher version of the Langlands global program
\[
K^{2i}(G^{(2n)}_{\rm red}(F_{\o v}\times F_v)) 
=\Pi _{2i}(\BGL(F_{\o v}\times F_{v})^+)\]
implies that the equivalence classes of $2i$-dimensional deformations of the Galois compact representations of the reducible bilinear semigroup $\GL_{2n}(F_{\o v}\times F_v)$ result from quantum homomorphisms of the global coefficient bisemiring $F_{\o v}\times F_v$\/}.
}
\vskip 11pt

\bpr
Referring to proposition~3.11, the ``plus'' classifying bisemispace
$\BGL(F_{\o v}\times F_{v})^+$ is the base bisemispace of the principal
$\GL^{(1)}(F_{\o v }\times F_{v})$-bibundle in such a way that the total bisemispace
$\EGL(F_{\o v }\times F_{v})^+$ verifies:
\[ \EGL(F_{\o v }\times F_{v})=\Pi _1(\BGL(F_{\o v }\times F_{v})^+)\]
and corresponds to quantum homomorphisms of the global coefficient bisemiring $F_{\o v }\times F_{v}$:
\[ 
Qh_{F_{\o {v+\ell } }\times F_{v+\ell }\to F_{\o v }\times F_{v}}: \qquad
F_{\o {v+\ell }} \times F_{v+\ell } \quad \To \quad F_{\o v }\times F_{v}\]
according to section~2.7.

Then, the $2i$-th homotopy bisemigroup $\Pi _{2i}(\BGL(F_{\o v }\times F_{v})^+)$, describing the equivalence classes of $2i$-dimensional deformations of the Galois representations of
$\GL_{2n}(F_{\o v }\times F_{v})$, implies the monomorphism:
\be
\Pi _1(\BGL(F_{\o v }\times F_{v})^+) \quad \To \quad
\Pi _{2i}(\BGL(F_{\o v }\times F_{v})^+)\;.\tag*{\eop}
\ee
\vskip 11pt

\subsection{Restricted Chern character}

The Chern character restricted to the class $C^i$ in the higher bilinear $K$-cohomology is given by the homomorphism:
\[
C^i(G^{(2n)}_{\rm red}(F_{\o v}\times F_v)) : \qquad
K^{2i}(G^{(2n)}_{\rm red}(F_{\o v}\times F_v)) \quad \To \quad
H^{2i}(G^{(2n)}_{\rm red}(F_{\o v}\times F_v),G^{(2i)}(F_{\o v}\times F_v))\]
where
$G^{(2n)}_{\rm red}(F_{\o v}\times F_v)$ is the reducible representation of 
$\GL_{2n}(F_{\o v}\times F_v)$, i.e. a compact bisemivariety decomposing into:
\[
G^{(2n)}_{\rm red}(F_{\o v}\times F_v)=G^{(2)}(F_{\o v}\times F_v)\oplus\dots\oplus
G^{(2i)}(F_{\o v}\times F_v)\oplus\dots\oplus
G^{(2n_s)}(F_{\o v}\times F_v)\;.\]
\vskip 11pt

\subsection{Proposition}

{\em The higher bilinear $K$-cohomology restricted to the class ``$2i$'' implies the ``higher'' bilinear semigroup homomorphisms of Hurewicz, i.e. {\bbf a higher restricted $\Pi $-cohomology\/}:
\[
hhH\RL : \qquad
\Pi _{2i} ( \BGL ( F_{\o v }\times F_{v})^+)\quad \To \quad
H^{2i}(G^{(2n)}_{\rm red}(F_{\o v}\times F_v),\zit\times_{(D)}\zit )\]
from the ``Galois'' higher bilinear homotopy $\Pi _{2i}(\cdot)$ into the entire ``higher'' bilinear cohomology $H^{2i}(\cdot)$.

This leads to the commutative diagram:
\vspace*{20pt}
\[
\begin{psmatrix}[colsep=1cm,rowsep=.6cm]
K^{2i}(G^{(2n)}_{\rm red}(F_{\o v}\times F_v)) && \Pi _{2i}(\BGL (F_{\o v}\times F_v)^+)\\[15pt]
& H^{2i}(G^{(2n)}_{\rm red}(F_{\o v}\times F_v),\zit\times_{(D)}\zit )
\psset{arrows=->,nodesep=5pt}
\everypsbox{\scriptstyle}
\ncline{1,1}{1,3}^{\Bsm \text{higher algebraic}\\\text{$K$-theory}\\ -\Esm}_-
\ncline{1,1}{2,2}<{\Bsm \text{``Chern'' higher}\\ \text{restricted character}\Esm}
\ncline{2,2}{1,3}>{\Bsm \text{inverse restricted}\\ \text{higher $\Pi $-cohomology}\Esm}
\end{psmatrix}
\]
in such a way that the classes of the entire bilinear cohomology
$H^{2i}(G^{(2n)}_{\rm red}(F_{\o v}\times F_v),\zit\times_{(D)}\zit)$
refer to a bisemilattice deformed by the homotopy classes of maps of
$\Pi _{2i}(\BGL(F_{\o v}\times F_v)^+)$, corresponding to lift of quantum deformations of the Galois representations of $\GL_{2n}^{(\rm red)}(\wt F_{\o v}\times \wt F_v)$.
}
\vskip 11pt

\bpr
The higher algebraic $K$-theory
\[
K^{2i}(G^{(2n)}_{\rm red}(F_{\o v}\times F_v))=\Pi _{2i}(\BGL (F_{\o v}\times F_v)^+)\;, \]
relative to the Langlands global program,

together with the restricted ``higher'' $K$-cohomology
\[
C^i(G^{(2n)}_{\rm red}(F_{\o v}\times F_v)) : \qquad
K^{2i}(G^{(2n)}_{\rm red}(F_{\o v}\times F_v)) \quad \To \quad
H^{2i}(G^{(2n)}_{\rm red}(F_{\o v}\times F_v),G^{(2i)}(F_{\o v}\times F_v))\]
implies the restricted higher $\Pi $-cohomology:
\be
\Pi _{2i}(\BGL (F_{\o v}\times F_v)^+)\quad \To \quad
H^{2i}(G^{(2n)}_{\rm red}(F_{\o v}\times F_v),\zit\times_{(D)}\zit )\;.\tag*{\eop}
\ee
\vskip 11pt

\subsection{Proposition}

{\em
{\bbf The total Chern character in the bilinear $K$-cohomology of the reducible representation
$G^{(2n)}_{\rm red}(F_{\o v}\times F_v)$ of $\GL_{2n}( F_{\o v}\times  F_v)$\/}:
\[
ch^*(G^{(2n)}_{\rm red}(F_{\o v}\times F_v)): \qquad
K^*(G^{(2n)}_{\rm red}(F_{\o v}\times F_v)) \quad \To \quad
H^*(G^{(2n)}_{\rm red}(F_{\o v}\times F_v),G^{*}(F_{\o v}\times F_v))\;, \]
where $*$ is the partition $2n=2+\dots+2i+\dot+2n_s$ of $2n$,

implies the total higher algebraic $K$-theory:
\[
K^*(G^{(2n)}_{\rm red}(F_{\o v}\times F_v))=\Pi_* ( \BGL (F_{\o v}\times F_v)^+ )\;.\]
}
\vskip 11pt

\bpr This results from the preceding sections.\epr
\vskip 11pt

\subsection{Corollary}

{\em The total higher algebraic $K$-theory associated with the reducible global program of Langlands is based on the commutative diagram:
\vspace*{20pt}
\[
\begin{psmatrix}[colsep=1cm,rowsep=.6cm]
K^{*}(G^{(2n)}_{\rm red}(F_{\o v}\times F_v)) && \Pi _{*}(\BGL (F_{\o v}\times F_v)^+)\\[15pt]
& H^{*}(G^{(2n)}_{\rm red}(F_{\o v}\times F_v),G^{*}(F_{\o v}\times F_v) )
\psset{arrows=->,nodesep=5pt}
\everypsbox{\scriptstyle}
\ncline{1,1}{1,3}^{\Bsm \text{total higher algebraic}\\\text{$K$-theory}\\ -\Esm}_-
\ncline{1,1}{2,2}<{\Bsm \text{``Chern'' total }\\ \text{higher character}\Esm}>{ch^*}
\ncline{2,2}{1,3}>{\Bsm \text{higher inverse}\\ \text{$\Pi $-cohomology}\Esm}
\end{psmatrix}
\]
}
\vskip 11pt

\subsection{Higher bilinear algebraic $K$-theory referring to cohomology}

As a lower bilinear (algebraic) $K$-theory referring to cohomotopy was introduced in corollary~2.22, {\bbf a higher bilinear algebraic $K$-theory relative to cohomotopy can be introduced by the equality\/}:
\[
K_{2i}(G^{(2n)}_{\rm red}(F_{\o v}\times F_v))=\Pi^{2i} ( \BGL (F_{\o v}\times F_v)^+ )\]
where
$\Pi^{2i} ( \BGL (F_{\o v}\times F_v)^+ )$ are the cohomotopy equivalence classes of $2i$-dimensional deformations of the Galois reducible representations of $\GL_{2n}(\wt  F_{\o v}\times  \wt F_v)$.
\vskip 11pt

\subsection{Proposition}

{\em
{\bbf The higher bilinear algebraic $K$-theory relative to cohomotopy implies the commutative diagram\/}:
\vspace*{20pt}
\[
\begin{psmatrix}[colsep=1cm,rowsep=.6cm]
K_{2i}(G^{(2n)}_{\rm red}(F_{\o v}\times F_v)) && \Pi ^{2i}(\BGL (F_{\o v}\times F_v)^+)\\[20pt]
& H_{2i}(G^{(2n)}_{\rm red}(F_{\o v}\times F_v),\zit\times_{(D)}\zit )
\psset{arrows=->,nodesep=5pt}
\everypsbox{\scriptstyle}
\ncline{1,1}{1,3}^{\Bsm \text{ higher algebraic}\\\text{$K$-theory}\\ -\Esm}_{\Bsm -\\ \text{referring to}\\ \text{cohomotopy}\Esm}
\ncline{1,1}{2,2}<{\Bsm \text{``Chern'' higher }\\ \text{restricted character}\\ \text{relative to homology}\Esm}
\ncline{2,2}{1,3}>{\Bsm \text{higher inverse}\\ \text{restricted $\Pi $-homology}\Esm}
\end{psmatrix}
\]
where:
\Bean
\item $C_{i}(G^{(2n)}_{\rm red}(F_{\o v}\times F_v)): K_{2i}(G^{(2n)}_{\rm red}(F_{\o v}\times F_v))\to
H_{2i}(G^{(2n)}(F_{\o v}\times F_v),\zit\times_{(D)}\zit)$ is the Chern higher restricted character relative to the $K$-homology where
$H_{2i}(G^{(2n)}_{\rm red}(F_{\o v}\times F_v),\zit\times_{(D)}\zit)$ is the entire bilinear homology of the reducible bisemivariety $G^{(2n)}_{\rm red}(F_{\o v}\times F_v)$ in the real bisemilattice deformed by the cohomotopy classes of maps of\linebreak 
$\Pi ^{2i}(\BGL (F_{\o v}\times F_v)^+)$, corresponding to lifts of inverse quantum deformations of the Galois representations of $\GL_{2n}^{(\rm red)}(\wt F_{\o v}\times \wt F_v)$;

\item $hhCH^{(2i)}\RL: \Pi ^{2i}(\BGL (F_{\o v}\times F_v)^+)\to
H_{2i}(G^{(2n)}_{\rm red}(F_{\o v}\times F_v),\zit\times_{(D)}\zit)$ is the Hurewicz higher homomorphism relative to cohomotopy.
\Ee
}
\vskip 11pt

\paragraph{Sketch of proof}: \hspace{5mm} 
The higher bilinear algebraic $K$-theory referring to cohomotopy given by the equality
\[
K_{2i} ( G^{(2n)}_{\rm red} (F_{\o v}\times F_v))=\Pi ^{2i}(\BGL (F_{\o v}\times F_v)^+)\]
together with the Chern higher restricted character relative to homology implies the Hurewicz higher homomorphism $hhCH\RL$.\epr
\vskip 11pt

\subsection{Corollary}

{\em
{\bbf The total higher bilinear algebraic $K$-theory relative to cohomotopy is given by the equality\/}:
\[
K_{*}(G^{(2n)}_{\rm red}(F_{\o v}\times F_v))=\Pi ^{*}(\BGL (F_{\o v}\times F_v)^+)\]
and implies the commutative diagram:
\vspace*{20pt}
\[
\begin{psmatrix}[colsep=1cm,rowsep=.6cm]
K_{*}(G^{(2n)}_{\rm red}(F_{\o v}\times F_v)) && \Pi ^{*}(\BGL (F_{\o v}\times F_v)^+)\\[20pt]
& H_{*}(G^{(2n)}_{\rm red}(F_{\o v}\times F_v),\zit\times_{(D)} \zit) )
\psset{arrows=->,nodesep=5pt}
\everypsbox{\scriptstyle}
\ncline{1,1}{1,3}^{\Bsm \text{total higher}\\\text{algebraic $K$-theory}\\ -\Esm}_{\Bsm -\\ \text{relative to}\\ \text{cohomotopy}\Esm}
\ncline{1,1}{2,2}<{\Bsm \text{``Chern'' total higher }\\ \text{character relative to homology}\Esm}>{ch_*}
\ncline{2,2}{1,3}>{\Bsm \text{higher inverse}\\ \text{$\Pi $-homology}\Esm}
\end{psmatrix}
\]
where:
\Bean
\item $Ch_{*}(G^{(2n)}_{\rm red}(F_{\o v}\times F_v)): K_{*}(G^{(2n)}_{\rm red}(F_{\o v}\times F_v))\to
H_{*}(G^{(2n)}(F_{\o v}\times F_v),\zit\times_{(D)}\zit)$ is the Chern character of the reducible bisemivariety $G^{(2n)}_{\rm red}(F_{\o v}\times F_v)$ in the higher bilinear $K$-homology;

\item $hhCH^{(*)}\RL: H_{*}(G^{(2n)}_{\rm red}(F_{\o v}\times F_v),\zit\times_{(D)}\zit)\to
 \Pi ^{*}(\BGL (F_{\o v}\times F_v)^+)$
 is the corresponding Hurewicz total higher inverse homomorphism relative to cohomotopy.
\Ee
}
\vskip 11pt

\paragraph{Sketch of proof}: \hspace{5mm} 
The framework of the total higher bilinear algebraic $K$-theory relative to cohomotopy is similar to that of homotopy handled in proposition~3.18 and corollary~3.20.\epr

\section{Mixed higher bilinear algebraic $KK$-theories related to the  Langlands dynamical bilinear global program}

The lower and higher versions of the Langlands dynamical global program refer respectively to dynamical lower and higher bilinear (algebraic) $K$-theories related to {\bbf the existence of $K_*K^*$ functors on the categories of elliptic bioperators and (reducible) bisemisheaves
$FG^{(2n)}_{\rm red}(F_{\o v}\times F_v)$, being (reducible) functional representation spaces\/} of the (algebraic)  general bilinear semigroups $\GL_{2n}(F_{\o v}\times F_v)$.
\vskip 11pt


\subsection{Bilinear contracting fibres of tangent bibundles}

Let then $FG^{(2i)}(F_{\o v}\times F_v)=\FREPSP(\GL_{2i}(F_{\o v}\times F_v))$ denote the functional representation space of $\GL_{2i}(F_{\o v}\times F_v)$, $i\le n\le\infty $, which splits into:
\begin{multline*}
\qquad \FREPSP (\GL_{2i}(F_{\o v}\times F_v ) ) \\
=\FREPSP (\GL_{2k}(F_{\o v}\times F_v))\oplus
\FREPSP (\GL_{2i-2k}(F_{\o v}\times F_v))\;, \quad k\le i\;,\qquad\end{multline*}
in such a way that $\FREPSP(\GL_{2k}(F_{\o v}\times F_v))$ is the functional representation space of geometric dimension $2k$ of the bilinear semigroup
$\GL_{2k}(F_{\o v}\times F_v)$ on which acts the elliptic bioperator $D_R^{2k}\otimes D_L^{2k}$.

Let then $D_R^{2k}\otimes D_L^{2k}$ be the product of a right linear differential (elliptic) operator
$D_R^{2k}$ acting on $2k$ variables by its left equivalent $D_L^{2k}$ \cite{Sat}, \cite{Kash}.

This bioperator $D_R^{2k}\otimes D_L^{2k}$ is defined by its biaction:
\[ D_R^{2k}\otimes D_L^{2k}: \qquad
FG^{(2i)}(F_{\o v}\times F_v) \quad \To \quad FG^{(2i[2k])}(F_{\o v}\times F_v)\]
where 
$FG^{(2i[2k])}(F_{\o v}\times F_v)$ is the functional representation space of
$\GL_{2i}(F_{\o v}\times F_v)$ shifted in $(2k\times_{(D)}2k)$ dimensions, i.e. bisections of a $2i$-dimensional bisemisheaf shifted in $(2k\times_{(D)}2k)$ dimensions, of differentiable bifunctions on the abstract bisemivariety $G^{2i[2k]}(F_{\o v}\times F_v)$.

Referring to chapter~3 of \cite{Pie6}, the shifted bisemisheaf
$FG^{(2i[2k])}(F_{\o v}\times F_v)$ decomposes into:
\[
FG^{(2i[2k])}(F_{\o v}\times F_v)=(\Delta ^{2k}_R\times\Delta _L^{2k})\oplus FG^{(2i-2k)}(F_{\o v}\times F_v)\]
where:
\Bean
\item 
$\Delta ^{2k}_R\times\Delta _L^{2k}
\begin{aligned}[t]
&\simeq \Ad\fREPSP (\GL_{2k}(\rit\times \rit))\times \fREPSP(\GL_{2k}(F_{\o v}\times F_v))\\
&\simeq \fREPSP(\GL_{2k}(F_{\o v}\times \rit)(F_v\times \rit))\end{aligned}$

with $\Ad\fREPSP(\GL_{2k}(\rit\times \rit))$ being the adjoint functional representation space of
$\GL_{2k}(\rit\times \rit)$ corresponding to the biaction of the bioperator
$(D ^{2k}_R\otimes D _L^{2k})$ on the bisemisheaf $ FG^{(2k)}(F_{\o v}\times F_v)$ which is the functional representation space of $ \GL_{(2k)}(F_{\o v}\times F_v)$;

\item $ FG^{(2i-2k)}(F_{\o v}\times F_v)=\FREPSP(\GL_{2i-2k}(F_{\o v}\times F_v))$ is the $(2i-2k)$ (geometric) dimensional bisemisheaf being the functional representation space of 
$\GL_{2i-2k}(F_{\o v}\times F_v)$.
\Ee

In fact, $(\Delta ^{2k}_R\times\Delta _L^{2k})$ is the total bisemispace of the tangent (bi)bundle\linebreak
$\TAN(FG^{(2k)}(F_{\o v}\times F_v))$ to the bisemisheaf
$FG^{(2k)}(F_{\o v}\times F_v)$ of which bilinear fibre
\[
\Fs^{2k}\RL(\TAN)=(\AdF)\REPSP ( \GL_{2k}(\rit\times \rit ))\]
is isomorphic to the adjoint functional representation space of 
$\GL_{2k}(\rit\times \rit)$.

And, $\Aut(\TAN_e(FG^{(2k)}(F_{\o v}\times F_v)))$ is an open subset of the bilinear vector semispace of endomorphisms of
$\TAN_e(FG^{(2k)}(F_{\o v}\times F_v))$ at the identity element ``$e$'' in order to define differentials on it.

If the bisemisheaf
$FG^{(2k)}(F_{\o v}\times F_v)=\FREPSP(\GL_{2k}(F_{\o v}\times F_v))$ is 
the functional representation space over the 
abstract bisemivariety $G^{(2k)}(F_{\o v}\times F_v)$, then $\Delta ^{2k}_R\times\Delta _L^{2k}$ is in one-to-one correspondence with a Galois (bisemi)group(oid) of which isomorphism is generated by the bilinear fibre
$(\AdF)\REPSP(\GL_{2k}(\rit\times \rit))$.

Indeed, the Galois (bisemi)group(oid) associated with the shifted bisemisheaf
$(\Delta ^{2k}_R\times\Delta _L^{2k})$ refers essentially to the bilinear semigroup
$\GL(\Delta ^{2k}_R\times\Delta _L^{2k})$ of ``Galois'' automorphisms of
$(\Delta ^{2k}_R\times\Delta _L^{2k})$, i.e. to the bilinear semigroup of shifted ``Galois'' automorphisms of the base bisemisheaf
$G^{(2k)}(F_{\o v}\times F_v)=\fREPSP(\GL_{2k}(F_{\o v}\times F_v))$ endowed with a nontrivial fundamental bisemigroup $\Pi _1(G^{(2k)}(F_{\o v}\times F_v))$.
\vskip 11pt

\subsection{Proposition}

{\em
{\bbf The existence of a bilinear contracting fibre $\Fs^{2k}\RL(\TAN)$ in the tangent bibundle
$\TAN(FG^{(2k)}(F_{\o v}\times F_v))$ implies the homology\/}:
\[
H_{2k} (F G^{(2i[2k])}(F_{\o v}\times F_v),\Fs^{2k}\RL (\TAN ) )
\simeq \Ad\fREPSP (\GL_{2k}(\rit\times \rit ))\]
in such 
a way that  the cohomoloby of the shifted bisemisheaf
$FG^{(2i[2k])}(F_{\o v}\times F_v)$ be given by:
\begin{multline*}
H^{2k} (FG^{(2i[2k])}(F_{\o v}\times F_v),\Delta ^{2k}\RL )\\
\begin{aligned}
&=H_{2k} (FG^{(2i[2k])} (F_{\o v}\times F_v),\Fs^{2k}\RL (\TAN ))
\times H^{2k}(FG^{(2i[2k])}(F_{\o v}\times F_v),FG^{(2k)}(F_{\o v}\times F_v))\\
&= \FREPSP ( \GL_{2k}(F_{\o v}\times \rit )\times (F_v\times \rit ))
\end{aligned}
\end{multline*}
where $\Delta ^{2k}\RL=\Delta ^{2k}_R\times\Delta _L^{2k}$.
}
\vskip 11pt

\bpr
The cohomology $H^{2k}(FG^{(2i[2k])}(F_{\o v}\times F_v),\Delta ^{2k}\RL)$ of the bisemisheaf
$FG^{(2i[2k])}(F_{\o v}\times F_v)$ shifted under the action of the bioperator
$D ^{2k}_R\otimes D _L^{2k}$ must be expressed by means of the homology
$H_{2k}(FG^{(2i[2k])}(F_{\o v}\times F_v),\Fs^{2k}\RL(\TAN))$ with value in the bilinear fibre
$\Fs^{2k}\RL(\TAN)$ as developed in  chapter~2 of \cite{Pie6}.

As \[
H^{2k}(FG^{(2i[2k])}(F_{\o v}\times F_v),FG^{(2k)}(F_{\o v}\times F_v))=
\FREPSP(\GL_{2k}(F_{\o v}\times F_v))\]
 and as
 \begin{align*}
 H_{2k}(FG^{(2i[2k])} ( F_{\o v}\times F_v),\Fs^{2k}\RL (\TAN ))
 &= \FREPSP (\GL_{2k}( \rit\times  \rit ))\\
 &\simeq  \Ad\FREPSP (\GL_{2k}( \rit\times  \rit ))
\;, \end{align*}
we get the thesis.\epr
\vskip 11pt

\subsection{Proposition}

{\em
The bilinear cohomolgy of the shifted bisemisheaf (also called bilinear mixed cohomology)
$FG^{(2n[2k])}(F_{\o v}\times F_v)$  is given by the functional representation space of the bilinear general semigroup
$\GL_{2i}(F_{\o v}\times F_v)$ shifted in $2k$ real geometric dimensions according to:
\[
H^{2i-2k}(FG^{(2n[2k])}(F_{\o v}\times F_v),FG^{(2i[2k])}(F_{\o v}\times F_v))
=\FREPSP(\GL_{2i[2k]}(F_{\o v}\otimes \rit)\times (F_{v}\otimes \rit ))\]
where
$\FREPSP(\GL_{2i[2k]}(F_{\o v}\otimes \rit)\times (F_{v}\otimes \rit))$ is a condensed notation for\linebreak
$\FREPSP(\GL_{2k}(\rit\times \rit)) \times \FREPSP(\GL_{2i}(F_{\o v}\times F_v)) $.
}
\vskip 11pt

\bpr
Referring to section~4.1 giving the decomposition of the shifted bisemisheaf\linebreak
$FG^{(2i[2k])}(F_{\o v}\times F_v)$ into:
\[
FG^{(2i[2k])}(F_{\o v}\times F_v)= (\Delta ^{2k}_R\times\Delta _L^{2k})
\oplus FG^{(2i-2k)}(F_{\o v}\times F_v)\;,\]
we see that the cohomology $H^{2i-2k}(FG^{(2n[2k])}(F_{\o v}\times F_v),FG^{(2i[2k])}(F_{\o v}\times F_v))$ must similarly decompose into:
\begin{multline*}
H^{2i-2k}(FG^{(2n[2k])}(F_{\o v}\times F_v),FG^{(2i[2k])}(F_{\o v}\times F_v))\\
=H_{2k}(FG^{(2n[2k])} (F_{\o v}\times F_v),\Fs^{2k}\RL (\TAN ))\hspace*{3cm}\\
\times \L[
H^{2k}(FG^{(2n[2k])}(F_{\o v}\times F_v),FG^{(2k)}(F_{\o v}\times F_v))\R]\\
\oplus \L[H^{2i-2k}(FG^{(2n[2k])}(F_{\o v}\times F_v),FG^{(2i-2k)}(F_{\o v}\times F_v))\R]\;.
\end{multline*}
Taking into account that:
\begin{multline*}
H_{2k}(FG^{(2i[2k])}(F_{\o v}\times F_v),\Fs^{2k}\RL (\TAN ))
\times H^{2k}(FG^{(2n[2k])}(F_{\o v}\times F_v),FG^{(2k)}(F_{\o v}\times F_v))\\
=\FREPSP (\GL_{2k}(F_{\o v}\times \rit )\times (F_{v}\times \rit ))
\end{multline*}
and that
\[
H^{2i-2k}(FG^{(2n[2k])}(F_{\o v}\times F_v),FG^{(2i-2k)}(F_{\o v}\times F_v))
=\FREPSP (\GL_{2i-2k}(F_{\o v}\times F_v))\;, \]
we get the thesis.\epr
\vskip 11pt

\subsection{Cohomotopy resulting from the action of a differential bioperator}

In order to introduce a mixed homotopy bisemigroup in relation to a mixed Hurewicz homomorphism to be defined, we have to precise what must be the cohomotopy bisemigroup corresponding to the homology
$H_{2k}(FG^{(2i[2k])}(F_{\o v}\times F_v),\Fs^{2k}\RL(\TAN))$ with coefficients in the bifibre
$\Fs^{2k}\RL(\TAN)$.

As $H_{2k}(FG^{(2i[2k])}(F_{\o v}\times F_v),\Fs^{2k}\RL(\TAN))=
\fREPSP(\GL_{2k}(\rit\times \rit))$ and as a cohomotopy bisemigroup refers, according to section~2.16, to classes resulting from inverse deformations of Galois representations, under the circumstances of the shifted general bilinear semigroup
$\GL_{2i[2k]}(\wt F_{\o v}\times \wt F_v)$, {\bbf the searched cohomotopy bisemigroup
$\Pi ^{2k}(FG^{(2i[2k])}(F_{\o v}\times F_v))$ must be described by classes\/}:
\Bean
\item resulting from inverse deformations
\[
(GD^{(2k)})^{-1}: \qquad FG^{(2k)} (\rit\times\rit )\quad \To \quad F G^{(2k)} (\rit\times\rit )\]
of the differential Galois representation \cite{Car} of
$\GL_{2k}(\rit\times\rit)$;

\item depending on the classes of deformations of the Galois representation of
$\GL_{2k}(\wt F_{\o v}\times \wt F_v)$, i.e. the homotopy semigroup
$\Pi _{2k}(FG^{(2i[2k])}(F_{\o v}\times F_v))$.
\Ee

Consequently, {\bbf the mixed homotopy bisemigroup $\Pi _{2i[2k]}(FG^{(2n[2k])}(F_{\o v}\times F_v))$ of the shifted bisemisheaf
$FG^{(2n[2k])}(F_{\o v}\times F_v)$,
 will be defined by the product\/}:
\[
\Pi _{2i[2k]}(FG^{(2n[2k])}(F_{\o v}\times F_v))
= \Pi ^{2k}(FG^{(2i[2k])}(F_{\o v}\times F_v))
\times \Pi _{2i}(FG^{(2n[2k])}(F_{\o v}\times F_v))\]
of the cohomotopy
$\Pi ^{2k}(FG^{(2i[2k])}(F_{\o v}\times F_v))$ resulting from the action of a differential bioperator
$(D^{2k}_R\otimes D^{2k}_L)$ on the bisemisheaf
$FG^{(2n)}(F_{\o v}\times F_v)$ by the homotopy\linebreak
$\Pi _{2i}(FG^{(2n[2k])}(F_{\o v}\times F_v))$ of the bisemisheaf
$FG^{(2n)}(F_{\o v}\times F_v)$ shifted in $2k$ geometric dimensions.
\vskip 11pt

\subsection{Proposition}

{\em
The mixed bilinear semigroup homomorphism of Hurewicz will be given by:
\[
mhH: \qquad 
\Pi _{2i[2k]}(FG^{(2n[2k])}(F_{\o v}\times F_v))\quad \To \quad
H^{2i-2k}(FG^{(2n[2k])}(F_{\o v}\times F_v),\zit\times_{(D)}\zit )\;, \]
i.e. {\bbf a restricted $\Pi $-homology-$\Pi $-cohomology\/}.}
\vskip 11pt

\bpr
Indeed, the classes of the entire mixed bilinear cohomology
$H^{2i-2k}(FG^{(2n[2k])}(F_{\o v}\times F_v),\zit\times_{(D)}\zit )$ are the classes of the mixed homotopy bisemigroup
$\Pi _{2i[2k]}(FG^{(2n[2k])}(F_{\o v}\times F_v))$ and correspond to the bisemilattice
$\Lambda ^{(2i)}_{\o v}\otimes \Lambda ^{(2i)}_v \subset G^{(2i)}(F_{\o v}\times F_v)$ deformed by the mixed deformations of the considered Galois representation associated with\linebreak 
$\Pi _{2i[2k]}(FG^{(2n[2k])}(F_{\o v}\times F_v))$.\epr
\vskip 11pt

\subsection{Proposition (Chern mixed restricted character)}

{\em Let
\[
C_{[k]}(FG^{(2i[2k])}(F_{\o v}\times F_v)): \;
K _{2k}(FG^{(2i[2k])}(F_{\o v}\times F_v)) \To 
H _{2k}(FG^{(2i[2k])}(F_{\o v}\times F_v),\Fs^{2k}\RL (\TAN ))
\]
denote the restricted Chern character in the  operator bilinear $K$-homology in such a way that $C_{[k]}$ corresponds to the Todd class $J(FG^{(2i[2k])}(F_{\o v}\times F_v))\equiv J (\TAN   (FG^{(2k)} (F_{\o v}\times F_v )))$ according to chapter~3 of \cite{Pie6}

and let
\[
C^i(FG^{(2n)}(F_{\o v}\times F_v)): \quad
K^{2i}(FG^{(2n)}(F_{\o v}\times F_v)) \To 
H^{2i}(FG^{(2n)}(F_{\o v}\times F_v),FG^{(2i)}(F_{\o v}\times F_v))\]
be the Chern restricted character in the bilinear $K$-cohomology.

Then, {\bbf the Chern mixed restricted character in the $K$-homology-$K$-cohomology corresponds to a bilinear version of the index theory \cite{A-S} and is given by\/}:
\[
C_{[k]}\cdot C^i(FG^{(2n[2k])}(F_{\o v}\times F_v)): \qquad
K^{2i-2k}(FG^{(2n)}(F_{\o v}\times F_v))\quad \To \quad
H^{2i-2k}(FG^{(2n)}(F_{\o v}\times F_v))\]
where the mixed topological (bilinear) $K$-theory 
$K^{2i-2k}(FG^{(2n)}(F_{\o v}\times F_v))$ is given by:
\[
K^{2i-2k}(FG^{(2n)}(F_{\o v}\times F_v))
=K _{2k}(FG^{(2i[2k])}(F_{\o v}\times F_v))
\times K ^{2i}(FG^{(2n)}(F_{\o v}\times F_v))\]
where $K _{2k}(FG^{(2i[2k])}(F_{\o v}\times F_v))$ is the topological bilinear {\bbf contracting\/} $K$-theory of contracting tangent-bibundles (with contracting bifibres).
}
\vskip 11pt

\paragraph{Sketch of proof}: \hspace{5mm} 
The differential bioperator
$(D^{2k}_R\otimes D^{2k}_L)$ defined by its biaction
\[D^{2k}_R\otimes D^{2k}_L: \qquad
FG^{(2i)}(F_{\o v}\times F_v)\quad \To \quad
FG^{(2i[2k])}(F_{\o v}\times F_v)\]
leads to the Chern mixed restricted character
\[
C_{[k]}\cdot C^i(FG^{(2n[2k])}(F_{\o v}\times F_v))
=J(FG^{(2i[2k])}(F_{\o v}\times F_v))\times
 C^i(FG^{(2n)}(F_{\o v}\times F_v))\]
 corresponding to a bilinear version of the index theorem according to chapter~3 of \cite{Pie6},
 
 in such a way that the restricted Chern character
 $C_{[k]}(FG^{(2i[2k])}(F_{\o v}\times F_v))$ in the operator bilinear $K$-homology corresponds to the Todd class
 \begin{align}
J(FG^{(2i[2k])}(F_{\o v}\times F_v))&=
J ( \TAN (FG^{(2k)} (F_{\o v}\times F_v )))\notag\\
&= C_{[k]}(D^{2k}_R\otimes D^{2k}_L)\;.\tag*{\eop}
 \end{align}
 \vskip 11pt

 
\subsection{Proposition}

{\em If the classes of mixed bilinear cohomology 
$H^{2i-2k} (FG^{(2n)}(F_{\o v}\times F_v),\zit\times_{(D)} \zit )$ are classes deformed by the mixed deformations of the Galois representation associated with\linebreak  
$\Pi _{2i[2k]} (FG^{(2n)}(F_{\o v}\times F_v))$,
then {\bbf we can define a mixed lower bilinear (algebraic) $K$-theory by the equality\/} 
\resp{homomorphism}:
\[
K^{2i-2k}(FG^{(2n)}(F_{\o v}\times F_v))\underset{(\to)}{=}
\Pi _{2i[2k]}(FG^{(2n)}(F_{\o v}\times F_v))\]
implying that the mixed topological bilinear $K$-theory 
$K^{2i-2k}(FG^{(2n)}(F_{\o v}\times F_v))$ is \resp{corresponds to} a bisemigroup of deformed vector bibundles.
}
\vskip 11pt

\bpr
The thesis results from the commutative diagram:
\vspace*{20pt}
\[
\begin{psmatrix}[colsep=1cm,rowsep=.6cm]
K^{2i-2k}(FG^{(2n)}(F_{\o v}\times F_v)) && \Pi _{2i[2k]}(FG^{(2n)} (F_{\o v}\times F_v))\\[15pt]
& H^{2i-2k}(FG^{(2n)}(F_{\o v}\times F_v) )
\psset{arrows=->,nodesep=5pt}
\everypsbox{\scriptstyle}
\ncline{1,1}{1,3}^{\Bsm \text{mixed lower bilinear}\\\text{$K$-theory}\\ -\Esm}_-
\ncline{1,1}{2,2}<{\Bsm \text{inverse Chern mixed}\\ \text{restricted character}\Esm}
\ncline{2,2}{1,3}>{\Bsm \text{mixed homomorphism}\\ \text{of Hurewicz}\Esm}
\end{psmatrix}
\]
implying that
$K^{2i-2k}(FG^{(2n)}(F_{\o v}\times F_v))$ is a bisemigroup of vector bibundles deformed by mixed deformations of the Galois representations associated with
$\Pi _{2i[2k]}(FG^{(2n)}(F_{\o v}\times F_v))$.\epr
\vskip 11pt

\subsection{Proposition}

{\em
{\bbf The total Chern mixed character in the $K$-homology-$K$-cohomology\/}
\[
ch^*_*(FG^{(2n[2k])}(F_{\o v}\times F_v)): \quad
K_{2*}K^{2*}(FG^{(2n)}(F_{\o v}\times F_v))   \To  
H_{2*}H^{2*}(FG^{(2n)}(F_{\o v}\times F_v)) \]
where:
\Bean
\item $ch^*_*(FG^{(2n[2k])}(F_{\o v}\times F_v))=\sum\limits_k\sum\limits_i C_{[k]}
C^i(FG^{(2n[2k])}(F_{\o v}\times F_v))$,

\item $K_{2*}K^{2*}(FG^{(2n)}(F_{\o v}\times F_v))=\sum\limits_k\sum\limits_i
K^{2i-2k}(FG^{(2n)}(F_{\o v}\times F_v))$,

\item $H_{2*}H^{2*}(FG^{(2n)}(F_{\o v}\times F_v))=\sum\limits_k\sum\limits_i
H^{2i-2k}(FG^{(2n)}(F_{\o v}\times F_v))$,
\Ee
as well as the total mixed bilinear semigroup homomorphism of Hurewicz:
\[
mhH^*: \qquad
\Pi _{2*[2*]}(FG^{(2n)}(F_{\o v}\times F_v))\quad \To \quad
H_{2*}H^{2*}(FG^{(2n)}(F_{\o v}\times F_v))\]
where $\Pi _{2*[2*]}(FG^{(2n)}(F_{\o v}\times F_v))=
\sum\limits_k\sum\limits_i\Pi _{2i[2k]}(FG^{(2n)}(F_{\o v}\times F_v))$ implies the total mixed lower bilinear (algebraic) $K_*K^*$-theory given by the equality:
\[
K_{2*}K^{2*}(FG^{(2n)}(F_{\o v}\times F_v))=
\Pi _{2*[2*]}(FG^{(2n)}(F_{\o v}\times F_v))\;.\]
}
\vskip 11pt

\bpr This results evidently from the preceding sections.\epr
\vskip 11pt

\subsection{Proposition}

{\em
Let
\[ f : \qquad
FG^{(2n)}_Y(F_{\o v}\times F_v)\quad \To \quad
FG^{(2n)}_X(F_{\o v}\times F_v)\]
be a morphism between two compact bisemisheaves and let
\[ f_{!!} \qquad
K_{2*}K^{2*}(FG^{(2n)}_Y(F_{\o v}\times F_v))\quad \To \quad
K_{2*}K^{2*}(FG^{(2n)}_X(F_{\o v}\times F_v))\]
be the homomorphism between the corresponding total mixed
$K_{2*}K^{2*}$-theories.

Then, {\bbf the bilinear version of the Riemann-Roch theorem\/} {\em \cite{B-S}, \cite{A-H}, \cite{Gil1}} asserts that the diagram:
\[
\begin{psmatrix}[colsep=1cm,rowsep=.6cm]
K_{2*}K^{2*}(FG^{(2n)}_Y(F_{\o v}\times F_v)) & &
 K_{2*}K^{2*}(FG^{(2n)}_X(F_{\o v}\times F_v))\\[15pt]
H_{2*}H^{2*}(FG^{(2n)}_Y(F_{\o v}\times F_v)) & &
 H_{2*}H^{2*}(FG^{(2n)}_X(F_{\o v}\times F_v))
\psset{arrows=->,nodesep=5pt}
\everypsbox{\scriptstyle}
\ncline{1,1}{1,3}^{f_{!!}}
\ncline{2,1}{2,3}^{f^*_*}
\ncline{1,1}{2,1}>{ch^*_{*Y}}
\ncline{1,3}{2,3}>{ch^*_{*X}}
\end{psmatrix}
\]
is commutative,

or that:
\[ f^*_*\circ ch^*_{*Y}= ch^*_{*X}\circ f_{!!}\;.\]
}
\vskip 11pt

\bpr
\Bi
\item The linear classical version \cite{B-S}, \cite{Gil1}, \cite{A-H} of the Riemann-Roch theorem is
\[ 
f_*(ch(y)\cdot T(Y))=ch(f_!(y))\cdot T(X)\]
for any proper morphism $f:Y\to X$ between nonsingular, irreducible quasiprojective varieties where:
\Bi
\item $f_!: K(Y)\to K(X)$, $y\in K(Y)$,
\item $f_*: H^*(Y,\qit)\to H^*(X,\qit)$, 
\item $T(Y)$ is the Todd class of the tangent bundle to $Y$.
\Ei

\item The mixed bilinear version of the Riemann-Roch theorem
\[ f^*_*\circ ch^*_{*Y}= ch^*_{*X}\circ f_{!!}\]
then corresponds to the linear classical version if the total Todd class
$J(FG^{(2*[2k])}(F_{\o v}\times F_v))=
\sum\limits_k\sum\limits_ic_{[k]}(FG^{(2i[2k])}(F_{\o v}\times F_v))$ is the total Chern character $ch_*$ in the operator bilinear $K$-homology:
\[
ch_*: \qquad
K _{2*}(FG^{(2n)}(F_{\o v}\times F_v))\quad \To \quad
H_{2*}(FG^{(2n)}(F_{\o v}\times F_v))\;.\]

Remark that this way of envisaging the Riemann-Roch theorem by operator $K$-homology-$K$-cohomology is much more natural than the classical one working only in the form of $K$-cohomology.\epr
\Ei
\vskip 11pt

\subsection{Operator on the functional representation space of the infinite general bilinear semigroup}

In order to develop a bilinear version of the algebraic mixed $KK$-theory relative to the dynamical global program of Langlands, we have to introduce a higher bilinear algebraic operator $K$-theory relative to cohomotopy.

In this respect, we have to introduce the classical (and quantum) infinite bilinear semigroup acting by the biactions of the differential bioperator
$\{(D^{2k)}_R\otimes D^{2k}_L)\}_k$ on the infinite bilinear semisheaf
$\FGL(F_{\o v}\times F_v)$
over  the infinite bilinear semigroup
$\GL(F_{\o v}\times F_v)=
\underset{i}{\underrightarrow\lim}\GL_{2i}(F_{\o v}\times F_v)$
corresponding to the reducible functional representation space
$\fREPSP (\GL_{2n=2+\dots+2i+\dots+2n_s}(F_{\o v}\times F_v))$ of the reducible bilinear semigroup
$\GL_{2n}(F_{\o v}\times F_v)$ with $n\to\infty $.

Referring to the preceding sections, it appears {\bbf that the searched operator infinite bilinear semisheaf must be $\FGL(\rit\times\rit )=
\underset{k}{\underrightarrow\lim}\FGL_{2k}(\rit\times \rit )$ acting on
$\FGL(F_{\o v}\times F_v)$ by the biaction\/}
\[
\FGL_{\rit\times\rit}: \quad
\FGL (F_{\o v}\times F_v) \To
\FGL (\rit\times \rit )\times\FGL(F_{\o v}\times F_v)=
\FGL((F_{\o v}\otimes \rit )\times (F_{v}\otimes \rit ))\]
in such a way that:
\Bean
\item each factor of $\FGL(\rit\times \rit)$ acts on a factor of $\FGL(F_{\o v}\times F_v)$;

\item the factor ``$2k$'' $\FGL_{2k}(\rit\times \rit)\subset \FGL(\rit\times \rit)$ acts on the factor ``$2i$'' $\FGL_{2i}(F_{\o v}\times F_v)\subset \FGL(F_{\o v}\times F_v)$ in such a way that the (adjoint) functional representation space of
$\GL_{2i[2k]} ((F_{\o v}\otimes \rit )\times (F_{v}\otimes \rit ))$ is given by the shifted bisemisheaf
$FG^{(2i[2k])}(F_{\o v}\times F_v)=
\FREPSP (\GL_{2i[2k]}(F_{\o v}\otimes \rit )\times (F_{v}\otimes \rit ))$ with bilinear fibre
$\FREPSP (\GL_{2k}(\rit \times  \rit ))$ in the sense of proposition~4.3.
\Ee
\vskip 11pt

\subsection{The classifying bisemispace $\BFGL(\rit\times\rit )$}

The classifying bisemisheaf $\BFGL(\rit\times\rit )$ of $\FGL(\rit\times\rit )$
is the quotient of a weakly contractible bisemisheaf
$\EFGL(\rit\times\rit )$ by a free action of $\FGL(\rit\times\rit )$ in such a way that the continuous mapping:
\[
\GD_\rit : \qquad
\EFGL(\rit\times\rit )\quad \To \quad \BFGL(\rit\times\rit )\]
of the principal $\FGL(\rit\times\rit )$-bibundle over
$\BFGL(\rit\times\rit )$ is {\bbf a cohomotopy map corresponding to inverse Galois deformations of the Galois differential representation of $\BFGL(\rit\times\rit )$\/}.

The classifying bisemisheaf $\BFGL(\rit\times\rit )$ is then the base bisemisheaf of all equivalence classes of inverse deformations of the Galois differential representation of
$\FGL(\rit\times\rit )$ in the one-to-one correspondence with the kernels
$\FGL(\delta F_{\o{v+\ell }}\times \delta F_{v+\ell })$ of the maps
\[ \GD_\ell : \qquad 
\FGL( F_{\o{v+\ell }}\times  F_{v+\ell })\quad \To \quad
\FGL(F_{\o v}\times F_v)\;, \qquad 1\le\ell \le \infty \;, \]
introduced in proposition~3.8 and corollary~3.9.

The mixed classifying bisemisheaf
$\BFGL((F_{\o v}\otimes \rit)\times (F_v\otimes\rit ))$ then results from the (bi)action of $\BFGL(\rit\times\rit )$ on $\BFGL(F_{\o v}\times F_v )$.
\vskip 11pt


\subsection{The plus construction of $\BFGL(\rit\times\rit )$}

The ``plus'' construction of the mixed bilinear case of the Langlands dynamical global program is based on the map:
\[ \BFG_\rit(1) : \qquad
\BFGL(\rit\times\rit )\quad \To \quad
\BFGL(\rit\times\rit )^+\;, \]
unique up to cohomotopy, in such a way that:
\Bena
\item the kernel of the fundamental cohomotopy bisemigroup
$\Pi ^1(\BFG_\rit(1 ))$ be one-\linebreak dimensional inverse deformations of the Galois differential representation of\linebreak $\FGL(\rit\times\rit )$;

\item the cohomotopy fibre of $\BFG_\rit(1 )$ has the same integral homology as a (bi)point.
\Ee

The classifying bisemisheaf  $\BFGL(\rit\times\rit )^+$ is the base bisemisheaf of all equivalence classes of one-dimensional inverse deformations of the Galois differential representation of
$\FGL(\rit\times\rit )$ in one-to-one correspondence with the one-dimensional deformations of the Galois representation of
$\GL(\wt F_{\o v}\times \wt F_v )$ given by the kernel 
$\{\GL^{(1)}(\delta F_{\o {v+\ell }}\times\delta F_{v+\ell })\}_\ell $ of the maps
$\GD(1)_\ell $ according to proposition~3.11.
\vskip 11pt

\subsection{Proposition}

{\em
{\bbf The bilinear version of the mixed higher algebraic $KK$-theory of Quillen adapted to the Langlands dynamical bilinear global program\/} is:
\[
K_{2k}(FG^{(2n)}_{\rm red}(\rit\times \rit ))\times
K^{2i}(FG^{(2n)}_{\rm red}(F_{\o v}\times F_v))
=\Pi ^{2k}(\BFGL (\rit\times\rit )^+)\times\Pi _{2i}(\BFGL (F_{\o v}\times F_v)^+)\]
written in condensed form according to:
\[
K^{2i-2k}(FG^{(2n)}_{\rm red}((F_{\o v}\times\rit)\times (F_{ v}\times \rit)))
= \Pi _{2i[2k]} ( \BFGL ((F_{\o v}\otimes\rit )\times (F_v\otimes\rit ))^+ )\]
in such a way that the bilinear contracting $K$-theory $K_{2k}(FG^{(2n)}_{\rm red}(\rit\times \rit ))$, responsible for  a differential biaction, acts on the $K$-theory
$K^{2i}(FG^{(2n)}_{\rm red}(F_{\o v}\times F_v))$ of the reducible functional representation space
$FG^{(2n)}_{\rm red}(F_{\o v}\times F_v)$ of the bilinear semigroup
$\GL_{2n}(F_{\o v}\times F_v)$ in one-to-one correspondence with the biaction of the cohomotopy bisemigroup
$\Pi ^{2k}(\BFGL(\rit\times\rit )^+)$ of the ``plus'' classifying bisemisheaf
$\BFGL(\rit\times\rit )^+$.
}
\vskip 11pt

\bpr This mixed higher bilinear (algebraic) $KK$-theory is directly related to the commutative diagram:
\vspace*{20pt}
\[
\begin{psmatrix}[colsep=1cm,rowsep=.6cm]
K^{2i-2k} (FG^{(2n)}_{\rm red} (F_{\o v}\otimes\rit )\times (F_v\otimes\rit )) && 
\Pi _{2i[2k]}(\BFGL ((F_{\o v}\otimes\rit )\times (F_v\otimes\rit )^+))\\[15pt]
& \raisebox{-5mm}{$\hspace*{-3cm}H^{2i-2k}(FG^{(2n)}_{\rm red}(F_{\o v}\otimes\rit )\times (F_v\otimes\rit ) )\hspace*{-3cm}$}
\psset{arrows=->,nodesep=5pt}
\everypsbox{\scriptstyle}
\ncline{1,1}{1,3}^{\Bsm \text{mixed higher}\\ -\Esm}_{\Bsm - \\ \text{bilinear $K$-theory}\Esm}
\ncline{1,1}{2,2}<{\Bsm \text{inverse Chern mixed}\\ \text{higher restricted character}\Esm}
\ncline{2,2}{1,3}>{\Bsm \text{inverse restricted higher}\\ \text{$\Pi $-homology-$\Pi $-cohomology}\Esm}
\end{psmatrix}
\]
referring to the preceding sections.\epr
\vskip 11pt

\subsection{Proposition}

{\em
{\bbf The bilinear version of the total mixed higher algebraic $KK$-theory of Quillen adapted to the Langlands dynamical reducible global program of Langlands\/} is:
\[
K_* (FG^{(2n)}_{\rm red}(\rit\times \rit ))\times K^* (FG^{(2n)}_{\rm red}(F_{\o v}\times F_v))=
\Pi ^*(\BFGL(\rit\times \rit)^+)\times \Pi _*(\BFGL (F_{\o v}\times F_v)^+)\;.\]
}
\vskip 11pt

\bpr
This total mixed higher algebraic $KK$-theory has to be related to the commutative diagram:
\[
\begin{psmatrix}[colsep=1cm,rowsep=.6cm]
K^{2*-2*}(FG^{(2n)}_{\rm red}(F_{\o v}\otimes\rit )\times (F_v\otimes\rit )) && 
\Pi _{2x[2x]}(\BFGL ((F_{\o v}\otimes\rit )\times (F_v\otimes\rit ))^+) \\[15pt]
& \raisebox{-5mm}{$\hspace*{-3cm}H^{2*-2*}(FG^{(2n)}_{\rm red}(F_{\o v}\otimes\rit)\times (F_v\otimes\rit ) ))\hspace*{-3cm}$}
\psset{arrows=->,nodesep=5pt}
\everypsbox{\scriptstyle}
\ncline{1,1}{1,3}^{\Bsm \text{total mixed higher}\\ -\Esm}_{\Bsm - \\ \text{bilinear $KK$-theory}\Esm}
\ncline{1,1}{2,2}<{\Bsm \text{inverse Chern total}\\ \text{mixed higher character}\Esm}
\ncline{2,2}{1,3}>{\Bsm \text{inverse  higher}\\ \text{$\Pi $-homology-$\Pi $-cohomology}\Esm}
\end{psmatrix}
\]
\mbox{}\epr

\vfill

C. Pierre\\
Universit\'e de Louvain\\
Chemin du Cyclotron, 2\\
B-1348 Louvain-la-Neuve,  Belgium\\
pierre.math.be@gmail.com
}
\end{document}